\documentclass[11pt]{amsart}   

\usepackage{amssymb,amsthm,amsmath}
\usepackage{color}
\definecolor{darkblue}{rgb}{0, 0, .4}
\definecolor{grey}{rgb}{.7, .7, .7}
\usepackage[breaklinks]{hyperref}
\hypersetup{
	colorlinks=true,
	linkcolor=darkblue,
	anchorcolor=darkblue,
	citecolor=darkblue,
	pagecolor=darkblue,
	urlcolor=darkblue,
	pdftitle={},
	pdfauthor={}
}
	
\newtheorem{theorem}{Theorem}[section]
\newtheorem{lemma}[theorem]{Lemma}
\theoremstyle{definition}
\newtheorem{definition}[theorem]{Definition}
\newtheorem{example}[theorem]{Example}

\theoremstyle{remark}
\newtheorem{remark}[theorem]{Remark}
\numberwithin{equation}{section}

\theoremstyle{theorem}
\newtheorem{corollary}[theorem]{Corollary}

\newtheorem{proposition}[theorem]{Proposition}

\newtheorem{question}[theorem]{Question}

\newcommand{\N}[0]{\mathbb{N}}
\newcommand{\Z}[0]{\mathbb{Z}}

\newcommand{\s}[0]{\sigma}
\newcommand{\f}[0]{\varphi}
\newcommand{\ul}[1]{\underline{#1}}

\newcommand{\union}{\cup}

\newcommand{\sfc}{sub}
\newcommand{\afc}{add}

\newcommand{\gn}{ \bullet }  

\newcommand{\hs}{}  
\newcommand{\hd}{\diamond}  
\newcommand{\hz}{\circ}  
\newcommand{\hf}{\bullet}  
\newcommand{\hb}{{\color{grey} \bullet}}  
\newcommand{\hv}{\ast}  
\newcommand{\ha}{{\color{grey} \star}}  
\newcommand{\hh}{\star}  

\newcommand{\stringlessheap}{ \xymatrix @=-5pt @! } 
\newcommand{\heap}{ \xymatrix @=-9pt @! } 
\def\SDSize{6}  
\def\SDSizeTANRIGHT{4}  
\def\SDSizeTANLEFT{2}  
\def\SDMidpt{3}  
\def\SDColor{blue}
\def\SDEColor{black}

\newcommand{\StringLXD}[1]{{\color{\SDColor}\xy (0, \SDSize)*{}; (\SDSize, 0)*{}; **\crv{~**\dir{.}(\SDMidpt,\SDMidpt)};  (\SDMidpt, \SDMidpt)*{{\color{\SDEColor}#1}}; \endxy }}
\newcommand{\StringL}[1]{{\color{\SDColor}\xy (\SDSize, \SDSize)*{}; (0, 0)*{}; (0, \SDSize)*{}; **\crv{(\SDSizeTANLEFT,\SDMidpt)};  (\SDMidpt, \SDMidpt)*{{\color{\SDEColor}#1}}; \endxy}}
\newcommand{\StringR}[1]{{\color{\SDColor}\xy (0,0)*{}; (\SDSize, 0)*{}; (\SDSize, \SDSize)*{}; **\crv{(\SDSizeTANRIGHT,\SDMidpt)}; (\SDMidpt, \SDMidpt)*{{\color{\SDEColor}#1}}; \endxy }}
\newcommand{\StringLX}[1]{{\color{\SDColor}\xy (0, \SDSize)*{}; (\SDSize, 0)*{}; **\crv{(\SDMidpt,\SDMidpt)};  (\SDMidpt, \SDMidpt)*{{\color{\SDEColor}#1}}; \endxy }}
\newcommand{\StringRX}[1]{{\color{\SDColor}\xy (\SDSize, \SDSize)*{}; (0, 0)*{}; **\crv{(\SDMidpt,\SDMidpt)};  (\SDMidpt, \SDMidpt)*{{\color{\SDEColor}#1}}; \endxy }}
\newcommand{\StringLR}[1]{{\color{\SDColor}\xy (0, 0)*{}; (0, \SDSize)*{}; **\crv{(\SDSizeTANLEFT,\SDMidpt)};  (\SDSize, 0)*{}; (\SDSize, \SDSize)*{}; **\crv{( \SDSizeTANRIGHT,\SDMidpt)}; (\SDMidpt, \SDMidpt)*{{\color{\SDEColor}#1}}; \endxy }}
\newcommand{\StringLRX}[1]{{\color{\SDColor}\xy (0, \SDSize)*{}; (\SDSize, 0)*{}; **\crv{(\SDMidpt,\SDMidpt)};  (\SDSize, \SDSize)*{}; (0, 0)*{}; **\crv{(\SDMidpt,\SDMidpt)}; (\SDMidpt, \SDMidpt)*{{\color{\SDEColor}#1}}; \endxy }}

\usepackage[all, dvips, knot]{xy} 
\usepackage[all, knot, dvips, color]{xypic} 
\xyoption{arc}

\newcommand{\w}{\mathsf{w}}
\newcommand{\y}{\mathsf{y}}
\renewcommand{\u}{\mathsf{u}} 
\renewcommand{\v}{\mathsf{v}} 
\newcommand{\p}{\mathsf{p}}

\newcommand{\block}[1]{\begin{bmatrix} #1 \end{bmatrix}}

\newcommand{\bcone}{\operatorname{Cone_{\wedge}}}
\newcommand{\ucone}{\operatorname{Cone^{\vee}}}
\newcommand{\chs}[2]{ \left(\begin{smallmatrix} #1 \\#2 \end{smallmatrix} \right)}

\begin{document}

\title{Embedded factor patterns for Deodhar elements \linebreak in
Kazhdan--Lusztig theory}

\begin{abstract}
The Kazhdan--Lusztig polynomials for finite Weyl groups arise in the
geometry of Schubert varieties and representation theory.  It was
proved very soon after their introduction that they have nonnegative
integer coefficients, but no simple all positive interpretation
for them is known in general.  Deodhar \cite{d} has given a framework
for computing the Kazhdan--Lusztig polynomials which generally
involves recursion.  We define embedded factor pattern avoidance for
general Coxeter groups and use it to characterize when Deodhar's
algorithm yields a simple combinatorial formula for the
Kazhdan--Lusztig polynomials of finite Weyl groups.  Equivalently, if
$(W, S)$ is a Coxeter system for a finite Weyl group, we classify the
elements $w \in W$ for which the Kazhdan--Lusztig basis element $C'_w$
can be written as a monomial of $C'_s$ where $s \in S$.  This work
generalizes results of Billey--Warrington \cite{b-w} that identified
the Deodhar elements in type $A$ as 321-hexagon-avoiding permutations,
and Fan--Green \cite{f-g} that identified the fully-tight Coxeter
groups.
\end{abstract}

\author{Sara C. Billey}
\address{Department of Mathematics Box 354350, University of Washington, Seattle, WA 98195}
\email{\href{mailto:billey@math.washington.edu}{\texttt{billey@math.washington.edu}}}
\urladdr{\url{http://www.math.washington.edu/\~billey/}}

\author{Brant C. Jones}
\address{Department of Mathematics Box 354350, University of Washington, Seattle, WA 98195}
\email{\href{mailto:brant@math.washington.edu}{\texttt{brant@math.washington.edu}}}
\urladdr{\url{http://www.math.washington.edu/\~brant/}}

\thanks{NSF grant DMS-9983797 supported both authors and the second
author also received support from NSF VIGRE grant DMS-0354131.}


\date{\today}

\maketitle



\bigskip
\section{Introduction}

The Kazhdan--Lusztig polynomials for finite Weyl groups \cite{k-l}
arise as Poincar\'e polynomials for intersection cohomology of
Schubert varieties \cite{K-L2} and as a $q$-analogue of the
multiplicities for Verma modules \cite{BeilBern,BryKash}.  They are
defined to be the coefficients in the transition matrix for expanding
the Kazhdan--Lusztig basis elements in the Hecke algebra associated to
the Weyl group into the standard basis.  Several algorithms exist,
formulas for special cases, and interesting properties are known for
these polynomials; see for example
\cite{h,Deodhar94,MW02,polo,l-s,brenti04,b-b}.  In particular, these
polynomials have nonnegative integer coefficients but no simple all
positive formula for the coefficients is known in general for all Coxeter
groups.

Deodhar \cite{d} proposes a combinatorial framework for determining
the Kazhdan--Lusztig polynomials of an arbitrary Coxeter group. The
algorithm he describes is shown to work for all Coxeter groups where
the Kazhdan--Lusztig polynomials are known to have nonnegative integer
coefficients which includes Weyl groups and the Coxeter groups
associated to crystallographic Kac--Moody groups.  Under certain
conditions, Deodhar's algorithm for determining the Kazhdan--Lusztig
polynomials turns out to be a beautiful combinatorial formula.  These
conditions are also equivalent to the Kazhdan--Lusztig basis element
$C'_{w}$ being equal to a product of $C'_{s}$'s indexed by generators
of the Coxeter group.  We say that $w$ is \em Deodhar \em when it
satisfies these conditions.  In 1999, Billey and Warrington \cite{b-w}
gave an efficient characterization of the Deodhar elements in the
symmetric group as \textit{321-hexagon avoiding} permutations.  Their
results extend to finite \textit{linear} Weyl groups, types $A,B,F,G$.
Our goal is to give a similar characterization for all finite Weyl
groups.

In this paper we give two characterizations of the Deodhar elements
for all finite Weyl groups.  One characterization is given in terms of
1-line pattern avoidance in analogy with the type $A$ result.  This
characterization gives a polynomial time algorithm to test for the
Deodhar status of an element, but involves a long list of patterns.
The second characterization is in terms of a new type of pattern
called an \textit{embedded factor}.  These patterns are defined in
terms of reduced expressions, and generalize containment in the
2-sided weak Bruhat order.  Theorem~\ref{t:main} states that the
Deodhar elements of Weyl groups can be characterized by avoiding
embedded factors from the following list, as well as an additional
1-line pattern for type $D$.

\bigskip
\noindent
\begin{figure}[h]
\begin{tabular}{|p{0.45in}|p{2in}|p{2.7in}|} 
	\hline
	Type & Coxeter Graph & Embedded Factor Patterns \\
	\hline
	$I_2(m)$, $m \geq 3$ & \xymatrix @-1pc {
\gn_1 \ar@{-}[r]^m & \gn_2 \\
} & $s_1 s_2 s_1$, \  $s_2 s_1 s_2$ \ (short braids) \\
	\hline
	$A_7$ &  \xymatrix @-1pc {
\gn_1 \ar@{-}[r] & \gn_2 \ar@{-}[r] & \gn_3 \ar@{-}[r] & \gn_4 \ar@{-}[r] & \gn_5 \ar@{-}[r] & \gn_6 \ar@{-}[r] & \gn_7 \\
}
& $s_5 s_6 s_7 s_3 s_4 s_5 s_6 s_2 s_3 s_4 s_5 s_1 s_2 s_3$ ($HEX$) \\
	\hline
	$B_7 / C_7$ &  \xymatrix @-1pc {
\gn_0 & \ar@{-}[l]_{4} \gn_1 \ar@{-}[r] & \gn_2 \ar@{-}[r] & \gn_3 \ar@{-}[r] & \gn_4 \ar@{-}[r] & \gn_5 \ar@{-}[r] & \gn_6 \\
}
& $s_4 s_5 s_6 s_2 s_3 s_4 s_5 s_1 s_2 s_3 s_4 s_0 s_1 s_2$ ($BHEX$) \\
	\hline
	$D_6$ &  \xymatrix @-1pc {
{\gn}_{\tilde{1}} &        &                &                            &  \\
\gn_1 & \ar@{-}[l] \ar@{-}[ul] \gn_2 \ar@{-}[r] & \gn_3 \ar@{-}[r] & \gn_4 \ar@{-}[r] & 
 \gn_5\\
}
& $s_3 s_4 s_5 s_1 s_2 s_3 s_4 s_{\tilde{1}} s_2 s_3 s_1$ ($HEX_4$)  \\
	\hline
	$D_7$ &  \xymatrix @-1pc {
{\gn}_{\tilde{1}} &        &                &                            &  \\
\gn_1 & \ar@{-}[l] \ar@{-}[ul] \gn_2 \ar@{-}[r] & \gn_3 \ar@{-}[r] & \gn_4 \ar@{-}[r] & \gn_5 \ar@{-}[r] & \gn_6\\
}

& { $s_3 s_4 s_5 s_6 s_2 s_3 s_4 s_5 s_{\tilde{1}} s_2 s_3 s_4 s_1 s_2 s_3$ ($HEX_2$)

$s_4 s_5 s_6 s_2 s_3 s_4 s_5 s_1 s_2 s_3 s_4 s_{\tilde{1}} s_2 s_1$ ($HEX_{3a}$)

$s_1 s_4 s_5 s_6 s_2 s_3 s_4 s_5 s_{\tilde{1}} s_2 s_3 s_4 s_1 s_2$ ($HEX_{3b}$) } \\

	\hline
	$D_8$ &  \xymatrix @-1pc {
{\gn}_{\tilde{1}} &        &                &                            &  \\
\gn_1 & \ar@{-}[l] \ar@{-}[ul] \gn_2 \ar@{-}[r] & \gn_3 \ar@{-}[r] & \gn_4 \ar@{-}[r] & 
\gn_5 \ar@{-}[r] & \gn_6 \ar@{-}[r] & \gn_7 \\
}
& { $s_4 s_5 s_6 s_7 s_3 s_4 s_5 s_6 s_2 s_3 s_4 s_5 s_{\tilde{1}} s_1
s_2 s_3 s_4$ \ \ \ \ \ \ (diamond, to be avoided as a 1-line
pattern) } \\
	\hline
	$E_6$ & 
	\xymatrix @-1pc {
    & & \gn_{5} &        &                &                            &  \\
    \gn_0 \ar@{-}[r] & \gn_1 & \ar@{-}[l] \ar@{-}[u] \gn_2 \ar@{-}[r] & \gn_3 \ar@{-}[r] & \gn_4 & \\
    }
    & { $s_0 s_1 s_2 s_5 s_3 s_4 s_2 s_3 s_1 s_2 s_5 s_0 s_1$ 
	 $s_5 s_1 s_2 s_3 s_0 s_1 s_2 s_5 s_4 s_3 s_2 s_1 s_0$ 
	 $s_1 s_2 s_5 s_3 s_4 s_2 s_3 s_1 s_2 s_5 s_0 s_1 s_2$ 
	 $s_2 s_5 s_1 s_2 s_3 s_0 s_1 s_2 s_5 s_4 s_3 s_2 s_1$ } \\
	\hline
	$E_7$ & \xymatrix @-1pc {
    & & \gn_{5} &        &                &                            &  \\
    \gn_0 \ar@{-}[r] & \gn_1 & \ar@{-}[l] \ar@{-}[u] \gn_2 \ar@{-}[r] & \gn_3 \ar@{-}[r] & \gn_4 \ar@{-}[r] & \gn_5  \\
    }
    & { $s_0 s_1 s_2 s_3 s_4 s_6 s_5 s_2 s_3 s_4 s_1 s_2 s_3 s_0 s_1$
	$s_3 s_4 s_6 s_1 s_2 s_3 s_0 s_1 s_2 s_5 s_4 s_3 s_2 s_1 s_0$ 
	$s_1 s_2 s_3 s_4 s_6 s_5 s_2 s_3 s_4 s_1 s_2 s_3 s_0 s_1 s_2$
	$s_2 s_3 s_4 s_6 s_1 s_2 s_3 s_0 s_1 s_2 s_5 s_4 s_3 s_2 s_1$
	$s_5 s_2 s_3 s_4 s_6 s_1 s_2 s_5 s_3 s_4 s_2 s_3 s_0 s_1 s_2 s_5$ }  \\
	\hline
\end{tabular}
\caption{Minimal non-Deodhar patterns}
\label{f:patterns}
\end{figure}
\bigskip

The embedded factor patterns take into account different ways of
embedding one Weyl group into another as a parabolic subgroup.  For
example, the Weyl group of type $E_{8}$ has parabolic subgroups of
types $A_{2}, A_{7}, D_{6}, D_{7}, E_{6}$, and $E_{7}$ from this list.
Therefore, a Deodhar element in the Weyl group of type $E_8$
cannot have any embedded factors in the form of a short braid,
hexagon, $HEX_i$, or any of the $E_{6}$ or $E_{7}$ patterns.

In type $D_n$, the Deodhar elements must also avoid the ``diamond''
pattern $[\bar{1} 6 7 8 \bar{5} 2 3 4]$ as a 1-line
pattern.  This single 1-line pattern encapsulates an infinite
antichain of type $D$ embedded factor patterns further discussed in
Example~\ref{e:infantichain}.

We also provide a finite test in
Theorem~\ref{t:one.iff.embedded.finite} to determine when it is
possible to translate between classical 1-line pattern avoidance and
the embedded factor pattern avoidance of
Definition~\ref{d:embedded.factor}.  This result generalizes a fact
which is implicit in \cite{b-w} that avoiding the hexagon embedded
factor pattern can be characterized by avoiding 4 classical
permutation patterns when we restrict to the fully commutative
elements.  Theorem~\ref{t:one.iff.embedded.finite} also justifies
using the methods in this paper to study certain classical pattern
classes, and has been extended for type $A$ in \cite{j2}.

In Section~\ref{s:background}, we recall the basic definitions of
Kazhdan--Lusztig polynomials and basis elements.  We define the Deodhar
elements and recall the theorem that inspired this name.  In
Section~\ref{s:heaps}, we recall the heap of a reduced expression for
Weyl groups of types $A,B,D$, and in
Section~\ref{s:classical.patterns} we review the definition of
classical pattern avoidance.  The heaps will be the main tool for
proving the characterization in Theorem~\ref{t:main} presented in
Section~\ref{s:Deodhar.elements}.  In Section~\ref{s:d.heaps} we
reduce the proof of the main theorem to short braid avoiding elements.
In Sections~\ref{s:convex} and~\ref{s:classification.d}, we define the
convex elements and give a complete classification of convex Deodhar
elements.  Then in Section~\ref{s:d.proof}, we characterize the
non-convex Deodhar elements in type $D$ which completes the proof for
type $D$.  In Section~\ref{s:classification.e}, we complete the
classification of Deodhar elements for the remaining finite Weyl
groups.  Section~\ref{s:convex.patterns} gives a pattern comparison
result, which shows that the Deodhar property for type $D$ can be
characterized by avoiding finitely many 1-line patterns.  
Finally, we close with some open problems
and enumerative data in Section~\ref{s:open}.


\bigskip
\section{Background and notation}\label{s:background}

In this section we will set up our notation and review some of the
motivation for our main theorems.  For a reader unfamiliar with Coxeter
groups, we recommend either the classic text by Humphreys \cite{h} or
the recent text by Bj\"orner and Brenti \cite{b-b} for a more
combinatorial treatment.

Let $W$ be a Coxeter group with generating set $S$ and relations of
the form $(s_{i}s_{j})^{m(i,j)} = 1$.  The Coxeter graph for $W$ is
the graph on the generating set $S$ with edges connecting $s_{i}$ and
$s_{j}$ labeled $m(i,j)$ for all pairs $i,j$ with $m(i,j)>2$.  For
example, the table in Figure~\ref{f:patterns} shows the Coxeter
graphs for the finite Weyl groups that contain minimal non-Deodhar
patterns.  Note that if $m(i,j)=3$ it is customary to leave the
corresponding edge unlabeled.  

An \em expression \em is any product of generators from $S$.  The \em
length \em $l(w)$ of an element $w \in W$ is the minimum length of
any expression for the element $w$.  Such a minimum length expression
is called \em reduced\em.  Each element $w \in W$ can have several
different reduced expressions that represent it.  Given $w \in W$, we
represent a reduced expressions for $w$ in sans serif font, say
$\w=\w_{1} \w_{2}\cdots \w_{p}$ where each $\w_{i} \in S$.  

It is a theorem of Tits \cite{t} that every reduced expression for an
element $w$ of a Coxeter group can be obtained from any other by
applying a sequence of braid moves of the form 
\[
{\underbrace{s_i s_j s_i s_j \cdots }_{m(i,j)} } \mapsto
{\underbrace{s_j s_i s_j s_i \cdots}_{m(i,j)}}
\]
where $s_i$ and $s_j$ are generators in $S$ that appear in the reduced
expression for $w$, and each factor in the move has $m(i,j)$ letters.
Let the \em support \em of an element $w \in W$, denoted $supp(w)$, be
the set of all generators appearing in any reduced expression for $w$,
which is well-defined by Tits' theorem.  We say that the element $w$
is \em connected \em if $supp(w)$ is connected in the Coxeter graph
of $W$.

We define an equivalence relation on the set of reduced expressions
for a fixed Coxeter element where two reduced expressions are in the
same \em commutativity class \em if one can be obtained from the other
by \textit{commuting moves} of the form $s_i s_j \mapsto s_j s_i$,
where $m(i,j) = 2$.  In particular, if every reduced expression for
$w$ can be obtained from any other by commuting moves then we say $w$
is \em fully commutative\em.  By Tits' theorem, an element $w$ is
fully commutative if and only if no reduced expression for $w$
contains a consecutive subexpression of the form $s_i s_j s_i s_j
\cdots$ of length $m(i,j) \geq 3$.  

We call any expression of the form $s_i s_j s_i$ for $m(i,j) \geq 3$ a
\em short braid\em.  This name reflects the fact that we are not
considering any longer braid, even if $m(i,j) > 3$.  We caution the
reader that some authors have used the term short braid to refer to a
commuting move between two entries $s_i$ and $s_j$ where $m(i,j) = 2$.
The \textit{short braid avoiding} elements of a Coxeter group are
those with no reduced expression containing a factor $s_i s_j s_i$
where $s_i$ and $s_j$ are any pair of generators that do not commute
(i.e.  $m(i,j) \ne 2$).  Hence, a short braid avoiding element is also
fully commutative.

Given a Coxeter group $W$, we can form the \em Hecke algebra
$\mathcal{H}$ \em over $\Z[q^{1/2}, q^{-1/2}]$ with basis $\{ T_w
: w \in W \}$ and relations:
\begin{align}\label{e:hecke.def}
T_s T_w = & T_{s w} \text{ for } l(s w) > l(w)  \\
 (T_s)^2 = & (q - 1) T_{s} + q T_1 
\end{align}
where $T_1$ corresponds to the identity element.  In particular, this implies
that \[ T_w = T_{\w_1} T_{\w_2} \dots T_{\w_p} \] whenever $\w_1 \w_2
\dots \w_p$ is a reduced expression for $w$.

Kazhdan and Lusztig \cite{k-l} described another basis for
$\mathcal{H}$ that is invariant under the Hecke algebra involution
mapping
\begin{align*}
&  q \mapsto q^{-1}  \\
&  T_s \mapsto (T_s)^{-1}. 
\end{align*}
This basis, denoted $\{ C_w' : w \in W \}$, has important
applications in representation theory and algebraic geometry
\cite{k-l,K-L2}.  The Kazhdan--Lusztig polynomials $P_{x,w}(q)$ arise as the
``change of basis'' matrix between these two bases of $\mathcal{H}$:
\begin{equation*}
C_w' = q^{-\frac{1}{2} l(w)} \sum_{x \in W} P_{x,w}(q) \ T_x. 
\end{equation*}
The $C_{w}'$ are defined uniquely to be the Hecke algebra elements
that are invariant under the involution and have expansion
coefficients as above where $P_{x,w}$ is a  polynomial in $q$ with
\begin{equation}\label{e:max.degree}
\text{degree } P_{x,w}(q) \leq \frac{(l(w)-l(x)-1)}{2}
\end{equation}
and $P_{w,w}(q)=1$.  We use the notation $C'_{w}$ to be consistent
with the literature because there is already a related basis denoted
$C_{w}$.

For $w \in W$ and $s \in S$ with $l(sw)>l(w)$, the Kazhdan--Lusztig
basis elements multiply according to the rule
\begin{equation*}
 C_s' C_w' = C_{s w}' + \sum_{s z < z < w} \mu(z,w) C_z'
\end{equation*}
where $\mu(z, w)$ is the coefficient of $q^{ {1 \over 2} (l(w)-l(z)-1)
}$ (the term of highest possible degree) in the Kazhdan--Lusztig
polynomial $P_{z, w}(q)$.  This is appreciably more complicated than
the corresponding multiplication formula \eqref{e:hecke.def} in the
$\{ T_w \}$ basis.

Deodhar \cite{d} studied the case when $C'_{w}$ can be written simply
as a product of $C'_{s_i}$'s.  In this case, he also gives nice
combinatorial formulas for all the polynomials $P_{x,w}(q)$.  We will
describe Deodhar's defect statistic and his theorem in terms of masks
on reduced expressions.

Fix a reduced expression $\w = \w_{1} \cdots \w_{k}$.  Define a \em
mask \em $\s$ associated to the reduced expression $\w$ to be any
binary word $\s_1 \cdots \s_k$ of length $k = l(w)$.  Every mask
corresponds to a subexpression of $\w$ defined by $\w^\s =
\w_{1}^{\s_1} \cdots \w_{k}^{\s_k}$ where
\[
\w_{j}^{\s_j}  =
\begin{cases}
\w_{j}  &  \text{ if  }\s_j=1\\
\text{id}  &  \text{ if  }\s_j=0.
\end{cases}
\]
Each $\w^\s$ is a product of generators in a subsequence of $\w_{1}
\cdots \w_{k}$ so it determines an element of $W$ that is less than
$\w$ in \textit{Bruhat order}.  For $1\leq j\leq k$, we also consider initial
sequences of masks denoted $\s[j] = \s_1 \cdots \s_j$, and the
corresponding initial subexpressions $\w^{\s[j]} = \w_{1}^{\s_1}
\cdots \w_{j}^{\s_j}$.  For example, we have $\w^{\s[k]} = \w^\s$.  A
mask $\s$ is \em proper \em if it has at least one zero.

We say that a position $j$ (for $2 \leq j \leq k$) of $\w$ is a \em
defect \em with respect to the mask $\s$ if
\begin{equation}\label{e:deodhar.length}
l(\w^{\s[j-1]} \w_{j}) < l(\w^{\s[j-1]}).
\end{equation}
Note that a defect occurs in position $j$ if $\w^{\s[j-1]}$ satisfies the
length condition above; the value of $\s_j$ is irrelevant when
determining if $j$ is a defect.  Let $d_{\w}(\s)$ denote the number of
defects of $\w$ for a mask $\s$.  We will use the notation
$d(\s) = d_{\w}(\s)$ when the reduced word $\w$ is fixed.  

We are now in a position to define Deodhar's condition.

\begin{definition}\label{d:deodhar}
Let $w \in W$ be a Coxeter element with reduced expression $\w$, and
let $\s$ be a proper mask for $\w$.  We say that a position $j$ is a
\em zero-defect \em if $\s_j =0$ and $j$ is also a defect in $\w$.  We
say that position $j$ in $\w$ is a \em plain-zero \em if $\s_j =0$ and
$j$ is not a defect in $\w$.  Then, the mask $\s$ is \em Deodhar \em if
\begin{equation}\label{e:deodhar.ineq}
\text{\# of zero-defects of } \s < \text{\# of plain-zeros of } \s.
\end{equation}
Moreover, a reduced expression $\w$ is \em Deodhar \em if every
proper mask $\s$ on $\w$ is Deodhar.  
\end{definition}

\begin{example} Assume $s_{2} s_{1}s_{3}s_{2}$ is a reduced expression
in a Coxeter group.  The word/mask pair
\[ \w = \block{ s_2 & s_1 & s_3 & s_2 \\ 1 & 0 & 0 & 0 } \]
has a zero-defect in position 4 and plain-zeros in
positions 2 and 3.  One can verify that \eqref{e:deodhar.ineq} holds for
all proper masks on $\w$, so $\w$ is Deodhar.
\end{example}

For certain Coxeter groups, Deodhar has shown that when a reduced
expression $\w \in W$ satisfies this condition, the Kazhdan--Lusztig
polynomials $P_{x,w}$ can be obtained as the generating function that
counts masks on $\w$ with respect to the defect statistic.
Equivalently, $C'_{w}$ can be written as a product of $C'_{s_i}$'s.  He
also shows that the notion of being Deodhar is well-defined on Coxeter
group elements and is independent of the choice of reduced word used to
verify the condition.  Note that Deodhar actually used a slightly
different condition which is equivalent to the one given here in
\eqref{e:deodhar.ineq} \cite[Lemma 2]{b-w}.  The original condition is
\[ d(\s) \leq {1 \over 2}(l(\w) - l(\w^{\s}) -1) \] which comes
directly from the maximum degree bound of the Kazhdan--Lusztig
polynomial $P_{\w^{\s}, \w}(q)$ in \eqref{e:max.degree}.

\begin{theorem}\cite{d}\label{t:deodhar}
Let $W$ be any Coxeter group where the Kazhdan--Lusztig
polynomials are known to have nonnegative coefficients, and let
$\w=\w_{1}\cdots \w_{k}$ be a reduced expression for some $w \in W$.
Then the following are equivalent:
\begin{enumerate}
\item The expression $\w$ is Deodhar.

\item The element $w$ is Deodhar.

\item The Kazhdan--Lusztig basis element $C_w'$ is given by
	\[ C_w' = q^{-{1 \over 2} l(w)} \sum q^{d(\s)} T_{ \w^{\s} } \]
where the sum is over all masks $\s$ on $\w$.

\item For all $x \in W$, the Kazhdan--Lusztig polynomial $P_{x,w}$ is
given by
\[ P_{x,w}(q) = \sum q^{d(\s)} \]
where the sum is over all masks $\s$ on $\w$ such that $\w^\s = x$. 

\item The Kazhdan--Lusztig basis element $C_w'$ satisfies $C_w' =
C_{\w_1}'\cdots C_{\w_k}'$.\label{itfour} 

\item The Bott--Samelson resolution of the corresponding Schubert
variety $X_w$ is small.\label{itfive}

\item \label{itthree} The Poincar\'e polynomial for the full intersection 
cohomology group of the Schubert variety $X_w$ is 
$$
\sum_i \dim(\mathcal{IH}^{2i}(X_w))q^i = (1+q)^{l(w)}.
$$  
\end{enumerate}
\end{theorem}

\begin{remark}
The equivalence of (1) through (6) are implicit in Deodhar \cite{d}.  The equivalence
of (4) and (7) is proved explicitly by Billey and Warrington \cite{b-w}.
Lusztig \cite{lusztig} and Fan and Green \cite{f-g} have studied those elements
for which (5) holds.  These elements are called ``tight'' in the terminology of
those papers.
\end{remark}

The main goal of this paper is to give an efficient way to identify
Deodhar elements.  In the case when $W$ is the symmetric group, Billey
and Warrington \cite{b-w} gave a concise description of the Deodhar
elements as those that are 321-hexagon avoiding, where the term
``hexagon'' comes from the notion of a heap on a permutation.  We will
describe the heap construction and classical pattern avoidance in the
next two sections and return to the study of Deodhar elements in
Section~\ref{s:Deodhar.elements}.

\bigskip
\section{Heaps and string diagrams}\label{s:heaps}

Each reduced expression can be associated with a partial order called
the heap that we define below.  This partial order allows us to
visualize a reduced expression as a set of lattice points while
maintaining the pertinent information about the relations among the
generators.  Cartier and Foata \cite{cartier-foata} were among the
first to study heaps of dimers, and these were generalized to other
settings by Viennot \cite{viennot}.  More recently, Stembridge has
studied enumerative aspects of heaps \cite{s1,s2} in the context of
fully commutative elements.  We will use the heaps to visualize the
reduced expressions that appear in the table on Page 2 and prove our
characterization in type $D$.

\begin{definition}\label{d:heap.poset}
Suppose $\w = \w_1 \cdots \w_k$ is a fixed reduced expression, and
define a partial ordering on the indices $\{1, \cdots, k\}$ by the
transitive closure of the relation $i \lessdot j$ if $i < j$ and
$m(\w_i,\w_j) \neq 2$.  In particular, $i \lessdot j$ if $i < j$ and
$\w_i = \w_j$.  This partial order is called the \em heap \em of
$\w$.  We label the element $i$ of the poset by the corresponding
generator $\w_i$.
\end{definition}

\begin{remark}
Observe that heaps are well defined up to commutativity class, so if
$\u$ and $\v$ are two reduced expressions for $w$ in the same
commutativity class then the labeled heaps of $\u$ and $\v$ are
equal.  In particular, if $w$ is fully commutative then there is a
unique labeled heap poset for the element $w$, regardless of which
reduced expression is used to generate it.
\end{remark}

Let $G$ denote the Coxeter graph for $W$.  We can embed the heap poset
as a set of lattice points in $G \times \N$.  To do this, begin by
reading the reduced expression $\w$ from left to right, and drop a
point in the column representing each generator $\w_i$.  We can
envision each point as being ``fat'' and under the influence of
``gravity,'' in the sense that the point must fall to the lowest
possible position in the column over the generator corresponding to
$\w_i$ in the Coxeter graph without passing any previously placed
points in adjacent columns.  Here, we say two columns are adjacent
when they correspond to adjacent vertices in the Coxeter graph.  Since
generators that are adjacent in the Coxeter graph do not commute, we
must place the point representing $\w_i$ at a level that is \em above
\em the level of any other adjacent points that have already been
placed.  Because generators that are not adjacent in the Coxeter
graph do commute, points that lie in non-adjacent columns can slide
past each other or land at the same level.  

\begin{definition}\label{d:heap.lp}
Let $\w$ be a reduced expression for a Coxeter element.  We let
$Heap(\w)$ denote the lattice representation of the heap poset in $G
\times \N$ constructed as described in the preceding paragraph.
(We will amend this definition in Example~\ref{ex:type.d} to account
for the fork in the Coxeter graph of type $D$.)
\end{definition}

We will give several examples of heaps in different Coxeter groups in
the next three examples and also introduce some useful terminology for
permutations and signed permutations.
\begin{example}\label{ex:typeA}
 The Coxeter graph of type $A_{n-1}$ is the following:
\[
\xymatrix @-1pc {
\gn_1 \ar@{-}[r] & \gn_2 \ar@{-}[r] & \gn_3  & \ar@{-}[r] & \dotsb &  \gn_{n-1}\\
} .
\]
The corresponding Coxeter group is the symmetric group $S_{n}$.  We
may refer to elements in the symmetric group by the \textit{1-line
notation} $w=[w_{1}w_{2}w_{3}\dotsb w_{n}]$ where $w$ is the bijection
mapping $i$ to $w_{i}$ written in italic font.  The generators
$s_{1},s_{2},\dotsb, s_{n-1}$ are the adjacent transpositions where $s_{i}$
interchanges $i$ and $i+1$.  For example, when multiplying a
permutation on the right by $s_2$, we interchange the entries in
positions $2$ and $3$ of the 1-line notation, so for $w=[3412] \in S_{4}$ we
have $ws_{2}=[3142]$.  Dually, when multiplying on the left by $s_2$,
we interchange the digits $2$ and $3$ in the 1-line notation for the
element, so $s_{2}w=[2413]$.  One reduced expression for $w$ is $s_2
s_3 s_1 s_2$.  We build up the heap one generator at a time for the
reduced expression $\w = s_2 s_3 s_1 s_2$ in type $A_3$ as shown
below.
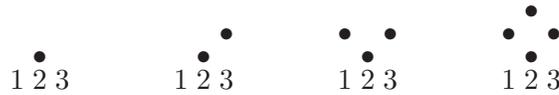
\begin{center}
\begin{tabular}{llll}
	\stringlessheap {
	\hs & \hs & \hs & \hs & \hs \\
	& \hs & \hs & \hs & \hs & \hs \\
	\hs & \hs & \hs & \hs & \hs \\
	& \hf & \hs & \hs & \hs & \hs \\
	1 & 2 & 3 \\
	} & 
	\stringlessheap {
	\hs & \hs & \hs & \hs & \hs \\
	& \hs & \hs & \hs & \hs & \hs \\
	\hs & \hs & \hf & \hs & \hs \\
	& \hf & \hs & \hs & \hs & \hs \\
	1 & 2 & 3  \\
	} & 
	\stringlessheap {
	\hs & \hs & \hs & \hs & \hs \\
	& \hs & \hs & \hs & \hs & \hs \\
	\hf & \hs & \hf & \hs & \hs & \hs \\
	& \hf & \hs & \hs & \hs & \hs \\
	1 & 2 & 3  \\
	} & 
	\stringlessheap {
	\hs & \hs & \hs & \hs & \hs \\
	& \hf & \hs & \hs & \hs & \hs \\
	\hf & \hs & \hf & \hs & \hs & \hs \\
	& \hf & \hs & \hs & \hs & \hs \\
	1 & 2 & 3  \\
	} \\
\end{tabular}
\end{center}
We can view the points in the lattice as the vertices in the Hasse
diagram for the heap poset where the edges are implied by the Coxeter
graph.  Note that the other reduced expression $s_2 s_1 s_3 s_2$ for
$[3412]$ corresponds to a different linear extension of the heap
above.  
\end{example}

\begin{remark}
In the lattice representation of a heap poset, all of the entries of
the reduced expression that correspond to the same generator lie in a
column over the given generator in the Coxeter graph.  Each entry will
have a certain \em level \em in the heap, but the poset is not ranked.
In the example below, the reduced expression $\w = s_1 s_4 s_2 s_3$
has a heap where the rank of the $s_3$ entry is not well defined and
the level of the $s_4$ entry is an artifact of the way we imposed
``gravity'' in the construction.
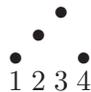
\[
        \stringlessheap {
        & \hs & \hs & \hs & \hs & \hs \\ 
        \hs & & \hf & \hs & \hs & \hs & \hs \\
        & \hf & \hs & \hs & \hs & \hs & \hs \\
        \hf & \hs & \hs & \hf & \hs & \hs \\
        1 & 2 & 3 & 4 \\ 
        }
\]
In Section~\ref{s:convex}, we will further refine the lattice
representation of the heap by coalescing the connected components so
entries that are connected in $G \times \N$ satisfy the covering
relation in the heap poset.
\end{remark}
 
\bigskip

In type $A$, the heap construction can be combined with another
combinatorial model for permutations in which the entries from the
1-line notation are represented by strings.  The points where two
strings cross can be viewed as adjacent transpositions of the 1-line
notation.  Hence, we can overlay strings on top of a heap diagram to
recover the 1-line notation for the element, by drawing the strings
from bottom to top so that they \textit{cross} at each entry in the
heap where they meet and \textit{bounce} at each lattice point not in
the heap.  Conversely, each permutation string diagram corresponds
with a heap by taking all of the points where the strings cross as the
``fat'' points of the heap and letting them ``fall'' according to the
relations given by the Coxeter graph.

For example, we can overlay strings on the two heaps of $[3214]$.
Note that the labels in the picture below refer to the strings, not
the generators.
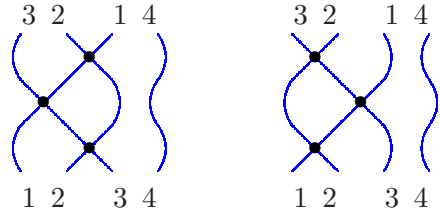
\begin{center}
\begin{tabular}{ll}
	\heap{
	& 3 \ \ 2 &  & 1 \ \ 4 & \\
	\StringR{\hs} & \hs & \StringLRX{\hf} & \hs & \StringL{\hs} \\
	\hs & \StringLRX{\hf} & \hs & \StringLR{\hs} & \hs \\
	\StringR{\hs} & \hs & \StringLRX{\hf} & \hs & \StringL{\hs} \\
	& 1 \ \ 2 &  & 3 \ \ 4 & \\
	} & 
	\heap{
	& 3 \ \ 2 &  & 1 \ \ 4 &  \\
	\hs & \StringLRX{\hf} & \hs & \StringLR{\hs} \\
	\StringR{\hs} & \hs & \StringLRX{\hf} & \hs & \StringL{\hs} \\
	\hs & \StringLRX{\hf} & \hs & \StringLR{\hs} \\
	& 1 \ \ 2 &  & 3 \ \ 4 &  \\
	} \\
\end{tabular}
\end{center}

The string diagram helps us to visualize the relationship between the
1-line notation for a permutation and the corresponding heap.  Fixing
a reduced expression $\w$ for a permutation, there exists a string
diagram that is obtained from $Heap(\w)$ by adding strings that cross
at each lattice point of $Heap(\w)$.  Observe that we are able to read
off the 1-line notation for an element by labeling the strings $1,
\dots, n$ along the bottom and then reading the corresponding labels
from the top.  We can also obtain reduced expressions from any string
diagram by reading the string crossings as generators in any order
that is consistent with the implied heap poset structure.  

\begin{example}\label{e:b_n}
The Coxeter graph of type $B_{n}$ is of the form
\[
\xymatrix @-1pc {
\gn_0 & \ar@{-}[l]_{4} \gn_1 \ar@{-}[r] & \gn_2 \ar@{-}[r] & \gn_3 \ar@{-}[r] 
& \ldots - \gn_{n-1}\\
} .
\]
From this graph, we see that the symmetric group $S_{n}$ is a
parabolic subgroup of this Coxeter group.  That is, $S_n$ is generated
by a subset of the generators of $B_n$.  Because of this, the elements of this
group have a standard 1-line notation in which a subset of the
entries are \textit{barred}.  We often think of the barred entries as
negative numbers, and this group is referred to as the group of
\textit{signed permutations} or the \textit{hyperoctahedral group}.
The action of the generators on the 1-line notation is the same for
$\{s_1, s_2, \dots, s_{n-1} \}$ as in type $A$ in which $w s_i$
interchanges the entries in positions $i$ and $i+1$ in the 1-line
notation for $w$.  The $s_0$ generator acts on the right of $w$ by
changing the sign of the first entry in the 1-line notation for $w$.
For example, $w=[\bar{4}2\bar{3}1]$ is an element of $B_{4}$ and
\begin{align*}
ws_{0}=[42\bar{3}1]\\
ws_{1}=[2\bar{4}\bar{3}1].
\end{align*}
Note that because the edge in the $B_{n}$ Coxeter graph connecting
$s_{0}$ and $s_{1}$ is labeled 4, we have that $s_0 s_1 s_0 s_1$ and
$s_1 s_0 s_1 s_0$ are reduced expressions for the same element denoted
$[\bar{1}\bar{2}34\ldots n]$ in 1-line notation.

The heap for a type $B$ reduced expression will look  like the
heap of a type $A$ expression because its Coxeter graph is a path.  As
in type $A$, we can adorn these heap constructions with strings that
represent the digits of the 1-line notation for the element.  If we
label the strings at the bottom of the diagram with the numbers from
$1, \dots, n$ then the $s_0$ generator has the effect of bouncing the
string back in the direction from which it came, while changing the
sign of the label for the string.  All other generators cross the
strings as in type $A$.  
\end{example}

\begin{example}\label{ex:type.d}
The Coxeter graph for type $D_{n}$ is shown below.
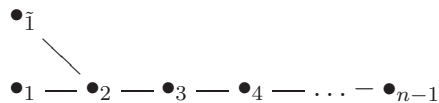
\[
\xymatrix @-1pc {
{\gn}_{\tilde{1}} &        &                &                            &  \\
\gn_1 & \ar@{-}[l] \ar@{-}[ul] \gn_2 \ar@{-}[r] & \gn_3 \ar@{-}[r] & \gn_4 \ar@{-}[r] & \ldots - \gn_{n-1} \\
}
\]
The elements of type $D$ can be viewed as the subgroup of $B_{n}$
consisting of signed permutations with an even number of barred
entries.  The action of the generators on the 1-line notation is the
same for $\{s_1, s_2, \dots \}$ as in type $A$ in which $w s_i$
interchanges the entries in positions $i$ and $i+1$ in the 1-line
notation for $w$.  The $s_{\tilde{1}}$ generator acts on the right of
$w$ by marking the first two entries in the 1-line notation for $w$
with bars and interchanging them.  For example, $v=[42\bar{3}\bar{1}]$
and $w=[\bar{4}2\bar{3}1]$ are elements of $D_{4}$ and
\begin{align*}
v s_{\tilde{1}} =& [\bar{2}\bar{4}\bar{3}\bar{1}]\\
w s_{\tilde{1}} =& [\bar{2}4\bar{3}1] .
\end{align*}
Although the Coxeter graph for type $D$ has a fork, we will draw the
heap of type $D$ elements in a linearized way by allowing entries in
the first column to consist of either generator or both:
\begin{align*}
s_{1}=& \hf\\
s_{\tilde{1}}=& \tilde{\hf} \\
s_{1}s_{\tilde{1}}=& \hf \tilde{\hf} .
\end{align*}
We denote this linearized lattice point representation by $Heap(\w)$.
Hence, if $w \in D_n$ then $Heap(\w)$ is a subset of $[n-1] \times
\N$ rather than $G \times \N$, where $[n-1] = \{1, 2, \dots, n-1\}$.

As in type $A$, we can adorn $Heap(\w)$ with strings that represent
the digits of the 1-line notation for the element.  If we label the
strings at the bottom of the diagram with the numbers from $1, \dots,
n$, then the $s_{\tilde{1}}$ generator crosses the strings that
intersect it and changes the sign on the labels for both strings.
All other generators simply cross the strings as in type $A$.  For
example, the heap of the reduced expression
$s_{\tilde{1}}s_{2}s_{3}s_{1}s_{2}s_{\tilde{1}}s_{1} =[\bar{3}\bar{4}\bar{2}\bar{1}]$
is given below.
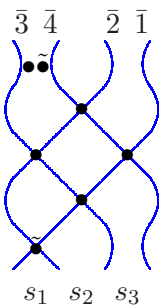
\[
	\heap {
	& \bar{3} \ \ \bar{4} & & \bar{2} \ \ \bar{1} &  &  \\
	& \StringLR{\hf \tilde{\hf}} & \hs & \StringLR{\hs} & \hs & \hs & \hs \\
	\StringR{\hs} & & \StringLRX{\hf} & \hs & \StringL{\hs} & \hs & \hs & \hs \\
	& \StringLRX{\hf} & {\hs} & \StringLRX{\hf} & \hs & \hs & \hs \\
	\StringR{\hs} & & \StringLRX{\hf} & \hs & \StringL{\hs} & \hs & \hs & \hs \\
	& \StringLRX{\tilde{\hf}} & \hs & \StringLR{\hs} & \hs & \hs & \hs \\
	& s_1 & s_2 & s_3 \\
	}
\]
\end{example}

\bigskip
\section{Classical pattern avoidance}\label{s:classical.patterns}

The 1-line notations for types $A, B, D$ carry a notion of pattern
containment that generalizes the following classical definition.
\begin{definition}
Let $w=[w_{1} \ldots w_{n}]$ be a permutation in $S_{n}$ written in
1-line notation as described in Example~\ref{ex:typeA}.  Let $p=[p_1
\ldots p_k]$ be another permutation in $S_{k}$ for $k\leq n$.  Then we
say $w$ \textit{contains the permutation pattern} $p$ if there exists a
subsequence $1\leq i_{1}<i_{2}<\ldots<i_{k}\leq n$ such that
\begin{equation*}
w_{i_{a}} < w_{i_{b}}  \  \iff  \  p_{a} < p_{b}
\end{equation*}
for all $1 \leq a < b \leq k$.  If $w$ does not contain $p$ then we
say that $w$ \em avoids the permutation pattern \em $p$.
\end{definition}
In other words, $w$ contains $p$ if there exist $k$ rows and columns
in the permutation matrix for $w$ whose common entries are the
permutation matrix for $p$.  For example, $w=[\ul{5}32\ul{4}\ul{1}]$
contains the pattern $p=[321]$ in several ways including the
underlined subsequence, and $w$ avoids $q=[1234]$.

One of the earliest uses of permutation patterns occurred in computer
science \cite{tarjan}.  A good introduction to enumerative methods in
pattern avoidance can be found in \cite{bona.book}.  Several
interesting properties of Schubert varieties, Kazhdan--Lusztig
polynomials and Bruhat order can be characterized by pattern avoidance
\cite{la-s,b2,BiBr,gasharov97,w-y,w-y2,t3,bmb}.

A property of permutations can be characterized by pattern avoidance
if there is no $w$ having the property and $p$ not having the property
such that $w$ contains $p$ as a pattern.  Equivalently, the subset of
permutations that do not have the property must be an upper order
ideal in the poset on $S_{\infty} = \bigcup_{n>0} S_{n}$ ordered by
pattern containment.  The property is then characterized by avoiding
the minimal elements of this upper order ideal.

For example, the property of a permutation being short braid avoiding
in type $A$ is characterized by avoiding the pattern $[321]$, first
noted by \cite{b-j-s}.  Hence, the fully commutative elements in
$S_{n}$ are enumerated by the Catalan numbers \cite{simion-schmidt}.  Also,
checking a permutation of length $n$ for a subsequence of length $3$
can be done in $O(n^{3})$ time while looking at all reduced words for
a typical permutation takes an exponential amount of time in $n$.  A
pattern avoidance characterization is most useful when the set of
minimal patterns is finite, but this need not be the case
\cite{Bona-Spielman}.

When $w=[w_1 \ldots w_n]$ is a signed permutation from type $B$ or
type $D$, then we say that $w$ \em contains the 1-line pattern \em
$p=[p_1 \ldots p_k]$ if the underlying permutation for $w$, obtained by
ignoring the bars in the 1-line notation, contains the underlying
permutation for $p$, and assuming the pattern instance occurs in
positions $i_1 < i_2 < \dots < i_k$, we further require that $w_{i_j}$
is barred only if $p_j$ is barred.  For example, the type $D$ element
$w = [\bar{5} 3 2 \bar{4} 1]$ contains the pattern $p = [\bar{3}
\bar{2} 1]$ as the underlined subsequence $[\ul{\bar{5}} 3 2
\ul{\bar{4}} \ul{1}]$, but \em not \em the underlined subsequence
$[\ul{\bar{5}} 3 \ul{2} \bar{4} \ul{1}]$ because the pattern of bars
does not match.

Some early applications of pattern avoidance in types $B$ and $D$
occurred in \cite{billey,billey-lam}.  Type $B$ enumerative results
have been obtained by \cite{beck,simion,mansour-west}, and there are
also some extensions to colored permutations \cite{mansour2,
mansour3}.  Classical pattern avoidance extends to all Coxeter groups
using the notion of root subsystems described in \cite{b2, BiBr}.

\bigskip
\section{Deodhar elements of Coxeter groups}\label{s:Deodhar.elements}

Let $W$ be an arbitrary Coxeter group with $w,y \in W$.  We say that
$w$ contains $y$ \textit{as a factor} if there exist elements $a$ and
$b$ in $W$ such that $w=ayb$ and $l(w)=l(a)+l(y)+l(b)$.  Equivalently,
$w$ contains $y$ as a factor if some reduced expression for $w$
contains some reduced expression for $y$ as a consecutive subword,
i.e. $\w=\w_1\w_2...\w_p$ and $y = \w_{i}\w_{i+1}\cdots \w_{j}$ for
some $1\leq i \leq j\leq p$.  The induced partial order on Coxeter
group elements is known as the \textit{two-sided weak Bruhat order}
\cite{b-b}.  We will use the following lemma to show that the Deodhar
elements form a lower order ideal in the two-sided weak Bruhat order.

\begin{lemma}\label{l:w.winverse}
Let $W$ be a Coxeter group.  If $w \in W$ is Deodhar then $w^{-1}$ is
Deodhar.
\end{lemma}
\begin{proof}
Let $\w=\w_1 \ldots \w_p$ be a reduced expression for $w$, and consider the two products
\begin{align}
 C'_{\w_1} C'_{\w_2}\cdots C'_{\w_p} = & q^{-\frac{p}{2}}\sum R_{u}(q) T_{u} \\
 C'_{\w_p}\cdots C'_{\w_2} C'_{\w_1} = & q^{-\frac{p}{2}}\sum S_{u}(q) T_{u}
\end{align}
It follows directly from the symmetry in the multiplication rule \eqref{e:hecke.def} that 
\[
R_{u}(q) = S_{u^{-1}}(q). 
\]
By Deodhar's Theorem~\ref{t:deodhar}, $w$ is Deodhar if and
only if
\begin{align}
C'_{w} & = C'_{\w_1} C'_{\w_2}\cdots C'_{\w_p} \\
& \implies R_{u}(q) = P_{u,w}(q) \hspace{.2in} \text{  for all $u$ } \\
& \implies S_{u^{-1}}(q) = P_{u,w}(q) \hspace{.2in} \text{for all $u$ } \\
& \implies \mathrm{deg}(S_{u^{-1}}(q)) \leq \frac{l(w^{-1})-l(u^{-1})-1}{2} \text{for all $u$ and } S_{w^{-1}}(q) =1 \\
& \implies C'_{\w_p} \cdots C'_{\w_2} C'_{\w_1} =C'_{w^{-1}}
\end{align}
by definition of $C'_{w^{-1}}$ since the product of $C'_{\w_i}$'s is
invariant under the involution.  Therefore, $w^{-1}$ is Deodhar.  
\end{proof}

The following proposition can easily be derived from Proposition 2.1.4
of \cite{f-g}.  We include an independent proof for completeness.
\begin{proposition}\label{p:deodhar.factor} 
Let $w, y \in W$ be Coxeter group elements.  If $y$ is not Deodhar and
$w$ contains $y$ as a factor, then $w$ is not Deodhar either.
\end{proposition}
\begin{proof}
Suppose $y$ is not Deodhar and $y=\y_1 \cdots \y_p$ is a reduced
expression.  Then by \eqref{e:deodhar.ineq} there exists a
proper mask $\s$ for $\y_1 \cdots \y_p$ with
\begin{equation*}
\text{\# of zero-defects of } \s \geq \text{\# of plain-zeros of } \s.
\end{equation*}
Consider multiplying $y$ on the right by a generator $s \in S$ such
that $l(y s)>l(y)$.  We can extend the mask $\s$ by placing a 1 in the
last position to obtain a mask for $\y_1 \cdots \y_p s$.  Since the
inequality above remains unchanged for the new mask, the element $y s$
is not Deodhar.

Next, consider multiplying $y$ on the left by a generator $s \in S$
such that $l(sy)>l(y)$ or equivalently $l(y^{-1}s)>l(y^{-1})$.  By
Lemma~\ref{l:w.winverse}, $y$ is not Deodhar implies $y^{-1}$ is not
Deodhar.  So by the argument above, $y^{-1}s$ is not Deodhar and so
neither is $sy$.

The theorem follows by induction on $l(w)-l(y)$.
\end{proof}

We immediately obtain a combinatorial proof of a result from
\cite{f-g} which shows that the Deodhar elements are all short braid
avoiding.

\begin{corollary}\label{c:short.braid}
Let $w \in W$ be a Coxeter group element.  If $w$ contains a short
braid then $w$ is not Deodhar.
\end{corollary}
\begin{proof}
Say $s_{i},s_{j}$ are noncommuting generators.  The reduced
expression/mask pair
\[ \block{ s_i & s_j & s_i \\ 1 & 0 & 0 } \]
has a zero-defect in the last position hence the number of
zero-defects equals the number of plain zeros.  Consequently,
$s_{i}s_{j}s_{i}$ is not Deodhar.  Therefore, any element that
contains a short braid as a factor in a reduced expression is not
Deodhar.
\end{proof}

Proposition~\ref{p:deodhar.factor} also shows that the non-Deodhar
elements form an upper order ideal in the two-sided weak Bruhat order.
In order to obtain an efficient generating set for this ideal, we
consider a refinement of factor containment.

\begin{definition}\label{d:coxeter.map}
Let $W, W'$ be Coxeter groups with associated Coxeter graphs $G, G'$
respectively.  Then, a \textit{Coxeter embedding of} $G'$ is an
injective map of the generators $f: G' \rightarrow G$ that restricts
to a labeled graph isomorphism onto its image.  

A Coxeter embedding induces an injection of $W'$ into $W$, and we will
abuse notation and call this map $f: W' \rightarrow W$ a Coxeter
embedding also.  When phrased algebraically, a Coxeter embedding $f:
W' \rightarrow W$ is an injective map of generators for the Coxeter
group such that $m(s_i, s_j) = m(f(s_i), f(s_j))$ for all $s_i, s_j
\in W'$.  Since $f$ is a map of generators, we can extend it to a map
of Coxeter group elements by treating it as a word homomorphism on any
reduced expression in $W'$.
\end{definition}

\begin{definition}\label{d:embedded.factor}
Suppose $W$ is a Coxeter group, and $w, y \in W$.  Let $W_{y}$ be the
parabolic subgroup whose generators are determined by the support of
$y$.  If there exists a Coxeter embedding $f:W_{y} \rightarrow W$ such
that $w$ contains $f(y)$ as a factor, then we say that $w$ contains
$y$ \textit{as an embedded factor}.
\end{definition}

This definition yields a stronger reformulation of
Proposition~\ref{p:deodhar.factor}, which enables us to characterize
Deodhar elements with a shorter list of non-Deodhar patterns.

\begin{corollary}
Let $y$ be a Coxeter element that is not Deodhar.  If $w$ is a Coxeter
element that contains $y$ as an embedded factor, then $w$ is not
Deodhar either.
\end{corollary}
\begin{proof}
If $\s$ is a mask for a reduced expression $\y$ then it follows from
Definition~\ref{d:coxeter.map} that $\s$ is also a mask for $f(\y)$.
This mask has defects in exactly the same positions as it does when it
is applied to $\y$.  Hence, $y$ is Deodhar if and only if $f(y)$ is
Deodhar, for any Coxeter embedding $f: W_{y} \rightarrow W$.
\end{proof}

\begin{example}
In the Coxeter group of type $B_{4}$,
$[2\bar{4}51\bar{3}]=
s_{1}s_{2}s_{3}s_{1}s_{0}s_{1}s_{2}s_{1}s_{0}s_{4}s_{3}s_{4}s_{1}$
contains the factor $s_{4}s_{3}s_{4}$ in the parabolic subgroup
generated by $\{ s_{3},s_{4} \}$.  This subgroup is isomorphic
to $S_{3}$ and $s_{4}s_{3}s_{4}$ maps to $s_2 s_1 s_2 = [321] \in S_{3}$, so
$[2\bar{4}51\bar{3}]$ contains $[321]$ as an embedded factor.
Consequently, it is not Deodhar.
\end{example}

\begin{example}\label{ex:embeddings}
In type $A$, the Coxeter embeddings of connected subgraphs are simply
shifts of the generators along the linear Coxeter graph or reversed
shifts.  In particular, if the generators of $S_{k}$ are labeled
$s_1$, $s_2$, \dots, $s_{k-1}$ in its Coxeter graph, then the image of
the generators under a Coxeter embedding $f:S_{k} \rightarrow S_{n}$
are either of the form $s_{1 + j}, s_{2 + j},\ldots, s_{k-1 + j}$ or
$s_{k-1 + j},\ldots,s_{2+j} s_{1+j}$ for some $0\leq j \leq n-k$.
\end{example}

The \textit{hexagon avoiding} elements of a Coxeter group are the ones
that avoid the element
\begin{equation*}
  \label{hexuj}
  u = s_{3} s_{2}s_{1}s_{5}s_{4}s_{3}s_{2}
  s_{6}s_{5}s_{4}s_{3}s_{7}s_{6}s_{5} = [46718235]
\end{equation*}
of $A_7$ as an embedded factor.  The name arises from the shape
formed by the heap of $u$:
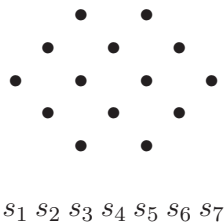
\[ \xymatrix @=-4pt @! {
& {\hs} & & {\hs} & & {\hs} & & {\hs} \\
& & {\hf} & & {\hf} & & {\hs} & \\
& {\hf} & & {\hf} & & {\hf} & & {\hs} \\
{\hf} & & {\hf} & & {\hf} & & {\hf} & \\
& {\hf} & & {\hf} & & {\hf} & & {\hs} \\
& & {\hf} & & {\hf} & & {\hs} & \\
& {\hs} & & {\hs} & & {\hs} & & {\hs} \\
s_1 & s_2 & s_3 & s_4 & s_5 & s_6 & s_7 & &  \\
 &  &  &  &  &  &  & &  \\
} \]

Billey and Warrington have given a complete characterization of the
Deodhar elements in \textit{linear Weyl groups} where the Coxeter
graph consists of a single path.  We can now state their theorem
precisely.

\begin{theorem}\cite{b-w}\label{r:typeahex}
In types $A_{n},B_{n}$, an element $w$ is Deodhar if and only if $w$
avoids short braids and hexagons as embedded factors.
\end{theorem}

Our main theorem below generalizes this theorem.  It is a concise
characterization of the Deodhar condition for the other finite Weyl
groups. 

\begin{remark}
The techniques used in \cite{b-w} do not easily extend to the
remaining finite Weyl groups of types $D$ and $E$.  Among the
connected short-braid avoiding heaps of type $A$, the notions of
coalesced heap containment and embedded factor containment are
essentially the same up to the orientation of the Coxeter graph.  This
is not the case in type $D$.
In addition, the non-Deodhar elements of type $D$ can
have heaps containing ``alcoves'' or ``holes'' so that the heap
lattice points do not form a laterally convex set in the sense of
\cite{b-w}.
\end{remark}

\begin{example}\label{e:infantichain}
In type $D$, there is an infinite antichain of non-Deodhar elements
(i.e. no pair of elements from the family contain each other as
embedded factors) whose heaps can contain ``alcoves:''

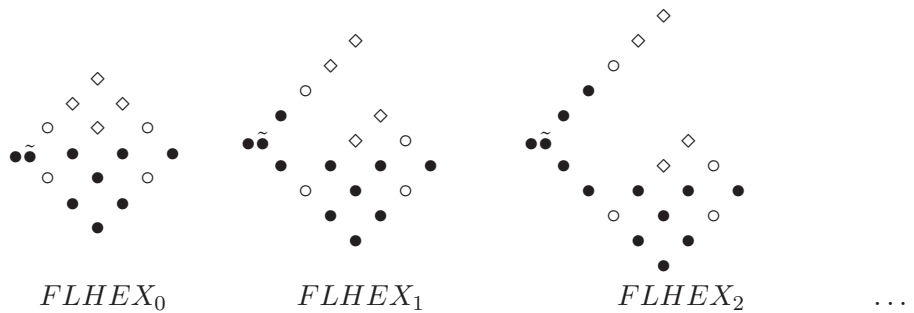
\begin{equation*}
\begin{matrix}
\begin{matrix}
	\xymatrix @=-8pt @! {
{\hs} & & {\hs} & & {\hs} & & {\hs} \\
& {\hs} & & {\hd} & & {\hs} &\\
{\hs} & & {\hd} & & {\hd} & & {\hs} \\
& {\hz} & & {\hd} & & {\hz}\\
{\hf \tilde{\hf} }&  & {\hf} & & {\hf} & & {\hf} & {\hs} \\
& {\hz} & & {\hf} & & {\hz} &\\
{\hs} & & {\hf} & & {\hf} & & {\hs} \\
& {\hs} & & {\hf} & & {\hs} &\\
    } 
\end{matrix}
&
\begin{matrix}
	\xymatrix @=-8pt @! {
 {\hs} & & {\hs} & & {\hd} & & {\hs} & \\
 & {\hs} & & {\hd} & & {\hs}& \\
 {\hs} & & {\hz} & & {\hs}& & {\hs}& & {\hs} & \\
 & {\hf} & & {\hs} & & {\hd} & \\
{\hf \tilde{\hf} } & & {\hs}& & {\hd} & & {\hz} & \\
  & {\hf} & & {\hf}& & {\hf} & & {\hf} \\
 {\hs} & & {\hz}& & {\hf} & & {\hz} & \\
 & {\hs} & & {\hf}& & {\hf} & \\
 & {\hs} & & {\hs}&  {\hf} & & \\
    } 
\end{matrix} 
&
\begin{matrix} 
	\xymatrix @=-8pt @! {
& {\hs} & & {\hs} & &  {\hd} &  & {\hs}  & & {\hs} &  \\
 {\hs} & & {\hs} & & {\hd} & & {\hs} & & {\hs} & & {\hs} \\
 & {\hs} & & {\hz} & & {\hs}& & {\hs} & & {\hs} \\
 {\hs} & & {\hf} & & {\hs}& & {\hs}& & {\hs} & & {\hs} & & {\hs} \\
 & {\hf} & & {\hs} & & {\hs} & & {\hs}& & {\hs} \\
{\hf \tilde{\hf}} & & {\hs}& & {\hs} & & {\hd} & & {\hs} & & {\hs} \\
  & {\hf} & & {\hs}& & {\hd} & & {\hz} & & {\hs} \\
 {\hs} & & {\hf}& & {\hf} & & {\hf} & & {\hf} & & {\hs} \\
 & {\hs} & & {\hz}& & {\hf} & & {\hz} & & {\hs} \\
 & {\hs} & & {\hs}&  {\hf} & &  {\hf} & & {\hs}  &\\
 & {\hs} & & {\hs}&  {\hs} & {\hf}& & {\hs} & & {\hs} \\
    } 
\end{matrix} \\
FLHEX_{0} & FLHEX_{1} & FLHEX_{2} & \ldots . \\
 & & & \\
\end{matrix}
\end{equation*}
We have drawn the heaps of these elements with decorations that
indicate a particular mask, as described in Section~\ref{s:d.heaps}.
The masks shown demonstrate that these elements are not Deodhar.
\end{example}

This example shows that the set of type $D$ Coxeter elements partially
ordered by embedded factor containment is not well quasi-ordered.
Also, there is a simple example showing that the permutations
partially ordered by embedded factor containment is not well
quasi-ordered.  See \cite{Bona-Spielman} for an analogous example in
the classical permutation pattern poset.

We can still obtain a finite characterization for the Deodhar
condition in type $D$ since all of the $FLHEX_k$ elements contain
$FLHEX_{0} = [\bar{1} 6 7 8 \bar{5} 2 3 4]$ as a classical
1-line pattern, in the manner described in
Section~\ref{s:classical.patterns}.

\begin{theorem}\label{t:main}
Let $w \in W$ be an element of a finite Weyl group.  Then, the following
are equivalent:
\begin{enumerate}
\item The element $w$ is Deodhar.
\item The element $w^{-1}$ is Deodhar.
\item The element $w$ avoids the short list of embedded factor
	patterns given in Figure~\ref{f:patterns}, as well as the
	$FLHEX_0$ 1-line pattern of type $D$.
\end{enumerate}
\end{theorem}
\begin{proof}
The equivalence of (1) and (2) follow from Lemma~\ref{l:w.winverse}.
The equivalence of (1) and (3) follows from Theorem~\ref{r:typeahex}
for types $A$ and $B$, Theorem~\ref{t:main.d} below for type $D$, and
Theorem~\ref{t:main.e} for the finite exceptional groups.  This
accounts for all irreducible finite Weyl groups. 
\end{proof}

The proof will occupy Sections~\ref{s:d.heaps}
through~\ref{s:classification.e}.  For the finite exceptional groups,
a brute-force search implemented on a computer suffices.  Our work is
simplified by the fact that we only need to check the short braid
avoiding elements of these groups.  For the infinite type $D$ family,
we need to show that our list of minimal non-Deodhar elements is
complete.  This is shown in Theorems~\ref{t:minimal.patterns} and
\ref{sbadh.th}.  The proof of Theorem~\ref{t:minimal.patterns}
involves a map from the heap of a type $D$ reduced expression to a
type $A$ heap where the Deodhar condition was already known by
Theorem~\ref{r:typeahex}.  We complete the classification by checking
that this map preserves the Deodhar property.


\bigskip
\section{Short-braid avoiding heaps in type $D$}\label{s:d.heaps}

This section develops the heap technology necessary to carry out the
classification of minimal non-Deodhar elements in type $D$ under
embedded factor containment.  For our work in this section, it
suffices by Corollary~\ref{c:short.braid} to consider only short braid
avoiding elements.  Short braid avoiding elements are fully
commutative so they have a unique heap poset.  We draw the heaps of
type $D$ elements in a linearized way as described in
Example~\ref{ex:type.d}, with entries corresponding to both $s_1$ and
$s_{\tilde{1}}$ generators in the first column, and denote this
lattice point representation of the heap by $Heap(w)$.  We will
consider masks to be assignments of 0's and 1's to the entries in the
heap instead of 0,1-sequences associated to a particular reduced
expression.  Of course, a reduced word/mask pair can be read off from
the heap by reading the entries in order of any linear extension of
the heap.  

We decorate $Heap(w)$ according to mask-value using the following
table: 
\bigskip
\begin{center}
\begin{tabular}{|p{0.7in}|p{3in}|} 
\hline
Decoration & Mask-value \\
\hline
$\hd$  &	zero-defect entry\\
$\hz$  &	plain-zero entry (not a defect)\\
$\hf$  &	mask-value 1 entry\\
$\ha$  &	entry of the heap with unknown mask-value\\
$\hb$  &	lattice point not necessarily in the heap\\
$\hv$  &	lattice point that is definitely not in the heap\\
$\hh$  &	lattice point that is highlighted for emphasis\\
\hline
\end{tabular}
\end{center}
\bigskip In the decorated heaps, we don't distinguish one-defects from
plain-ones since they don't contribute to the Deodhar bound
\eqref{e:deodhar.ineq}.

We can adorn our decorated heaps with strings representing the digits
of the standard 1-line notation for type $D$ elements as described in
Section \ref{s:heaps}.  Given a decorated heap, the strings will cross
only at entries corresponding with mask-value 1.  At mask-value 0
entries the strings ``bounce'' as if the entry were not in the heap.
If the decorated heap corresponds to a reduced expression $\w$ and
mask $\sigma$, then the resulting labels on the strings at the top of
the heap will be the 1-line notation for $\w^{\sigma}$.  Recall that a
string passing through a $\tilde{\hf}$ entry corresponds to a
$s_{\tilde{1}}$ generator so it changes sign (which is not shown
explicitly in our pictures).

We can use the strings to obtain a useful test for the defect status
of a particular entry in a decorated heap.  Note that at every entry
in the heap, two strings approach it from either side.  We will call
these the \textit{left string} and the \textit{right string} for that
entry.  

\begin{lemma}\label{l:defectlabellings}
Consider an entry $p$ in a heap.  Draw the left and right strings
emanating downward from $p$ and label the string that ends up
leftmost on the bottom by 1, and the string that ends up rightmost on
the bottom by 2.  The entry $p$ is a defect if and only if the strings
are labeled at the top according to the following table of patterns.
\bigskip
\begin{center}
\begin{tabular}{|p{2.3in}|p{2.8in}|} 
\hline
\text{If $p$ corresponds to the generator \dots} & \text{\dots then $p$ is a defect if and only
if the strings} \text{are  top-labeled } \\
\hline
$s_{\tilde{1}}$ & \vspace{0.02in} $1 \bar{2}$, $\bar{1} \bar{2}$, $\bar{2} 1$ \text{ or } $\bar{2} \bar{1}$ \\
\hline
\text{any other generator} & \vspace{0.02in} $1 \bar{2}$, $\bar{1} \bar{2}$, $2 1$ \text{ or } $2 \bar{1}$ \\
\hline
\end{tabular}
\end{center}
\bigskip
\end{lemma}

Note that when $p$ is not $s_{\tilde{1}}$, this test is just a signed
version of the usual type $A$ inversion.

\begin{proof}
This follows because the length formula in type $D$ for an even
signed permutation given in one line notation $w = [w_1 w_2 \dots
w_n]$ is $l(w) = \#\{i < j: w_i > w_j\} + \#\{i < j :
\overline{w_i} > w_j\}$, viewing the barred entries as negative.
\end{proof}

\begin{example}\label{ex:defect.status}
The following decorated heaps with strings demonstrate the defect
status of the top entry in the masked expressions below.
$$
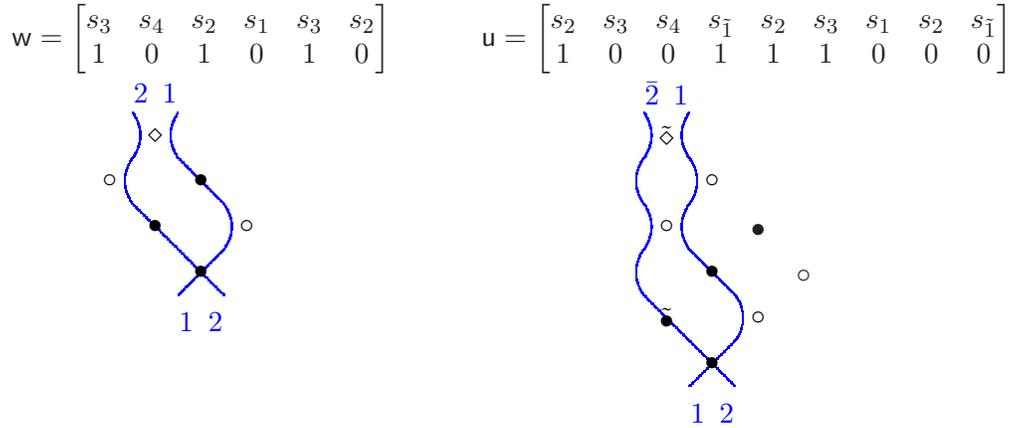
\begin{array}{ccc}
\w = \block{ s_3 & s_4 & s_2 & s_1 & s_3 & s_2 \\ 1 &
0 & 1 & 0 & 1 & 0 } & &
\u = \block{ s_2 & s_3 & s_4 & s_{\tilde{1}} & s_2 & s_3 & s_1 & s_2
& s_{\tilde{1}} \\ 1 & 0 & 0 & 1 & 1 & 1 & 0 & 0 & 0 }\\
 \heap {
	{\hs} & & {\color{blue} 2 {\hs} \ \ 1 } & & {\hs} & & {\hs} \\
	{\hs} & & \StringLR{\hd} & & {\hs} & & {\hs} \\
	& \StringR{\hz} & & \StringLX{\hf} & & {\hs} & \\
	& & \StringLX{\hf} & & \StringL{\hz} & & {\hs} \\
	& {\hs} & & \StringLRX{\hf} & & {\hs} & \\
	& & {\hs} & { \color{blue} 1 \ \ 2 } & {\hs} & & {\hs} \\
	} 
& \hspace{.2in}&
    \heap {
    	{\hs} & {\color{blue} \bar{2} \ \ 1 }  & {\hs} & & {\hs} & & {\hs} \\
	& \StringLR{\tilde{\hd}} & & {\hs} & & {\hs} & \\
	\StringR{\hs} & & \StringL{\hz} & & {\hs} & & {\hs} \\
	& \StringLR{\hz} & & {\hf} & & {\hs} & \\
	\StringR{\hs} & & \StringLX{\hf} & & {\hz} & & {\hs} \\
	& \StringLX{\tilde{\hf}} & & \StringL{\hz} & & {\hs} & \\
	& & \StringLRX{\hf} & & {\hs} & & {\hs} \\
	& & { \color{blue} 1 {\hs} \ \ 2 } & & {\hs} & & {\hs} \\
	} 
\end{array}
$$
\end{example}

Here is an observation that is used extensively in the classification.

\begin{corollary}\label{c:string.crossing}
The strings for a defect in any type $D$ heap must cross at least once strictly
below the defect.
\end{corollary}
\begin{proof}
If the strings for an entry never cross then the right string of the
entry will be labeled $2$, which does not match any of the defect
labelings in Lemma~\ref{l:defectlabellings}.
\end{proof}

\begin{example}\label{ex:not.critical}
In type $A$, an entry is a defect if and only if its strings cross an
odd number of times below.  By contrast, in type $D$ it is possible
for the strings of an entry to cross below, yet not form a defect.
For example,
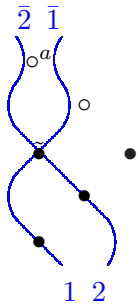
\[ \heap {
	{\hs} & {\color{blue} \bar{2} \ \ \bar{1}} & {\hs} & & {\hs} & & {\hs} \\
	& \StringLR{\hz^a} & & {\hs} & & {\hs} & \\
	\StringR{\hs} & & \StringL{\hz} & & {\hs} & & {\hs} \\
	& \StringLRX{\tilde{\hf}} & & {\hf} & & {\hs} & \\
	\StringR{\hs} & & \StringLX{\hf} & & {\hs} & & {\hs} \\
	& \StringLX{\hf} & & \StringL{\hs} & & {\hs} & \\
	& {\hs} & {\color{blue} 1 \ \ 2 } & {\hs} & & {\hs} & \\
} 
\] 
In this case, the strings are interchanged but both negatively signed,
so the values are in increasing order.  Hence, $a$ is not a defect by
Lemma~\ref{l:defectlabellings}.
\end{example}


Let $\w \in D_{n}$ be a reduced expression for a connected short-braid
avoiding element.  Suppose $x$ and $y$ are a pair of entries in
$Heap(w)$ that correspond to the same generator $s_i$, so they lie in
the same column $i$ of the heap (setting $i = 1$ in case the generator
is $s_{\tilde{1}}$).  Assume that $x$ and $y$ are a \em minimal pair
\em in the sense that there is no other entry between them in column
$i$.  Then, for $\w$ to be reduced there must exist at least one
non-commuting generator between $x$ and $y$, and for $\w$ to be
short braid avoiding there must actually be two non-commuting
generators that lie strictly between $x$ and $y$ in $Heap(w)$.  We
call these two non-commuting generators a \em resolution \em of the
pair $x,y$.  

\begin{definition}
If both of the generators in a resolution lie in column $i-1$ ($i+1$,
respectively), we call the resolution a \em left (right, respectively)
resolution\em.  If the generators lie in distinct columns, we call the
resolution a \em distinct resolution\em.
\end{definition}

In type $D$, every resolution of a minimal pair must be one of these
types.  Note that the $s_{1}$ and $s_{\tilde{1}}$ generators lie in
the same column of the type $D$ heap so although
$s_{2}s_{1}s_{\tilde{1}}s_{2}$ is a fully commutative element in
$D_{3}$, the pair of $s_{2}$ entries do not have a distinct
resolution.  On the other hand, the pair of $s_2$ entries in
$s_{1}s_{2}s_{\tilde{1}}s_{3}s_{2}s_{1}$ has a distinct resolution,
while the pair of $s_{1}$ generators has only a right resolution.

Recall that by Tits' theorem a permutation is fully commutative if and
only if every minimal pair has a distinct resolution.  Short braid
avoiding is equivalent to fully commutative in type $D$.  We establish
some structural lemmas about the resolutions in $Heap(w)$ when $w$ is
short braid avoiding.

\begin{lemma}\label{sbadh.1} 
Let $w\in D_{n}$ be a short braid avoiding element.  Then, no minimal
pair of generators in column $i \geq 2$ in $Heap(w)$ can have a right
resolution.
\end{lemma}
\begin{proof} Since $w$ is short braid avoiding, if $x$ and $y$ are
resolved by a pair of generators $z_1, z_2$ that lie to the right of
column $i$ then $z_1$ and $z_2$ necessarily correspond to the same
generator since the Coxeter graph of type $D$ is a path beyond column
1.  Choose $z_1, z_2$ to be a minimal pair in column $i+1$.  Since $x$
and $y$ form a minimal pair, we cannot backtrack when resolving $z_1$
and $z_2$, so $z_1$ and $z_2$ must be separated by another minimal
pair of generators to the right that are non-commuting with $z_1,
z_2$:
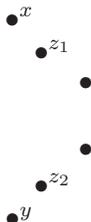
\[ \xymatrix @=-8pt @! {
& {\hs} & & {\hf^x} & & {\hs} & & {\hs} \\
{\hs} & & {\hs} & & {\hf^{z_1}} & & {\hs} \\
& {\hs} & & {\hs} & & {\hf} & & {\hs} \\
{\hs} & & {\hs} & & {\hs} & & {\hs} \\
& {\hs} & & {\hs} & & {\hf} & & {\hs} \\
{\hs} & & {\hs} & & {\hf^{z_2}} & & {\hs} \\
& {\hs} & & {\hf^y} & & {\hs} & & {\hs} \\
} \]
Since $D_{n}$ is finitely generated, eventually there exists a
minimal pair of entries in the rightmost column of $Heap(w)$ that
cannot be resolved.  Hence, every minimal pair of generators in column
$i \geq 2$ is resolved by two generators from distinct columns or by a
pair of generators to the left.
\end{proof}

\begin{lemma}\label{sbadh.2}
Let $w \in D_{n}$ be a short braid avoiding element.
\begin{enumerate}
\item There exists an $s_{\tilde{1}}$ between every minimal pair of
$s_1$ generators, and an $s_1$ between every minimal pair of
$s_{\tilde{1}}$ generators in $Heap(w)$.
\item If there is an entry in which both $s_1$ and $s_{\tilde{1}}$ lie
at the same level of $Heap(w)$, then column 1 contains no other
entries.
\item If there exists a pair of entries in column 1 corresponding to
the same generator then all of the entries in the first column must
be on distinct levels of the heap, and they must alternate between the
generators $s_1$ and $s_{\tilde{1}}$.
\end{enumerate}
\end{lemma}
\begin{proof}
Part (1):  Suppose that $x$ and $y$ are a minimal pair in column 1
corresponding to $s_1$.  Then they must be resolved by a pair $z_1,
z_2$ in column 2 since $s_{1}$ commutes with every other generator.  We
can choose $z_1, z_2$ to be a minimal pair.  Since $z_1, z_2$ cannot
have a right resolution by Lemma~\ref{sbadh.1}, there exists an entry
between $z_1$ and $z_2$ in column 1.  By the minimality of $x,y$, this
entry must correspond to $s_{\tilde{1}}$.  Moreover, $z_1, z_2$ cannot
use a left resolution without contradicting the minimality of the pair
$x$ and $y$.  Hence, $z_1, z_2$ have a distinct resolution with a heap
fragment of the form:
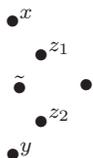
\[ \xymatrix @=-8pt @! {
{\hf^x} & & {\hs} & & {\hs} & & {\hs} \\
& {\hf^{z_1}} & & {\hs} & & {\hs} & & {\hs} \\
{\tilde{\hf}} & & {\hf} & & {\hs} & & {\hs} \\
& {\hf^{z_2}} & & {\hs} & & {\hs} & & {\hs} \\
{\hf^y} & & {\hs} & & {\hs} & & {\hs} \\
} \]
The same argument with the roles of $s_1$ and $s_{\tilde{1}}$ reversed
shows that there is an $s_1$ between every minimal pair of
$s_{\tilde{1}}$ generators.

Part (2):  Suppose $s_1$ and $s_{\tilde{1}}$ lie at the same level $k$
of the heap.  If column 1 contains another entry $p$, suppose
without loss of generality that $p$ forms a minimal pair with
the $s_1$ entry from level $k$ of the heap, and there are no other
entries in column 1 between $p$ and the $s_1$ entry from level $k$.
Then, in order for the heap to correspond with a reduced expression,
this pair must have a right resolution using a minimal pair of
$s_{2}$'s.  But resolution of the minimal pair of $s_2$ generators
includes an entry from column 1 by Lemma~\ref{sbadh.1}, contradicting
the minimality of our choice of $p$ from column 1.  Thus, there cannot
be other generators in the first column.  

Part (3) follows directly from Parts (1) and (2).
\end{proof}

\begin{lemma}\label{sbadh.3}
Let $w \in D_{n}$ be a short braid avoiding element.  If there exists
a minimal pair of entries $x,y$ in column $i\geq 2$ with only a left
resolution, then the part of the heap to the left of column $i$ has a
particular form shown below with exactly two entries in each of
columns $1,\dotsc, i$:

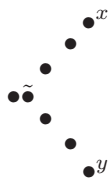
\[ \xymatrix @=-8pt @! {
& {\hs} & & {\hf^x} & & {\hs} & & {\hs} \\
{\hs} & & {\hf} & & {\hs} & & {\hs} \\
& {\hf} & & {\hs} & & {\hs} & & {\hs} \\
{\hf \tilde{\hf}} & & {\hs} & & {\hs} & & {\hs} \\
& {\hf} & & {\hs} & & {\hs} & & {\hs} \\
{\hs} & & {\hf} & & {\hs} & & {\hs} \\
& {\hs} & & {\hf^y} & & {\hs} & & {\hs} \\
} \]
\end{lemma}
\begin{proof}
Suppose $x$ and $y$ are separated by a minimal pair of generators
$z_1, z_2$ both lying in column $i-1$.  Applying reasoning similar
to that in Lemma~\ref{sbadh.1}, $z_1, z_2$ must have a left resolution
as well by minimality.  Continuing on, we eventually obtain a minimal
pair in column 2.  In contrast to the case of a right resolution, we
can resolve a minimal pair in column 2 to the left with $s_1
s_{\tilde{1}}$, and by minimality we must do so.  Since these
generators commute, they lie at the same level in the heap.

Observe that by minimality, there can be no entries between $x$ and
$y$, nor between any of the other minimal pairs in the left
resolution.  Note also that there cannot be any other entries in
column 1, by Lemma~\ref{sbadh.2}.  Therefore, if there existed any
other entries in columns $2, \ldots, i$ above or below the minimal
pairs, then they would create new minimal pairs with the existing
entries and would require a left resolution or a distinct resolution
by Lemma~\ref{sbadh.1}.  But resolving these implied minimal pairs
eventually requires additional entries in column 1, which contradicts
Lemma~\ref{sbadh.2}(2).  Thus, there can be no other elements in
columns $2,\ldots, i$.  
\end{proof}

\begin{definition}\label{d:critical_generator}
Observe that if position $k$ is a defect of $\w = \w_1 \dots \w_p$
with respect to the mask $\s$, then $\w^{\s[k-1]} \w_k$ is not
reduced.  Hence by Tits' theorem, there must be some entry $\w_j$ with
mask-value 1 that lies to the left of $\w_k$ in the reduced
expression and corresponds to the same generator as $\w_k$.  We call
the rightmost such position the \em critical generator \em for the
defect $k$.  
\end{definition}

The critical generator can be viewed as an element of the heap as
well.  The critical generator is always the first entry in the heap
below $k$ and in the same column as $k$.

\begin{lemma}\label{sbadh.4}
Suppose $w \in D_{n}$ is a connected, short braid avoiding,
non-Deodhar element and $Heap(w)$ contains a minimal pair of entries
with only a left resolution.  Then, we can construct an element
$\check{w}\in D_{k}$ with $k<n$ such that:
\begin{enumerate}
\item $w$ contains $\check{w}$ as a 1-line pattern,
\item every minimal pair of entries in $Heap(\check{w})$ has a
distinct resolution, and 
\item $\check{w}$ is a connected short-braid avoiding non-Deodhar element.  
\end{enumerate}
\end{lemma}
\begin{proof}
Suppose that the rightmost pair of entries that require a left
resolution lie in column $i$.  By Lemma~\ref{sbadh.3}, we have a
specific form for columns $1,\ldots, i$ in $Heap(w)$.  To construct
$Heap(\check{w})$ from $Heap(w)$, we remove columns 2 through $i$, and
shift columns $i+1, i+2, \dotsc $ to the left.  Hence, column 1 of
$Heap(\check{w})$ contains an $s_1s_{\tilde{1}}$ entry leftover from
the left resolution, and column $j\geq 2$ of $Heap(\check{w})$
contains the entries of column $j+i-1$ of $Heap(w)$.  In the example
below, $i=3$ and $Heap(\check{w})$ is obtained by removing the $4$
gray entries.

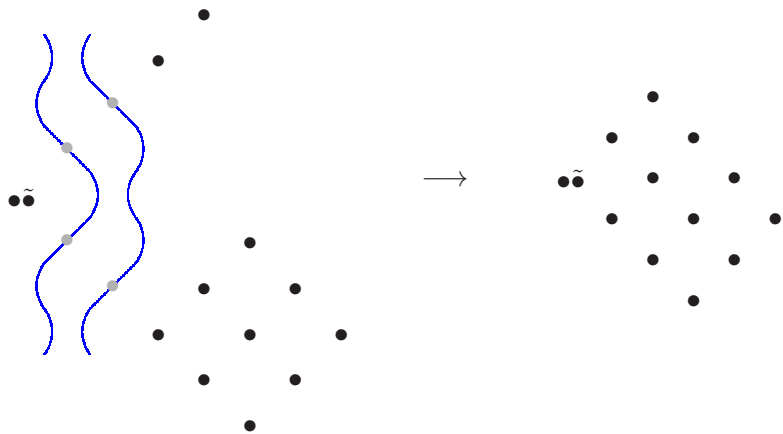
\[ \heap {
{\hs} & & {\hs} & & {\hf} & & {\hs} \\
& \StringLR{\hs} & & {\hf} & & {\hs} & & {\hs} \\
\StringR{\hs} & & \StringLX{\hb} & & {\hs} & & {\hs} \\
& \StringLX{\hb} & & \StringL{\hs} & & {\hs} & & {\hs} \\
{\hf \tilde{\hf}} & & \StringLR{\hs} & & {\hs} & & {\hs} \\
& \StringRX{\hb} & & \StringL{\hs} & & {\hf} & & {\hs} \\
\StringR{\hs} & & \StringRX{\hb} & & {\hf} & & {\hf} \\
& \StringLR{\hs} & & {\hf} & & {\hf} & & {\hf} \\
{\hs} & & {\hs} & & {\hf} & & {\hf} \\
& {\hs} & & {\hs} & & {\hf} & & {\hs} \\
} 
\parbox[t]{1in}{ \vspace{.8in} 
\hspace{.3in} $\longrightarrow$ }
\xymatrix @=-2pt @! {
{\hs} & & {\hs} & & {\hs} & & {\hs} \\
& {\hs} & & {\hs} & & {\hs} & & {\hs} \\
{\hs} & & {\hf} & & {\hs} & & {\hs} \\
& {\hf} & & {\hf} & & {\hs} & & {\hs} \\
{\hf \tilde{\hf}} & & {\hf} & & {\hf} & & {\hs} \\
& {\hf} & & {\hf} & & {\hf} & & {\hs} \\
{\hs} & & {\hf} & & {\hf} & & {\hs} \\
& {\hs} & & {\hf} & & {\hs} & & {\hs} \\
& {\hs} & & {\hs} & & {\hs} & & {\hs} \\
{\hs} & & {\hs} & & {\hs} & & {\hs} \\
& {\hs} & & {\hs} & & {\hs} & & {\hs} \\
} \]

By Lemma~\ref{sbadh.3}, the strings labeled $2,3,\ldots,i$ on the
bottom of the heap end up in the same positions at the top of the
heap.  Therefore, the strings of $Heap(\check{w})$ can be canonically
identified with strings $1,i+1,\ldots, n,$ of $Heap(w)$, so the
resulting string diagram for $\check{w}$ is well-defined with $w$
containing $\check{w}$ as a 1-line pattern.  Furthermore, since $w$ is
connected, so is $\check{w}$.

By the minimality assumption and Lemma~\ref{sbadh.1}, every minimal
pair of entries in $Heap(w)$ to the right of column $i$ has a
distinct resolution.  Therefore, by construction every pair in
$Heap(\check{w})$ has a distinct resolution.  In particular,
$\check{w}$ is short braid avoiding.

Finally, if $w$ is not Deodhar then we can find a proper non-Deodhar
mask $\s$ on $Heap(w)$.  We show this implies there exists a
non-Deodhar mask on $Heap(\check{w})$.  Recall from
\eqref{e:deodhar.ineq} that a mask is non-Deodhar whenever
\[ \text{\# zero-defects} \ge \text{\# plain-zeros} \]
which we refer to as the \em non-Deodhar bound \em throughout the
proof.

Consider the effect on the non-Deodhar bound of modifying
every mask value in $\sigma$ to be 1 in columns $2,\ldots, i$.  If we
set a plain-zero that is not involved with any defect to have
mask-value 1 then the mask remains non-Deodhar.

Say there exists a zero-defect at the top of column $h \leq i$.  Then
the strings for the defect must cross by
Corollary~\ref{c:string.crossing}.  By Lemma~\ref{sbadh.3}, the form
of the heap in columns $1, \dots, i$ is determined.  Hence, the left
string of the defect must travel southwest from the defect until it
hits a zero in column $1\leq g < h$, drop straight down until it hits
the next entry in the heap which must also be a zero, and then
continue southeast until it crosses the right string of the defect at
the bottom entry in column $h$.  Both of the entries in column $g$
must have mask-value 0 to facilitate the string crossing for the
defect.  Neither entry in column $g$ can be a zero-defect, because
every defect must have a critical generator below it in the same column.
As we already assumed that the mask values of both entries in column
$g$ were 0, the lower entry is not a critical generator.  Thus, we
find that setting the entries in columns $g$ and $h$ to have
mask-value 1 removes a zero-defect and two plain-zeros, which
preserves the non-Deodhar bound for the mask.

If there is no zero-defect in any column $j>i$ that has a
string passing through column $i$, then the mask $\sigma'$ obtained
from $\sigma$ by setting the mask-values of all entries in columns $1,
\dots, i$ to 1 will remain non-Deodhar, since the defect status
depends only on the string dynamics for the left and right strings of
an entry by Lemma~\ref{l:defectlabellings}.

On the other hand, if there exists a zero-defect in column $j>i$ whose left
string encounters column $i$, then in order for the strings of the
defect to cross, the path of the left string must follow a similar
course as described above.  That is, the left string travels southwest
to a plain-zero, say in column $f \leq i$, drops straight down to another
plain-zero in the same column, and continues southeast through column
$i$ again and beyond.  In this case, the mask $\sigma'$ is obtained
from $\sigma$ by setting the entries in column $1$ to have mask-value
0, and setting the mask-values of entries in columns $2, \dotsc, i$ to
1.  This will maintain the string dynamics for the left string of the
defect in column $j$, and effectively moves the plain zeros from
column $f$ to column 1.  Hence, the non-Deodhar bound is preserved.

In either of the definitions for $\sigma'$ above, the strings
$2,\ldots, i$ cross once to the right then once to the left so remain
in their original position.  In particular, all the plain-zeros and
zero-defects in $\sigma'$ applied to $w$ remain if we remove the
strings $2,\dotsc, i$.  Therefore, $\sigma'$ restricted to columns
$1,i+1,\dotsb, n$ determines a non-Deodhar mask for $\check{w}$.
\end{proof}


\bigskip
\section{Convex elements in type $D$}\label{s:convex} 

In this section, we restrict our attention to a subset of the
short braid avoiding elements that have the lateral convexity
property introduced for type $A$ in \cite{b-w}.

\begin{definition}
If $w \in D_{n}$ is short braid avoiding and every minimal pair of
entries in $Heap(w)$ has a distinct resolution then we say $w$ is
\em convex\em.
\end{definition}

\begin{remark}\label{r:oneright}
It follows immediately from the definition that in a convex type $D$
heap, there can only be a single generator in the rightmost column.
\end{remark}

For example, a permutation is short braid avoiding if and only if it
is convex.  The element $s_{2}s_{1}s_{\tilde{1}}s_{2} \in D_{4}$ is
not convex since it does not have a distinct resolution of the
$s_{2}$'s, while $s_{2}s_{1}s_{3} s_{\tilde{1}}s_{2} \in D_{4}$ is
convex.

It follows from \cite{simion-schmidt} that the number of short braid
avoiding elements in type $A_{n}$ is the Catalan number
$c_{n}=\frac{1}{n+1}\chs{2n}{n}$.  In types $B_{n}$ and $C_{n}$, the
short braid avoiding (equivalently, convex) elements are also counted
by Catalan numbers.  In type $D$, Fan and Stembridge have given an
explicit formula for the number of short braid avoiding elements
\cite{s2}.  Furthermore, Stembridge characterized the short braid
avoiding elements of type $D$ in terms of 1-line patterns.

\begin{theorem}\cite{s4}\label{t:stem}
An element $w$ in type $D$ is short braid avoiding if and only if $w$
avoids all 1-line patterns $[abc]$ where $|a| > b > c \ $ or $\ \bar{b} >
|a| > c$.  
\end{theorem}


In type $D$, there is a new sequence corresponding to the number of
convex elements in $D_{n}$ for $n\geq 1$.  The first 10 terms are
\[
1, 4, 13, 44, 154, 552, 2013, 7436, 27742, 104312
\]

\begin{remark}\label{r:convex.patterns}
This sequence is also characterized by 1-line patterns.  The patterns
are simply the patterns from Theorem~\ref{t:stem} and the single
additional pattern $[\bar{1}2\bar{3}]\in D_{3}$.  It would be
interesting to know more about this sequence.  The notion of a convex
element can be extended to other Coxeter groups where the Coxeter
graph can be linearized in a meaningful way.   
\end{remark}

\begin{definition}\label{d:coalescing}
Given a convex element $w \in D_{n}$, we define an operation called
\textit{coalescing}, which connects the lattice point components of
the heap as in \cite{b-w}.  Starting in column 1, if there exists an
entry in column 1 below an entry in column 2 with empty lattice points
between them, then allow the entry in column 1 to rise up until it is
blocked by the entry in column 2.  Then, allow all the entries below
the column 1 entry to rise up as well, until blocked by a
non-commuting generator.  Work to the right, continuing to apply the
same elevations until the heap is pushed together as much as possible.
The resulting collection of lattice points is called the
\textit{coalesced heap of} $w$.  
\end{definition}

Throughout the rest of the paper, whenever we refer to a
\textit{position}, we will mean the coordinates $(a,b)$ in the lattice
$\Z^2$ containing a coalesced heap for a convex element.  Here
$a$ is the column number and $b$ is the height of its row.  Each
position can be empty, contain a single generator of the heap, or
contain the two generators $s_1 s_{\tilde{1}}$ in the case when the
entry lies in the first column of a type $D$ heap.

Define two subsets of the heap as follows:
\begin{align*}
  \text{ The \textit{lower cone} of $(a,b)$: }\bcone(a,b) &= 
  \{(a, b) + i (-1,-1) + j (2,0)\in \mathbb{Z}^2: i \in [0, \infty), j
  \in [0, i] \}.\\    
  \text{ The \textit{upper cone} of $(a,b)$: }\ucone(a,b) &= 
  \{(a,b) + i (-1,1) + j(2,0) \in \mathbb{Z}^2: i \in [0, \infty), j
  \in [0,i] \}.    
\end{align*}

After coalescing the heap, the points in a heap that lie in
$\bcone(a,b)$ are precisely the elements of the heap poset that are
smaller than $(a,b)$.  Similarly, the points in a heap lying in
$\ucone(a,b)$ are the elements in the heap poset that are greater than
$(a,b)$.

\begin{definition}\label{d:latteral.convexity}
A collection of lattice points $L$ is \textit{laterally convex}
if $(a,b),(c,d) \in L$ implies $\bcone(a,b) \cap \ucone(c,d)
\subset L$.  
\end{definition}

\begin{lemma}[\textbf{Convexity Lemma}]\label{l:latterally.convex}
The coalesced heap of a connected convex element in type $A, B$ or $D$ is
laterally convex.
\end{lemma}
\begin{proof}
Every minimal pair in a convex element is resolved by two distinct
elements $x,y$.  If these elements appear on different levels in the
coalesced heap then there must be a chain of elements preventing the
lower element, say $x$, from rising without raising $y$.  But, this
would imply there is an alcove of the form
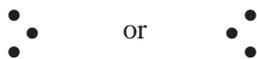
\begin{center}
\begin{tabular}{ccc}
	\stringlessheap {
	& \hf & \hs & \hs & \hs & \hs \\
	\hs & \hs & \hf & \hs & \hs \\
	& \hf & \hs & \hs & \hs & \hs \\
	} &
\parbox[t]{0.35in}{ \vspace{.05in} or } &
	\stringlessheap {
	& \hf & \hs & \hs & \hs & \hs \\
	\hf & \hs & \hs & \hs & \hs \\
	& \hf & \hs & \hs & \hs & \hs \\
	} 
\end{tabular}
\end{center}
consisting of a minimal pair without a distinct resolution,
contradicting that the heap is convex.  Therefore, locally every
minimal pair is only two rows apart $(a,b), (a,b+2)$, and its distinct
resolution is on adjacent points $(a-1,b+1), (a+1,b+1)$, given
pictorially as
\begin{equation}\label{e:four.corners}
\begin{array}{lll}
	\stringlessheap {
	\hs & \hs & \hs & \hs & \hs \\
	& \hf & \hs & \hs & \hs & \hs \\
	\hf & \hs & \hf & \hs & \hs \\
	& \hf & \hs & \hs & \hs & \hs \\
	}
\end{array}
\end{equation}

Suppose the positions $(a,b),(c,d)$ contain entries of the coalesced
heap, and $(a,b) \in \ucone(c,d)$ and $(c,d)\in \bcone(a,b)$.  Since
the heap is connected, there must exist a chain of adjacent
non-commuting generators connecting the two points in the coalesced
heap.  Every minimal pair along the path must have a distinct
resolution which by the argument above looks locally like
\eqref{e:four.corners}.  Filling out the diamond for each minimal pair
recursively implies that every point in $\bcone(a,b) \cap \ucone(c,d)$
is contained in the heap.
\end{proof}

\begin{remark}\label{r:convexity}
In the classification of Section~\ref{s:classification.d}, we will
consider only connected convex elements $w$, and we will assume that
they have been embedded in the lattice so that $Heap(w)$ is coalesced.
Then, these lattice points satisfy lateral convexity by
Lemma~\ref{l:latterally.convex}.  Hence, there is a \em minimal
resolution \em of any type $D$ heap fragment that arises in this
fashion, as a subset of lattice points of some $Heap(w)$.  In
particular, we can add lattice points to the heap fragment to resolve
minimal pairs of entries over columns $2, \dots, n$ distinctly as in
the classical lateral convexity of \cite{b-w}, and we may resolve
minimal pairs of entries in the first column as prescribed by
Lemma~\ref{sbadh.2}.  We frequently invoke the Convexity Lemma to
speak of \em the resolution by convexity \em of some collection of
entries inside a connected convex non-Deodhar type $D$ heap, and
deduce the existence of heap entries that were not explicitly given.
\end{remark}

\begin{lemma}\label{l:heap.cover}
Let $w$ be a connected convex element of type $A$, $B$, or $D$ with
coalesced $Heap(w)$.  Then, the point $y$ covers $x$ in the heap
poset for $w$ if and only if $x$ has heap coordinates $(i,j)$ and $y$
has heap coordinates $(i \pm 1, j+1)$.
\end{lemma}
\begin{proof}
Suppose $x$ has heap coordinates $(i,j)$ and $y$ has heap coordinates
$(i \pm 1, j+1)$.  Then, $x$ and $y$ do not commute, and since the
level difference is 1, we have that $y$ covers $x$ in the heap poset.

On the other hand, if $x$ has heap coordinates $(i,j)$ but $y$ has
heap coordinates $(i \pm \delta_1, j + \delta_2)$ for $\delta_2 \geq
\delta_1 \geq 1$ with $\delta_1 + \delta_2 > 2$, then all of the
points in $C = \bcone(i \pm \delta_1, j + \delta_2) \cap
\ucone(i,j)$ are in $Heap(w)$ by the lateral convexity proved in
Lemma~\ref{l:latterally.convex}.  In particular, $z = (i \pm 1, j+1)$
is in $C \subset Heap(w)$, so we have $x < z < y$ in the heap poset,
hence $y$ does not cover $x$ in the heap poset.

If $y$ has heap coordinates $(i \pm \delta_1, j + \delta_2)$ for $\delta_1
> \delta_2$, then $x$ and $y$ are unrelated in the heap poset.
\end{proof}

Let $M_{(x,y)}^{w}$ be the number of entries located at position
$(x,y)$ in the coalesced heap for $w$. These occupation numbers can
take the values $1$ or $2$ in type $D$, but the value $2$ can only
appear in column 1 and must obey the rules in Lemma~\ref{sbadh.2}.

\begin{definition}\label{d:saturated.subset}
Let $w$ and $p$ be connected convex elements of type $A$, $B$ or $D$
with coalesced heaps $Heap(w)$ and $Heap(p)$, respectively.  We say
that $Heap(w)$ contains $Heap(p)$ as a \em saturated subset \em of
lattice points, if there exist offsets $i, j \in \Z$ and a Coxeter
embedding $f$ such that for all occupied positions $(x,y)$ in
$Heap(f(p))$, we have
\[ M^{w}_{(x, y) + i(2, 0) + j(1, 1)} = M^{p}_{(x,y)}  \]
\end{definition}

\begin{example}
Let $p = s_2 s_1 s_3 s_2$ and $w = s_2 s_1 s_{\tilde{1}} s_3 s_2$.
Note that 
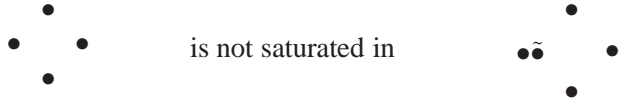
\[
\xymatrix @=1pt @! {
& {\hf} & & {\hs} \\
{\hf} & & {\hf} \\
& {\hf} & & {\hs} \\
} 
\parbox[t]{2in}{ \vspace{.1in} 
\hspace{.3in} is not saturated in }
\xymatrix @=-2pt @! {
& {\hf} & & {\hs} \\
{\hf \tilde{\hf}} & & {\hf} \\
& {\hf} & & {\hs} \\
}
\]
because $M^{p}_{(1,1)} = 1$ while $M^{w}_{(1,1)} = 2$.

On the other hand, $Heap(w)$ is saturated in the heap of any $w'$
containing $s_2 s_1 s_{\tilde{1}} s_3 s_2$ as a factor.
\end{example}

\begin{lemma}\label{l:ef.heap}
Let $w$ and $p$ be connected convex elements of type $A$, $B$, or $D$
with coalesced heaps $Heap(w)$ and $Heap(p)$, respectively.  Then,
the following are equivalent:
\begin{enumerate}
\item $w$ contains $p$ as an embedded factor.
\item $Heap(p)$ is contained in $Heap(w)$ as a saturated subset of
	lattice points.
\end{enumerate}
\end{lemma}
\begin{proof}
$(1) \implies (2)$.  Recall from \cite{ec1} that a subposet $Q$ of
some poset $R$ is called \em convex \em if $b \in Q$ whenever $a < b <
c$ in $R$ and $a, c \in Q$.  Take a reduced expression for $w$ of the
form $\w = \u f(\p) \v$ where $f$ is a Coxeter embedding for $\p$.  By
Definition~\ref{d:heap.poset}, the labeled heap poset of $\w = \u
f(\p) \v$ contains the labeled heap poset of $f(\p)$ as a convex
subposet.  Lemma~\ref{l:heap.cover} shows that the Hasse diagram for
the heap poset of $\w$ can be embedded into $[n] \times \N$ where $(i,
j) \lessdot (i \pm 1, j+1)$ is the corresponding cover relation, and 
this embedding is $Heap(w)$.  Hence, the convex subposet
corresponding to the labeled heap poset of $f(\p)$ appears in
$Heap(w)$ as a saturated subset of lattice points with shape
$Heap(f(\p))$.

$(2) \implies (1)$.  If $Heap(f(p))$ is contained as a saturated
subset of the lattice points in $Heap(w)$ for some Coxeter embedding
$f$, then we can build $Heap(w)$ from $Heap(f(p))$ by sequentially
adding lattice points.  Moreover we can only add points that are
maximal or minimal entries of the intermediate heap since
$Heap(f(p))$ is saturated in $Heap(w)$.  This operation is equivalent
to multiplying $f(p)$ on the left or right by the respective
generators.  Hence, $w$ contains $f(p)$ as a factor.
\end{proof}


We now generalize the notion from \cite{b-w} of a \em right critical
zero \em associated to a defect for convex heaps in type $D$.

\begin{definition}
The first zero encountered by the right string of a defect along the
southeast diagonal containing the defect is the \textit{right critical
zero}.
\end{definition}

For example, the defect $d$ has right critical zero $r$:
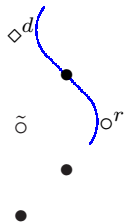
\[ \heap {
	\StringR{\hd^d} & & {\hs} & & {\hs} & \\
	& \StringLX{\hf} & & {\hs} & & {\hs} \\
	{\tilde{\hz}} & & \StringL{\hz^r} & & {\hs} & \\
	& {\hf} & & {\hs} & & {\hs} \\
	{\hf} & & {\hs} & & {\hs} & \\
	& {\hs} & & {\hs} & & {\hs} \\
} \]

\begin{lemma}[\textbf{Right Critical Zero Lemma}] \label{l:rcz}
Suppose $w\in D_{n}$ is convex, and $\s$ is a mask on $w$.  Then every
defect in $Heap(w)$ has a right critical zero.
\end{lemma}
\begin{proof}
Suppose the right string for a defect never encounters an entry with
mask-value 0 as it travels southeast.  It eventually turns southwest
when it hits the first entry $(a,b)$ not in the heap along that
diagonal.  As the string travels down, it cannot encounter another
heap element in column $a-1$ before encountering one in column $a$ by
convexity.  If an element is encountered in position $(a, b-i)$ as the
string travels down, then by lateral convexity $(a,b)$ must exist in
the heap contradicting our hypothesis.  If no element is
encountered in position $(a, b-i)$ as the string travels down, then it
cannot cross the left string below, contradicting the fact that the
original position was a defect by Corollary~\ref{c:string.crossing}.
\end{proof}

Note that defects in type $D$ need not have a corresponding left
critical zero, since a string leaving the heap on the left can
re-enter the heap without a mask-value 0 entry.  See
Example~\ref{ex:defect.status}.

\bigskip
\section{Classification of convex type $D$ patterns}\label{s:classification.d}

In this section, we show that among the set of convex elements of type
$D$, there are precisely six that are minimally non-Deodhar under the
partial order of embedded factor containment.  Recall that the convex
elements are short braid avoiding, hence they have unique heaps.
These minimal heaps have the property that if we remove either an
entry that is minimal or maximal with respect to the heap poset
structure, then the resulting heap corresponds to a Deodhar element.
Furthermore, we will always consider the heaps to be coalesced as
described in Definition~\ref{d:coalescing}.  

\begin{theorem}\label{t:minimal.patterns}
Below is the complete list of the heaps of convex minimal non-Deodhar
embedded factor patterns in type $D$ up to the Coxeter graph
isomorphism that interchanges $s_{1}$ and $s_{\tilde{1}}$:
\end{theorem}

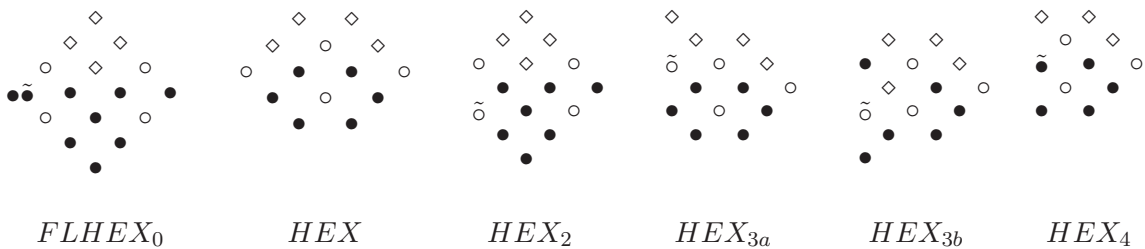
\begin{equation}\label{e:minimal.patterns}
\begin{matrix}
        \xymatrix @=-8pt @! {
{\hs} & & {\hs} & & {\hs} & & {\hs} \\
& {\hs} & & {\hd} & & {\hs} &\\
{\hs} & & {\hd} & & {\hd} & & {\hs} \\
& {\hz} & & {\hd} & & {\hz}\\
{\hf \tilde{\hf} }&  & {\hf} & & {\hf} & & {\hf} & {\hs} \\
& {\hz} & & {\hf} & & {\hz} &\\
{\hs} & & {\hf} & & {\hf} & & {\hs} \\
& {\hs} & & {\hf} & & {\hs} &\\
    } &
        \xymatrix @=-2pt @! {
& {\hs} & & {\hs} & & {\hs} \\
{\hs} & & {\hd} & & {\hd} & & {\hs} \\
& {\hd} & & {\hz} & & {\hd}\\
{\hz} & & {\hf} & & {\hf} & & {\hz} \\
& {\hf} & & {\hz} & & {\hf} \\
{\hs} & & {\hf} & & {\hf} & & {\hs} \\
& {\hs} & & {\hs} & & {\hs} \\
    }  &
	\stringlessheap {
{\hs} & & {\hs} & & {\hs} & & {\hs} \\
& {\hs} & & {\hd} & & {\hs} \\
{\hs} & & {\hd} & & {\hd} & & {\hs} \\
& {\hz} & & {\hd} & & {\hz}\\
{\hs} & & {\hf} & & {\hf} & & {\hf} \\
& {\tilde{\hz}} & & {\hf} & & {\hz} \\
{\hs} & & {\hf} & & {\hf} & & {\hs} \\
& {\hs} & & {\hf} & & {\hs} \\
    } &
	\stringlessheap {
{\hs} & & {\hs} & & {\hs} & & {\hs} \\
& {\hd} & & {\hs} & & {\hs} \\
{\hs} & & {\hd} & & {\hd} & & {\hs} \\
& {\tilde{\hz}} & & {\hz} & & {\hd}\\
{\hs} & & {\hf} & & {\hf} & & {\hz} \\
& {\hf} & & {\hz} & & {\hf} \\
{\hs} & & {\hf} & & {\hf} & & {\hs} \\
& {\hs} & & {\hs} & & {\hs} \\
    } & 
	\stringlessheap {
{\hs} & & {\hs} & & {\hs} & & {\hs} \\
& {\hs} & & {\hs} & & {\hs} \\
{\hs} & & {\hd} & & {\hd} & & {\hs} \\
& {\hf} & & {\hz} & & {\hd}\\
{\hs} & & {\hd} & & {\hf} & & {\hz} \\
& {\tilde{\hz}} & & {\hz} & & {\hf} \\
{\hs} & & {\hf} & & {\hf} & & {\hs} \\
& {\hf} & & {\hs} & & {\hs} \\
    } & 
	\stringlessheap {
 & {\hs} & & {\hs} \\
 {\hd} & & {\hd} & & {\hs} \\
 & {\hz} & & {\hd}\\
 {\tilde{\hf}} & & {\hf} & & {\hz} \\
 & {\hz} & & {\hf} \\
 {\hf} & & {\hf} & & {\hs} \\
 & {\hs} & & {\hs} \\
    } \\
    & &  &  &  &  \\
    FLHEX_0 & HEX & HEX_2 & HEX_{3a} & HEX_{3b} & HEX_4 \\
\end{matrix}
\end{equation}

Each heap above is decorated with a mask that demonstrates the
non-Deodhar condition using the notation described in
Section~\ref{s:d.heaps}.  It is straightforward to verify by computer
that these heaps are minimal among the convex non-Deodhar elements in
type $D$ using embedded factor containment.  Our goal is to show that
this list is complete.  Note that we could reduce the list further by
requiring both $w$ and $w^{-1}$ to avoid the pattern $HEX_{3a}$ since
$HEX_{3b}$ is its inverse and by Lemma~\ref{l:w.winverse}, $w$ is
non-Deodhar if and only if $w^{-1}$ is non-Deodhar.  All the other
heaps in \eqref{e:minimal.patterns} have $w=w^{-1}$.

Note that the pattern $HEX$ can embed on any 7 connected columns
corresponding to a path of length 7 in the Coxeter graph, including
$s_{\tilde{1}},s_{2},\dots , s_{7}$.  The other patterns are tied to
the initial segment of the type $D$ Coxeter graph, since they include
both of the generators $s_1, s_{\tilde{1}}$ from column 1.  These
patterns are the heaps for the reduced words that appear in the table
on Page~\pageref{f:patterns} in the type $D$ case.

\begin{remark}
Recall from Lemma~\ref{l:ef.heap} that factor containment is expressed
in the heap setting by deleting a sequence of entries that are minimal
or maximal entries of the intermediate heap at each step.  In
particular, $HEX$ is not a saturated subset of $FLHEX_{0}$ by
Definition~\ref{d:saturated.subset}.
\end{remark}

Let $D_n^{C}$ and $A_\infty^{C}$ denote the convex elements in type
$D_{n}$ and type $A$ on any finite number of generators, respectively.
We will prove Theorem~\ref{t:minimal.patterns} by constructing a pair
of projection maps
\begin{align*}
\pi_{NE}:& D_n^{C} \rightarrow A_\infty^{C}\\
\pi_{SE}:& D_n^{C} \rightarrow A_\infty^{C}
\end{align*}
that we can use to compare minimal patterns in type $D$ with the
hexagon in type $A$.  We will show that these maps
have have the following properties:
\smallskip
\begin{enumerate}
\item [a.] If $w \in D_{n}^{C}$ is non-Deodhar in type $D$, then at least one of the
two projections sends $w$ to a non-Deodhar element in type $A$.
\medskip

\item [b.] If $\pi_{NE}(w)$ or $\pi_{SE}(w)$ is non-Deodhar in type
$A$, then $w$ contains one of the patterns listed above as an embedded
factor.
\end{enumerate}
\smallskip Therefore, if $w$ is any minimal non-Deodhar element in
type $D$, then it appears on the list.  Hence, the list of convex
minimal non-Deodhar elements is complete.  In particular, all convex
minimal non-Deodhar patterns occur in $D_{8}$.

We will define the projections $\pi_{NE}(w)$, $\pi_{SE}(w)$ in terms
of three other maps.  First, define the map
\begin{equation*}
\f_i: D_n^{C} \rightarrow A_{n-1} \hspace{.1in} \text{for $i \in
\{1,\tilde{1} \}$ }
\end{equation*}
as the projection of $w$ onto a type $A$ heap obtained by removing all
the entries that correspond to the $\tilde{i}$ generator, where
$\tilde{i}$ is the generator $\{1,\tilde{1} \} \setminus \{i \}$.
That is, $\f_i$ assigns the $s_{i}$ generator from column 1 of the
type $D$ heap to $s_{1}$ in column 1 of the type $A$ heap, and removes
all the entries from column 1 that correspond to the other generator.
This generally leaves a heap fragment in type $A$ that is \em not \em
convex.

How might we construct an associated convex element in type $A$?  We
consider two constructions.   We could take an \textit{additive
approach}, denoted
\begin{equation*}
\mathbf{\afc}_i: D_{n}^{C} \rightarrow A_{\infty}^{C},
\end{equation*}
where we first project the element into type $A$ by the map $\f_i$,
then extend the Coxeter graph linearly to the left beyond the first
column and add entries as necessary in order to resolve all minimal
pairs of entries that do not have a distinct resolution.  By
Remark~\ref{r:convexity} restricted to type $A$, this operation is
unique.  For example, the composition shown below is
$\mathbf{\afc}_1$.
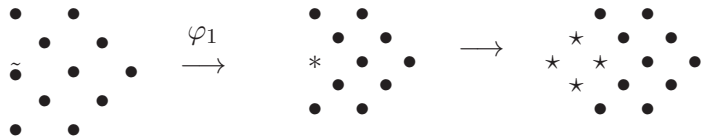
\[
\begin{matrix}
       \xymatrix @=-3pt @! {
{\hs} & {\hs} & {\hf} & {\hs} & {\hf} &  {\hs} & {\hs}\\
{\hs} & {\hs} & {\hs} & {\hf} & {\hs} &  {\hf} & {\hs}\\
{\hs} & {\hs} & {\tilde{\hf}} & {\hs} & {\hf} &  {\hs} & {\hf}\\
{\hs} & {\hs} & {\hs} & {\hf} & {\hs} &  {\hf} & {\hs}\\
{\hs} & {\hs} & {\hf} & {\hs} & {\hf} &  {\hs} & {\hs}\\
    }
\hspace{.1in}
\begin{array}{c}
\\ \\ 
\f_{1}\\
\longrightarrow  
\end{array}
\hspace{.1in}
       \xymatrix @=-3pt @! {
{\hs} & {\hs} & {\hf} & {\hs} & {\hf} &  {\hs} & {\hs}\\
{\hs} & {\hs} & {\hs} & {\hf} & {\hs} &  {\hf} & {\hs}\\
{\hs} & {\hs} & {\hv} & {\hs} & {\hf} &  {\hs} & {\hf}\\
{\hs} & {\hs} & {\hs} & {\hf} & {\hs} &  {\hf} & {\hs}\\
{\hs} & {\hs} & {\hf} & {\hs} & {\hf} &  {\hs} & {\hs}\\
    } 
\hspace{.1in}
\begin{array}{c}
\\
\\
\longrightarrow  
\end{array}
\hspace{.1in}
       \xymatrix @=-3pt @! {
{\hs} & {\hs} & {\hf} & {\hs} & {\hf} &  {\hs} & {\hs}\\
{\hs} & {\hh} & {\hs} & {\hf} & {\hs} &  {\hf} & {\hs}\\
{\hh} & {\hs} & {\hh} & {\hs} & {\hf} &  {\hs} & {\hf}\\
{\hs} & {\hh} & {\hs} & {\hf} & {\hs} &  {\hf} & {\hs}\\
{\hs} & {\hs} & {\hf} & {\hs} & {\hf} &  {\hs} & {\hs}\\
       }
\end{matrix}
\]
Here, $\hv$ denotes a point not in the heap, and $\hh$ denotes a
highlighted point in the heap.

The second construction is a \textit{subtractive approach} that
applies when there are at least two entries in the first column of the
type $D$ heap.  Define
\begin{align*}
&\mathbf{\sfc}_{SE}: D_n^C \rightarrow D_{n}^{C}\\
&\mathbf{\sfc}_{NE}: D_n^C \rightarrow D_{n}^{C}
\end{align*}
where $\mathbf{\sfc}_{NE}$ removes entries in the northeast (NE)
diagonal up to the first lattice point not in the heap.
Similarly, $\mathbf{\sfc}_{SE}$ removes the southeast (SE) diagonal.
This process is referred to as ``shaving off a diagonal'' in the
proofs that follow.  For example,
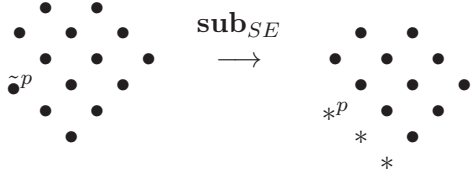
\[
\begin{matrix}
       \xymatrix @=-7pt @! {
 {\hs} & {\hf} & {\hs} &  {\hf} & {\hs} & {\hs}\\
 {\hf} & {\hs} & {\hf} & {\hs} &  {\hf} & {\hs}\\
 {\hs} & {\hf} & {\hs} &  {\hf} & {\hs} & {\hf}\\
 {\tilde{\hf}^p} & {\hs} & {\hf} & {\hs} &  {\hf} & {\hs}\\
 {\hs} & {\hf} & {\hs} &  {\hf} & {\hs} & {\hs}\\
 {\hs} & {\hs} & {\hf} & {\hs} &  {\hs} & {\hs}\\
   }
\hspace{.1in}
\begin{array}{c}
\\
\\
\mathbf{\sfc}_{SE}\\
\longrightarrow  
\end{array}
\hspace{.1in}
       \xymatrix @=-7pt @! {
 {\hs} & {\hs} & {\hs} & {\hs} &  {\hs} & {\hs}\\
 {\hs} & {\hf} & {\hs} &  {\hf} & {\hs} & {\hs}\\
 {\hf} & {\hs} & {\hf} & {\hs} &  {\hf} & {\hs}\\
 {\hs} & {\hf} & {\hs} &  {\hf} & {\hs} & {\hf}\\
 {{\hv}^p} & {\hs} & {\hf} & {\hs} &  {\hf} & {\hs}\\
 {\hs} & {\hv} & {\hs} &  {\hf} & {\hs} & {\hs}\\
 {\hs} & {\hs} & {\hv} & {\hs} &  {\hs} & {\hs}\\
       }
\end{matrix}
\]
Note that the result of applying $\mathbf{\sfc}_{NE}$ or
$\mathbf{\sfc}_{SE}$ is always again convex because we are removing
entries that are maximal or minimal with respect to the heap poset, so
the entries cannot be used in a distinct resolution.

We now define $\pi_{SE}(w)$ and $\pi_{NE}(w)$ in cases depending on
how many levels the coalesced $Heap(w)$ contains in the first column.
The rough idea that motivates the following technical definition is to
use the minimal patterns listed in \eqref{e:minimal.patterns} as a
guide, taking a subtractive approach when we have a heap fragment
that is Deodhar and taking an additive approach when we have a heap
fragment that is non-Deodhar.

\begin{definition}\label{d:pi}
Let $w \in D_{n}^{C}$.  Suppose that the first column of $Heap(w)$ has
entries in $k(w)$ distinct levels.  Let $i \in \{1, \tilde{1} \}$ be
the generator corresponding to the top entry in column 1.  Let $dir$
denote the direction SE or NE.  Then, we define $\pi_{dir}: D_n^{C}
\rightarrow A_\infty^{C}$ as follows:
$$
\pi_{dir}(w)=
\begin{cases}
w & \text{if $k(w)=0$.} \\
\mathbf{\afc}_{i}(w)  &   \text{if $k(w)=1$ and there is a single entry in column 1.} \\
\mathbf{\afc}_{1}(w)  &   \text{if $k(w) = 1$ and there are two entries in column 1, } \\
   & \text{\hspace{.2in}with at least 4 entries in column 4.} \\
\mathbf{\sfc}_{dir}(w)  &  \text{if $k(w)=1$ and there are two entries in column 1, } \\ 
   & \text{\hspace{.2in}with at most 3 entries in column 4. } \\
\mathbf{\afc}_1(w)   & \text{if $k(w)=2$ and there are at least 4 entries in column 3.} \\
\mathbf{\sfc}_{dir}(w) & \text{if $k(w)=2$ and there are at most 3 entries in column 3.} \\
\mathbf{\afc}_{i}(w)  &  \text{if $k(w)=3$.}  \\
\pi_{dir} \circ \mathbf{\sfc}_{dir}(w)  & \text{if $k(w)=4$.} \\ 
\mathbf{\afc}_{1}(w) & \text{if $k(w)\geq 5$.} 
\end{cases}
$$

\medskip 

\noindent 
This definition relies on the following observations.
\begin{enumerate}
\item When there are entries in more than one level of column 1, the
	entries of the first column must alternate by
	Lemma~\ref{sbadh.2}.  

\item In the case $k(w) = 1$, if there two entries in column 1 with
	less than 4 entries in column 4, then column 1 contains both
	$s_1$ and $s_{\tilde{1}}$ on the same level in column 1.  We
	define the map $\pi_{dir}$ to remove both of these generators
	in column 1, and shave any remaining entries from the chosen
	diagonal.

\item In the case $k(w) = 4$, observe that $\pi_{dir}$ is defined by
	shaving the top or bottom entry of column 1, and then reducing
	to the case when $k(w) = 3$.
\end{enumerate}
\end{definition}

Note that the map $\pi_{dir}$ always preserves the occupied lattice
points that are present in the type $D$ heap, with the possible
exception that they may remove the northeast diagonal of maximal
entries or southeast diagonal of minimal entries, starting from the
first column and working to the right.

The $\pi_{dir}$ projections satisfy a useful property with respect to
taking inverses that we will exploit in the classification.

\begin{lemma}\label{l:inverse.equivariant}
For all $w \in D_n^C$, we have
\[ \pi_{SE}(w^{-1}) = \pi_{NE}(w)^{-1} \]
\end{lemma}
\begin{proof}
Observe that the heap of $w^{-1}$ is obtained from the heap of $w$ by
flipping the diagram upside down.  In particular, the number of entries
in each column of $w$ and $w^{-1}$ is the same.  Therefore, $w^{-1}$
falls into the same case as $w$ when computing $\pi_{dir}$. Flipping the heap has
the effect of switching all the SE and NE diagonals, while the
additive resolutions remain the same.
\end{proof}

Consider the following partial heaps drawn with black dots called the $I$-\textit{shape} and the 4-\textit{stack}:
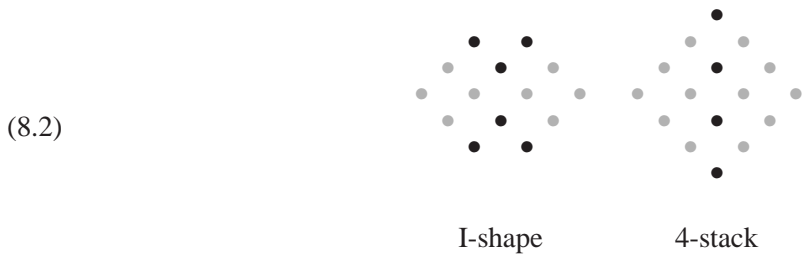
\begin{equation}\label{e:shapes}
\begin{matrix}
        \xymatrix @=-2pt @! {
& {\hs} & & {\hs} & & {\hs} \\
{\hs} & & {\hf} & & {\hf} & & {\hs} \\
& {\hb} & & {\hf} & & {\hb}\\
{\hb} & & {\hb} & & {\hb} & & {\hb} \\
& {\hb} & & {\hf} & & {\hb} \\
{\hs} & & {\hf} & & {\hf} & & {\hs} \\
& {\hs} & & {\hs} & & {\hs} \\
    }  &
         \xymatrix @=-2pt @! {
& {\hs} & & {\hf} & & {\hs} \\
{\hs} & & {\hb} & & {\hb} & & {\hs} \\
& {\hb} & & {\hf} & & {\hb}\\
{\hb} & & {\hb} & & {\hb} & & {\hb} \\
& {\hb} & & {\hf} & & {\hb} \\
{\hs} & & {\hb} & & {\hb} & & {\hs} \\
& {\hs} & & {\hf} & & {\hs} \\
    }  \\
    & \\
    \text{I-shape} & \text{4-stack} \\
\end{matrix}
\end{equation}
Some of the shaded dots will also appear in any convex heap containing
either the $I$-\textit{shape} or the 4-\textit{stack}.  However, if the
shape appears in the first 4 columns then some of the shaded dots on
the left hand side of the picture will get trucated in the heap.  Note that the
4-stack consists of 4 adjacent copies of the same generator if it
appears in column 1.  Therefore, a ``4-stack'' in column 1 of a convex
heap necessarily contains entries from 7 distinct levels.  The
following result is used frequently in the classification.

\begin{lemma}{\bf \em (Shape Lemma) \em}\label{l:shape}
Suppose $w \in D_{n}^{C}$, $Heap(w)$ contains an I-shape or a
4-stack, and the first column of $Heap(w)$ has entries from at least
two distinct levels.  Then, $w$ is non-Deodhar and there is a choice
of projection $\pi_{NE}(w)$ or $\pi_{SE}(w)$ that contains a hexagon,
and so preserves the non-Deodhar condition.
\end{lemma}
\begin{proof}
Consider the convex resolution of either the I-shape or the 4-stack.
By Remark~\ref{r:convexity}, the points in the convex resolution must
appear in $Heap(w)$.  This implies that all of the gray dots shown in
\eqref{e:shapes} to the right of the shape and a subset of the gray
dots to the left of the shape must be in the heap, depending on how far
the shape lies from the first column.  By considering the columns
where the shape lies, one can verify that no matter where the I-shape
or 4-stack appear, the convex resolution must contain one of the
patterns in \eqref{e:minimal.patterns} as a set of points in the
lattice.  Since $Heap(w)$ contains entries from at least two distinct
levels, we have that $Heap(w)$ contains one of the minimal heaps in
\eqref{e:minimal.patterns} as a saturated set of lattice points by
Lemma~\ref{sbadh.2}, hence $w$ also contains it as an embedded factor
by Lemma~\ref{l:ef.heap}.  Therefore, $w$ is non-Deodhar.

Furthermore, it can be verified that Definition~\ref{d:pi} has been
chosen precisely so that if $Heap(w)$ contains one of the heaps in
$\{ HEX, HEX_2, HEX_{3a}, HEX_{3b}, HEX_4 \}$ as an embedded factor,
then there is always a choice of direction NE or SE so that the
projection $\pi_{dir}$ contains a hexagon in type $A$.

To carry out the verification, suppose that $Heap(w)$ contains $HEX$
using column 1.  Then, if there are two entries in column 1, we either
project additively (if there is also a 4-stack in column 3), or we
choose to shave the entry of column 1 that is not being used in the
hexagon, and its diagonal.  If there are 3 entries in column 1, then
we always project additively, so the hexagon is preserved.  If there
are 4 entries in column 1, then we shave one of the extremal entries
that is not used in the hexagon, and project additively.  If there
are more than 4 entries in column 1, then we project additively so the
hexagon is preserved.  If the first column of the $HEX$ factor occurs
in column 2 of $Heap(w)$, or further to the right, and we are taking
a subtractive resolution in $\pi_{dir}(w)$, then there is always a
diagonal to shave that does not intersect the hexagon.  Otherwise, we
take an additive resolution, so the hexagon is preserved.

The shapes $\{ HEX_2, HEX_{3a}, HEX_{3b}, HEX_4 \}$ must all occur in
column 1, so it suffices to check that there is a choice of projection
that preserves the shape, regardless of how many entries column 1
contains.  Note that column 1 must contain at least two entries by
hypothesis, and if it contains more than 4 entries, then we resolve
additively, so the shapes are preserved.  The cases left to verify are
then that there exists a non-Deodhar projection $\pi_{NE}(w)$ or
$\pi_{SE}(w)$ if $Heap(w)$ contains: 
\begin{itemize}
	\item $HEX_2$ with 2, 3 or 4 entries in column 1, 
	\item $HEX_{3a}$ with 3 or 4 entries in column 1, 
	\item $HEX_{3b}$ with 3 or 4 entries in column 1, or 
	\item $HEX_4$ with 3 or 4 entries in column 1.
\end{itemize}
This verification is finite and straightforward, so we omit it.
\end{proof}

We now prove Property (b) for $\pi_{dir}$.

\begin{proposition}\label{t:part.a}
Let $w \in D_{n}^{C}$.  If $\pi_{dir}(w) \in A_{\infty}^{C}$ is
non-Deodhar, then $w$ contains a pattern from
\eqref{e:minimal.patterns} as an embedded factor.
\end{proposition}
\begin{proof}
Suppose $\pi_{dir}(w)$ is non-Deodhar.  Then, by Theorem
\ref{r:typeahex} the heap of $\pi_{dir}(w)$ contains a hexagon.  We
begin by considering the various ways in which the first column of the
type $D$ heap can be shifted under the $\pi_{dir}$ map.  Let $\s(w)$
be the number of the column that contains the image of column 1 under
$\pi_{dir}(w)$.  This can be computed by determining how many columns
must be added to the left when applying the $\mathbf{\afc}$ map 
or subtracted if the entire first column is removed.  

\begin{center}
\begin{tabular}{|p{4in}|p{0.8in}|} 
\hline
Case for $\pi_{dir}(w)$ & Shift $\s(w)$ \\
\hline
$k(w)=0$ & 1 \\
$k(w)=1$ and there is a single entry in column 1 & 1 \\
$k(w)=1$ and there are two entries in column 1, with at least 4 entries in column 4 & 1 \\
$k(w)=1$ and there are two entries in column 1, with less than 4 entries in column 4 & 0 \\
$k(w)=2$ and there are at least 4 entries in column 3 & 2 \\
$k(w)=2$ and there are less than 4 entries in column 3 & 1 \\
$k(w)=3$ & 3 \\
$k(w)=4$ & 3 \\ 
$k(w)\geq 5$ & $k(w)$ \\
\hline
\end{tabular}
\end{center}

If the hexagon appears weakly to the right of column $\s(w)$ then
$Heap(w)$ contains the shape $HEX$ in \eqref{e:minimal.patterns},
since $\pi_{dir}$ only adds points to the left of column $\s(w)$.  In
particular, if $\s(w) = 0$, then from the table above we see $k(w)=1$
and column 1 contains $\{ s_1, s_{\tilde{1}} \}$ on the same level,
with less than 4 entries in column 4.  In this case, $\pi_{dir}(w)$
removes all of the entries in the first column, and shaves a diagonal.
If there exists a hexagon in $\pi_{dir}(w)$, we have that there is a
hexagon in $Heap(w)$ that lies strictly to the right of column 1 in
$Heap(w)$, so $w$ contains $HEX$ as an embedded factor.

If $\s(w) = 1$ then from the table above we see that $k(w)=1$ or
$k(w)=2$ and there are less than 4 entries in column 3.  In each of
these cases the $\pi_{dir}$ projection doesn't add any points to the heap.
Therefore, if $\pi_{dir}(w)$ contains a hexagon that uses column 1, then
$Heap(w)$ contains either $HEX$ or $FLHEX_{0}$ as an embedded factor.

If $\s(w) = 2$ then $k(w)=2$ and there is a 4 stack in column 3, so by
Lemma~\ref{sbadh.2} and Remark~\ref{r:convexity}, $Heap(w)$ contains $HEX_{2}$ as
an embedded factor.

If $\s(w) = 3$, then $k(w)=3$ or $k(w)=4$.  In the case $k(w)=3$,
$\pi_{dir}(w)$ adds three new points to the left of $Heap(w)$.  If a
hexagon in $\pi_{dir}(w)$ uses column 1, then $Heap(w)$ must contain
$HEX_{4}$.  
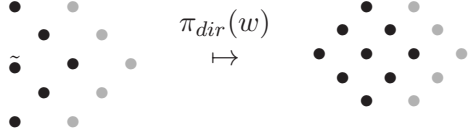
\[
\begin{matrix}
       \xymatrix @=-3pt @! {
{\hs} & {\hs} & {\hf} & {\hs} & {\hb} &  {\hs} & {\hs}\\
{\hs} & {\hs} & {\hs} & {\hf} & {\hs} &  {\hb} & {\hs}\\
{\hs} & {\hs} & {\tilde{\hf}} & {\hs} & {\hf} &  {\hs} & {\hb}\\
{\hs} & {\hs} & {\hs} & {\hf} & {\hs} &  {\hb} & {\hs}\\
{\hs} & {\hs} & {\hf} & {\hs} & {\hb} &  {\hs} & {\hs}\\
    }
\hspace{.1in}
\begin{array}{c}
\\ \\ 
\pi_{dir}(w)\\
\mapsto  
\end{array}
\hspace{.1in}
       \xymatrix @=-3pt @! {
{\hs} & {\hs} & {\hf} & {\hs} & {\hb} &  {\hs} & {\hs}\\
{\hs} & {\hf} & {\hs} & {\hf} & {\hs} &  {\hb} & {\hs}\\
{\hf} & {\hs} & {\hf} & {\hs} & {\hf} &  {\hs} & {\hb}\\
{\hs} & {\hf} & {\hs} & {\hf} & {\hs} &  {\hb} & {\hs}\\
{\hs} & {\hs} & {\hf} & {\hs} & {\hb} &  {\hs} & {\hs}\\
       }
\end{matrix}
\]
If a hexagon in $\pi_{dir}(w)$ uses column 2 (but not column 1) then
$Heap(w)$ contains $HEX_{3a}$ or $HEX_{3b}$ as an embedded factor.  If
a hexagon appears in column 3 or further to the right, then $Heap(w)$
contains $HEX$ as an embedded factor as well. 

In the case $\s(w)=3$ and $k(w)=4$ then $\pi_{dir}(w)$ shaves a maximal or
minimal diagonal and adds back three points as above.  The same
analysis of the placement of a hexagon in $\pi_{dir}(w)$ implies the
existence of a $HEX_{3a}$, $HEX_{3b}$ or a $HEX_{4}$ as embedded
factors.

If $\s(w) = s \geq 5$, we again must have at least 5 entries in
column 1 of $w$, so by convexity $w$ contains $HEX_4$.  Note that
$\s(w) = 4$ never occurs.
\end{proof}

We now turn to the proof of Property (a) for $\pi_{dir}$.  In preparation
for the proof, recall from Corollary~\ref{c:string.crossing} that the
strings for a defect must cross below the defect.  The first time the
two strings meet again will be called the \em initial string crossing
\em for the defect or just the \em string crossing\em, for short.

\begin{proposition}
Let $w \in D_{n}^{C}$.  If $w$ is non-Deodhar, then there is some
choice of projection $\pi_{NE}(w)$ or $\pi_{SE}(w)$ that is
non-Deodhar.
\end{proposition}
\begin{proof}
To prove this, we will consider the various cases given in
Definition~\ref{d:pi} for $\pi_{dir}$, breaking on the number of
levels in the first column of $w$ denoted $k(w)$.  Recall from
\eqref{e:deodhar.ineq} that a proper mask is non-Deodhar whenever
\[ \text{\# zero-defects} \ge \text{\# plain-zeros} \]
which we refer to as the \em non-Deodhar bound \em throughout the
proof.  The general strategy is to either show that $w$ contains a
4-stack or an I-shape from Lemma~\ref{l:shape} when $k(w) \geq 2$, or
to assume a proper non-Deodhar mask on $w$, and then show how to
adjust it to obtain a mask on $\pi_{dir}(w)$ that retains the
non-Deodhar bound by deleting at least as many plain-zeros as
zero-defects.

\bigskip

\noindent \textbf{Case $\mathbf{k(w)=0}$ or $\mathbf{1}$}.  First, suppose that $w$
has at most a single entry in the first column so $\pi_{dir}$ is just
the embedding map into type $A$.  Then we can interpret $w$ as an
element of type $A$ by projecting the unique generator in the first
column to the generator $s_1$ in type $A$.  Hence, $w$ is non-Deodhar
if and only if it contains a hexagon if and only if $\pi_{dir}(w)$
contains a hexagon, so $\pi_{dir}$ preserves the non-Deodhar condition
in this case.

Second, suppose that $w$ has two entries lying at the same level in
the first column and there are 4 entries in column 4.  Then $\pi_{NE}
= \pi_{SE} = \f_1$ just removes the single $s_{\tilde{1}}$ entry, and
by lateral convexity $\pi_{dir}(w)$ contains a hexagon.

Third, if there exist two entries in column 1 and there are fewer than
4 entries in column 4, then at least one of the points labeled $x$ or
$x'$ is not in the heap:
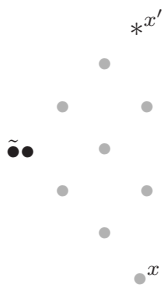
\[
\xymatrix @=-4pt @! {
& {\hs} & & {\hv}^{x'} & & {\hs} & & {\hs} \\
& & {\hb} & & {\hs} & & {\hs} & \\
& {\hb} & & {\hb} & & {\hs} & & {\hs} \\
\tilde{\hf} {\hf} &  & {\hb} & & {\hs} & \\
& {\hb} & & {\hb} & & {\hs} & & {\hs} \\
& & {\hb} & & {\hs} & & {\hs} & \\
& {\hs} & & {\hb}^x & & {\hs} & & {\hs} \\
} 
\]
By Lemma~\ref{l:w.winverse}, we have that $w$ is non-Deodhar if and
only if $w^{-1}$ is non-Deodhar.  Hence, by
Lemma~\ref{l:inverse.equivariant} we may assume that the entry $x'$ is
not in the heap, by taking a non-Deodhar mask on $w^{-1}$ if
necessary, whose heap is obtained from $Heap(w)$ by flipping the
diagram upside down.  
In this case, we choose the projection $\pi_{NE}(w)$ so
$Heap(\pi_{NE}(w))$ is obtained from $Heap(w)$ by shaving off the
maximal NE diagonal up to $x'$.  

We need to verify that the non-Deodhar bound holds on the mask
restricted to $\pi_{NE}(w)$. By convexity, there are at most three
occupied positions in the NE diagonal as in the figure above since
$x'$ is not in the heap.  There cannot be a defect for $w$ in column 1
since the defect would have no critical generator so there can be at
most two defects along the maximal NE diagonal.  If there exists a
defect whose left string touches either entry in the first column,
then both entries in the first column must have mask value 0, for
otherwise the left string of the defect would exit the heap without
crossing the right string, or be labeled negatively.  There could be a
single defect on the NE diagonal in column 3 whose left string does
not touch the entries in column 1, but then it must intersect a plain
zero along the NE diagonal in column 2.  Therefore, shaving off the NE
diagonal removes at least as many plain-zeros as zero-defects, so the
Deodhar bound is preserved.

\bigskip

\noindent \textbf{Case $\mathbf{k(w)=2}$}.
By Lemma~\ref{sbadh.2}, whenever the first column contains entries
on more than one distinct level, we must have that each level contains
a unique entry, and the entries alternate between the $s_1$ and
$s_{\tilde{1}}$ generators.  Without loss of generality, we can assume
assume the top entry in column 1 is $s_{\tilde{1}}$.

First, suppose $k(w)=2$ and there is a 4-stack in column 3.  By the
Shape Lemma~\ref{l:shape}, there is a hexagon in $Heap(\pi_{dir}(w))$.  

Second, suppose that there are at most 3 entries in column 3.  Then we
can further assume that the maximal NE diagonal from the
top entry in column 1 contains at most two entries by taking a
non-Deodhar mask on $w^{-1}$ if necessary using
Lemma~\ref{l:w.winverse} and Lemma~\ref{l:inverse.equivariant}.  The
projection $\pi_{NE}(w)$ shaves off the NE diagonal, so we need
to show that this removes at least as many plain-zeros as zero-defects
to preserve the bound.  There cannot be a defect in column 1, because
the defect would have no critical generator, so there can be at most
one defect in the maximal NE diagonal, located in position $a$:
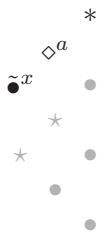
\[
\xymatrix @=-4pt @! {
& {\hs} & & {\hv} & & {\hs} & & {\hs} \\
& & {\hd}^a & & {\hs} & & {\hs} & \\
& {\tilde{\hf}^x} & & {\hb} & & {\hs} & & {\hs} \\
{\hs} &  & {\ha} & & {\hs} & \\
& {\ha} & & {\hb} & & {\hs} & & {\hs} \\
& & {\hb} & & {\hs} & & {\hs} & \\
& {\hs} & & {\hb} & & {\hs} & & {\hs} \\
} 
\]
If $x$ has mask value 1, then the left string of $a$ is negatively
signed while the right string has a positive label, in which case $a$
is not a defect by Lemma~\ref{l:defectlabellings}.  Hence, if $a$ is a
zero-defect, then $x$ must be a plain-zero so removing both $x$ and
$a$ from the heap will preserve the non-Deodhar bound.

\bigskip

\noindent \textbf{Case $\mathbf{k(w)=3}$}.  Suppose that $w$ has
entries on 3 distinct levels in the first column and assume without
loss of generality that the top entry corresponds to $s_{1}$.  The
projection $\pi_{NE}=\pi_{SE}$ adds the three entries marked as $\hh$
to the left of $Heap(w)$:
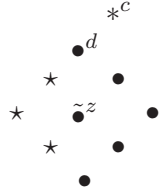
\[
\xymatrix @=-4pt @! {
& {\hs} & & {\hv}^c & & {\hs} & & {\hs} \\
& & {\hf}^d & & {\hs} & & {\hs} & \\
& {\hh} & & {\hf} & & {\hs} & & {\hs} \\
{\hh} & & \tilde{\hf}^z &  & {\hf} & & {\hs} & \\
& {\hh} & & {\hf} & & {\hs} & & {\hs} \\
& & {\hf} & & {\hs} & & {\hs} & \\
& {\hs} & & {\hs} & & {\hs} & & {\hs} \\
}
\]

If $w$ contains a 4-stack in column 2, then $\pi_{NE}(w)$ contains a
hexagon by the Shape Lemma~\ref{l:shape}.  Hence, we can assume there
are at most 3 entries in column 2.  By considering $w^{-1}$ if
necessary, and applying Lemma~\ref{l:w.winverse} and
Lemma~\ref{l:inverse.equivariant}, we can assume that the point marked
$c$ in column 2 is not in the heap.

Fix a non-Deodhar mask for $w$.  Since there are exactly three
alternating entries in the first column, there is at most a single
defect in the first column, and it must occur in position $d$ because
any defect requires a critical generator below it.  If $d$ is not a
defect then extend the mask to $\pi_{NE}(w)$ by setting the
mask-values of the three new points to 1:
\begin{equation}\label{e:pi.map.3}
\xymatrix @=-4pt @! {
& {\hs} & & {\hv}^c & & {\hs} & & {\hs} \\
& & {\ha}^d & & {\hs} & & {\hs} & \\
& {\hs} & & {\ha} & & {\hs} & & {\hs} \\
{\hs} & & \tilde{\ha}^z &  & {\ha} & & {\hs} & \\
& {\hs} & & {\ha} & & {\hs} & & {\hs} \\
& & {\ha} & & {\hs} & & {\hs} & \\
& {\hs} & & {\hs} & & {\hs} & & {\hs} \\
}
\hspace{.1in}
\begin{array}{c}
\\ \\
\pi_{NE}(w)\\
\longrightarrow
\end{array}
\hspace{.1in}
\xymatrix @=-4pt @! {
& {\hs} & & {\hv}^c & & {\hs} & & {\hs} \\
& & {\ha}^d & & {\hs} & & {\hs} & \\
& {\hf} & & {\ha} & & {\hs} & & {\hs} \\
{\hf} & & \tilde{\ha}^z &  & {\ha} & & {\hs} & \\
& {\hf} & & {\ha} & & {\hs} & & {\hs} \\
& & {\ha} & & {\hs} & & {\hs} & \\
& {\hs} & & {\hs} & & {\hs} & & {\hs} \\
} 
\end{equation}
If a string passes through point $z$ in $Heap(w)$, it changes sign and
so it cannot be the left string of any defect.  Therefore the mask
assignment in \eqref{e:pi.map.3} preserves all zero-defects and
plain-zeros, hence it preserves the non-Deodhar bound.

Now, suppose that $d$ is a defect in the first column.  Then we have
the following heap fragment where the critical generator below $d$
must have mask-value 1, but the mask-values of entries $x$, $y$ and
$z$ are variable:
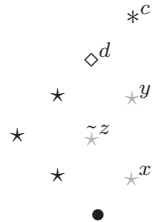
\[
\xymatrix @=-2pt @! {
& {\hs} & & {\hv}^c & & {\hs} & & {\hs} \\
& & {\hd}^d & & {\hs} & & {\hs} & \\
& {\hh} & & {\ha}^y & & {\hs} & & {\hs} \\
{\hh} & & {\tilde{\ha}}^z & & {\hs} & & {\hs} & \\
& {\hh} & & {\ha}^x & & {\hs} & & {\hs} \\
& & {\hf} & & {\hs} & & {\hs} & \\
& {\hs} & & {\hs} & & {\hs} & & {\hs} \\
} 
\]

Breaking into cases on the mask-values of $y$ and $z$, we demonstrate
how to extend the given non-Deodhar mask of $w$ to a non-Deodhar mask
of $\pi_{NE}(w)$.

If $z$ has mask-value 1, then $y$ must have mask-value 0 or else the
right string for $d$ will be labeled positively, while the left
string will be labeled negatively, contradicting that $d$ is a
defect.  In this case, the strings for $d$ cross at $z$, and both
become negatively labeled, so must cross once more below $z$.  Thus, we can set
the mask-values of the additional entries to the left as shown:
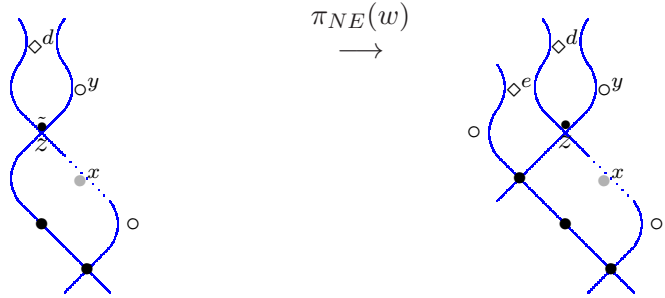
\[
\heap {
& {\hs} & & {\hs} & & {\hs} & & {\hs} \\
& & \StringLR{\hd^d} & & {\hs} & & {\hs} & \\
& \StringR{\hs} & & \StringL{\hz^y} & & {\hs} & & {\hs} \\
{\hs} & & \StringLRX{\stackrel{\tilde{\hf}}{z}} & & {\hs} & & {\hs} & \\
& \StringR{\hs} & & \StringLXD{\hb^x} & & {\hs} & & {\hs} \\
& & \StringLX{\hf} & & \StringL{\hz} & & {\hs} & \\
& {\hs} & & \StringLRX{\hf} & & {\hs} & & {\hs} \\
}
\hspace{.1in}
\begin{array}{c}
\\ \\
\pi_{NE}(w)\\
\longrightarrow
\end{array}
\hspace{.1in}
\heap {
& {\hs} & & {\hs} & & {\hs} & & {\hs} \\
& & \StringLR{\hd^d} & & {\hs} & & {\hs} & \\
& \StringLR{\hd^{e}} & & \StringL{\hz^y} & & {\hs} & & {\hs} \\
\StringR{\hz} & & \StringLRX{\stackrel{\hf}{z}} & & {\hs} & & {\hs} & \\
& \StringLRX{\hf} & & \StringLXD{\hb^x} & & {\hs} & & {\hs} \\
& & \StringLX{\hf} & & \StringL{\hz} & & {\hs} & \\
& {\hs} & & \StringLRX{\hf} & & {\hs} & & {\hs} \\
} \]
The string dynamics show that $d$ and $e$ are defects in
$Heap(\pi_{NE}(w))$.   The new mask has one additional plain-zero and
zero-defect so the non-Deodhar bound is maintained.

If both $y$ and $z$ have mask-value 0, then the strings for $d$ must
cross below $x$ in the heap poset.  In this case, we can set the
mask-values of the new entries as shown:
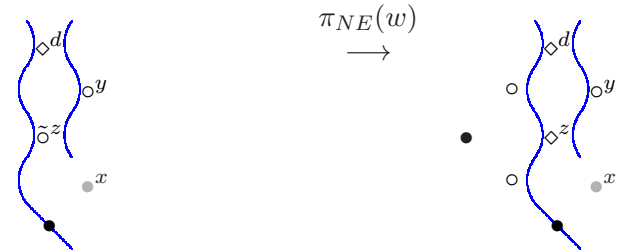
\[ 
\heap {
& {\hs} & & {\hs} & & {\hs} & & {\hs} \\
& & \StringLR{\hd^d} & & {\hs} & & {\hs} & \\
& \StringR{\hs} & & \StringL{\hz^y} & & {\hs} & & {\hs} \\
{\hs} & & \StringLR{\tilde{\hz}^z} & & {\hs} & & {\hs} & \\
& \StringR{\hs} & & {\hb}^x & & {\hs} & & {\hs} \\
& & \StringLX{\hf} & & {\hs} & & {\hs} & \\
& {\hs} & & {\hs} & & {\hs} & & {\hs} \\
}
\hspace{.1in}
\begin{array}{c}
\\ \\
\pi_{NE}(w)\\
\longrightarrow
\end{array}
\hspace{.1in}
\heap {
& {\hs} & & {\hs} & & {\hs} & & {\hs} \\
& & \StringLR{\hd^d} & & {\hs} & & {\hs} & \\
& \StringR{\hz} & & \StringL{\hz^y} & & {\hs} & & {\hs} \\
{\hf} & & \StringLR{\hd^z} & & {\hs} & & {\hs} & \\
& \StringR{\hz} & & {\hb}^x & & {\hs} & & {\hs} \\
& & \StringLX{\hf} & & {\hs} & & {\hs} & \\
& {\hs} & & {\hs} & & {\hs} & & {\hs} \\
} \]
The string dynamics for $d$ are preserved in the new mask, so $d$
remains a defect, and $z$ becomes a defect since it has exactly the
same string dynamics as $d$.  Note that $z$ is not a defect for $w$
since it corresponds to the $s_{\tilde{1}}$ generator in $Heap(w)$.
This mask is non-Deodhar because the original one was; we've changed a
plain zero to a defect and added two plain zero entries, maintaining
the non-Deodhar bound.

If $z$ has mask-value 0 and $y$ has mask-value 1, then either there
is another defect $e$ whose left string touches $z$, or there
is not.  If not, we can just move the zero at $z$ to the left to preserve
the string dynamics for $d$ and extend the mask to $\pi_{NE}(w)$ as shown:
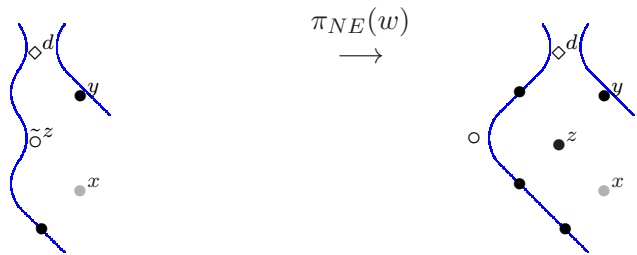
\[ 
\heap {
& {\hs} & & {\hs} & & {\hs} & & {\hs} \\
& & \StringLR{\hd^d} & & {\hs} & & {\hs} & \\
& \StringR{\hs} & & \StringLX{\hf^y} & & {\hs} & & {\hs} \\
{\hs} & & \StringL{\tilde{\hz}^z} & & {\hs} & & {\hs} & \\
& \StringR{\hs} & & {\hb}^x & & {\hs} & & {\hs} \\
& & \StringLX{\hf} & & {\hs} & & {\hs} & \\
& {\hs} & & {\hs} & & {\hs} & & {\hs} \\
}
\hspace{.1in}
\begin{array}{c}
\\ \\
\pi_{NE}(w)\\
\longrightarrow
\end{array}
\hspace{.1in}
\heap {
& {\hs} & & {\hs} & & {\hs} & & {\hs} \\
& & \StringLR{\hd^d} & & {\hs} & & {\hs} & \\
& \StringRX{\hf} & & \StringLX{\hf^y} & & {\hs} & & {\hs} \\
\StringR{\hz} & & {\hf}^z & & {\hs} & & {\hs} & \\
& \StringLX{\hf} & & {\hb}^x & & {\hs} & & {\hs} \\
& & \StringLX{\hf} & & {\hs} & & {\hs} & \\
& {\hs} & & {\hs} & & {\hs} & & {\hs} \\
} \]
This mask is non-Deodhar assuming the original one was since we have not
changed the number of zero-defects or plain-zeros.  

On the other hand, if $z$ is touched by the left string of a defect at
or above $e$, then we claim $Heap(w)$ also contains an entry $p$ below $x$ in
column 3.  Hence, the heap contains an I-shape with corners $c$, $d$,
$e$ and $p$, so the Shape Lemma~\ref{l:shape} implies that
$\pi_{NE}(w)$ is non-Deodhar.  For example, 
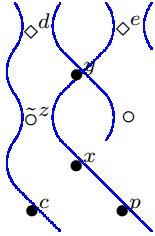
\[ \heap {
& {\hs} & & {\hs} & & {\hs} & & {\hs} \\
& & \StringLR{\hd^d} & & \StringLR{\hd^e} & & {\hs} & \\
& \StringR{\hs} & & \StringLRX{\hf^y} & & {\hs} & & {\hs} \\
{\hs} & & \StringLR{\tilde{\hz}^z} & & \StringL{\hz} & & {\hs} & \\
& \StringR{\hs} & & \StringLX{\hf^x} & & {\hs} & & {\hs} \\
& & \StringLX{\hf^c} & & \StringLX{\hf^p} & & {\hs} & \\
& {\hs} & & {\hs} & & {\hs} & & {\hs} \\
} 
\] 
To verify the existence of $p \in Heap(w)$, observe that the
highest possible crossing for the strings of $d$ is at the point $c$, and
consider two cases.  If the strings of $d$ do not cross at $c$,
they must cross at a point below $c$ and $p$ so by convexity $p\in Heap(w)$.
For the strings of $d$ to cross at $c$, the point $x$ must
have mask-value 1 and the right critical zero of $d$ is located as
shown, in which case the strings for the other defect $e$ cross at or
below the point $p$, so again $p \in Heap(w)$.

Thus, for all possible mask-values of $x,y$ the projection $\pi_{NE}$
preserves the non-Deodhar condition, finishing the case $k(w)=3$.

\bigskip

\noindent \textbf{Case $\mathbf{k(w)=4}$}.  Suppose that $w$ has entries on 4
distinct levels in the first column.  

First, we reduce to the case where $w$ contains a decorated heap
fragment of the form:
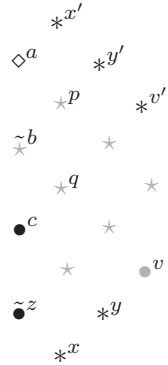
\begin{equation}\label{e:pi.map.4}
\xymatrix @=-4pt @! {
& {\hv}^{x'} & & {\hs} & & {\hs} & & {\hs} \\
{\hd}^a & & {\hv}^{y'} & & {\hs} & & {\hs} \\
& {\ha}^p & & {\hv}^{v'} & & {\hs} & & {\hs} \\
{\tilde{\ha}}^b & & {\ha} & & {\hs} & & {\hs} \\
& {\ha}^q & & {\ha} & & {\hs} & & {\hs} \\
{\hf}^c & & {\ha} & & {\hs} & & {\hs} \\
& {\ha} & & {\hb}^{v} & & {\hs} & & {\hs} \\
{\tilde{\hf}}^z & & {\hv}^{y} & & {\hs} & & {\hs} \\
& {\hv}^{x} & & {\hs} & & {\hs} & & {\hs} \\
} 
\end{equation}
Then we study cases corresponding to the mask values in the gray
star positions.  

If there is a point at $x$ or $x'$, then there is a 4-stack in the
second column.  If there is a point at $y$ or $y'$, we obtain an
I-shape with the other entries that exist by convexity. In either
case, the Shape Lemma~\ref{l:shape} implies
$\pi_{dir}(w)$ contains a hexagon.  Similarly, if both $v$ and $v'$ exist,
then we obtain an I-shape, so at least one of them is not in the heap.
By considering $w^{-1}$ if necessary, we can assume that $v'$ is not
in the heap.

By Definition~\ref{d:pi} for $\pi_{dir}$, we must first choose one of the
extremal diagonals to shave, and then additively resolve the remaining
heap fragment, in such a way that the non-Deodhar bound is maintained.
Hence, choosing a projection $\pi_{NE}$ or $\pi_{SE}$ really just
amounts to choosing the top or bottom entry from the first column
to remove in such a way that the non-Deodhar bound is maintained.
Once we choose the entry and shave it, we will have a non-Deodhar heap
with only three entries in the first column, so we can appeal
to the previous case $k(w) = 3$ when we apply $\pi_{dir}$ the second time.

If $z$ has mask-value 0, then it can play no critical role in any
defect: it cannot be a defect itself since there is no critical
generator below, and it cannot be a critical-zero for a defect since
both strings meeting at $z$ leave the heap below.  Therefore, we
choose to shave $z$ in this case which maintains the non-Deodhar bound.

Hence, we can assume that $z$ has mask-value 1 and that it is a string
crossing for some zero-defect, since we can again remove it if it is
not.  Note that if $a$ is not a zero-defect, we can choose to shave it
without changing the Deodhar bound.  Therefore, we can also assume
that $a$ is a zero-defect, which forces $c$ to have mask-value 1 since
it must be the critical generator for $a$.  Thus, we have a heap
fragment of the form in \eqref{e:pi.map.4}.

Suppose that $b$ is a plain zero.  Then, the strings that cross at
$z$ must both be labeled negatively everywhere above $z$ since there
is no other active $s_{\tilde{1}}$ generator to change the signs of
the strings.  Therefore, no defect above $z$ can be removed by
removing $z$ since the strings remain in increasing order as they pass
through $z$ and all defects in this case must correspond with one of
the type $A$ generators, $s_{1}, s_{2}, \dotsc, s_{n-1}$.  Hence, we
can shave $z$ without affecting the defect status of any defect whose
strings cross at $z$, and so we maintain the Deodhar bound in this
case.

Next, suppose $b$ is a zero-defect.  Then its strings must cross at
$z$ since this is its critical generator, and there are no lower
entries.  Moreover, the right critical zero for $a$ must occur at $p$,
or the strings for $a$ will never cross. Hence, the paths of the
strings from $b$ are completely prescribed and we have the heap
fragment:
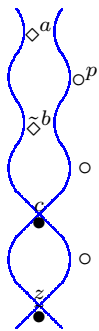
\[ \heap {
& \StringLR{\hd^a} & & {\hs} & & {\hs} & & {\hs} \\
\StringR{\hs} & & \StringL{\hz^p} & & {\hs} & & {\hs} & & {\hs} \\
& \StringLR{\tilde{\hd}^b} & & {\hs} & & {\hs} & & {\hs} \\
\StringR{\hs} & & \StringL{\hz} & & {\hs} & & {\hs} & & {\hs} \\
& \StringLRX{\stackrel{c}{\hf}} & & {\hs} & & {\hs} & & {\hs} \\
\StringR{\hs} & & \StringL{\hz} & & {\hs} & & {\hs} & & {\hs} \\
& \StringLRX{\stackrel{z}{\tilde{\hf}}} & & {\hs} & & {\hs} & & {\hs} \\
}  \]
 Since $p$ cannot be a defect
(because it lacks a critical generator), and it cannot enable a left
string crossing for any defect above $p$ (or we introduce an I-shape),
removing $a$ and changing the mask-value at $p$ to 1 maintains the
Deodhar bound.

Finally, suppose that $b$ has mask-value 1, and consider the possible
locations of the zero-defect whose strings cross at $z$; call it $d_{z}$.
If $d_{z}$ is located on the NE diagonal from $c$, then its left
string is labeled $\overline{2}$, which is a contradiction.
If $d_{z}$ is located at $p$, then the remaining mask values are
forced, and in particular, $q$ must have mask-value 1 as shown in
Figure~\ref{f:ap}(a).
\begin{figure}[h]
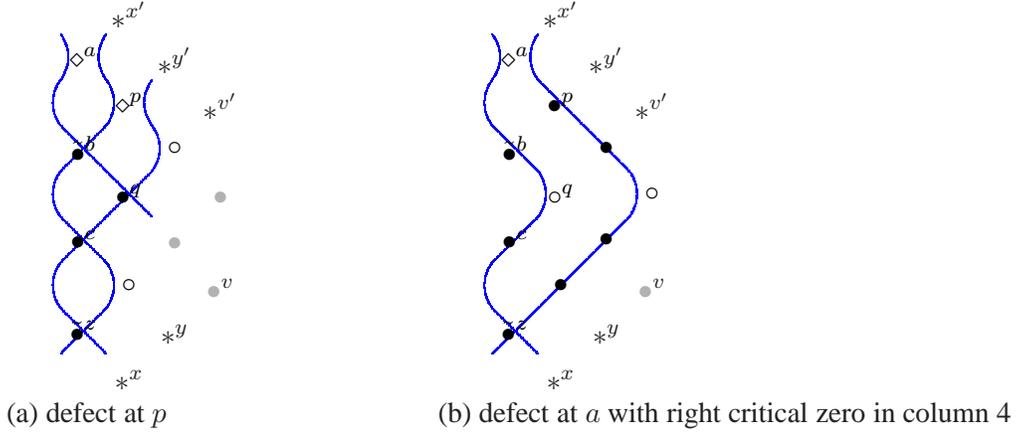

\begin{tabular}{ll}
\heap {
& & {\hv}^{x'} & & {\hs} & & {\hs} & & {\hs} \\
& \StringLR{\hd^a} & & {\hv}^{y'} & & {\hs} & & {\hs} \\
\StringR{\hs} & & \StringLR{\hd^p} & & {\hv}^{v'} & & {\hs} & & {\hs} \\
& \StringLRX{\tilde{\hf}^b} & & \StringL{\hz} & & {\hs} & & {\hs} \\
\StringR{\hs} & & \StringLRX{\hf^q} & & {\hb} & & {\hs} & & {\hs} \\
& \StringLRX{\hf^c} & & {\hb} & & {\hs} & & {\hs} \\
\StringR{\hs} & & \StringL{\hz} & & {\hb}^{v} & & {\hs} & & {\hs} \\
& \StringLRX{\tilde{\hf}^z} & & {\hv}^{y} & & {\hs} & & {\hs} \\
& & {\hv}^{x} & & {\hs} & & {\hs} & & {\hs} \\
} &
\heap {
& & {\hv}^{x'} & & {\hs} & & {\hs} & & {\hs} \\
& \StringLR{{\hd}^a} & & {\hv}^{y'} & & {\hs} & & {\hs} \\
\StringR{\hs} & & \StringLX{\hf^p} & & {\hv}^{v'} & & {\hs} & & {\hs} \\
& \StringLX{\tilde{\hf}^b} & & \StringLX{\hf} & & {\hs} & & {\hs} \\
{\hs} & & \StringL{\hz^q} & & \StringL{\hz} & & {\hs} & & {\hs} \\
& \StringRX{\hf^c} & & \StringRX{\hf} & & {\hs} & & {\hs} \\
\StringR{\hs} & & \StringRX{\hf} & & {\hb}^{v} & & {\hs} & & {\hs} \\
& \StringLRX{\tilde{\hf}^z} & & {\hv}^{y} & & {\hs} & & {\hs} \\
& & {\hv}^{x} & & {\hs} & & {\hs} & & {\hs} \\
} \\ 
(a) defect at $p$ & 
(b) defect at $a$ with right critical zero in column 4 \\ 
\end{tabular}
\caption{Locations of defects with strings crossing at $z$}
\label{f:ap}
\end{figure}
Then, we obtain the contradiction that $a$ cannot be a defect in this
case since its left string is labeled $\overline{2}$, while its right string is
labeled 1.  Thus $d_{z} =a$ is the only viable possibility.

Suppose $d_{z}=a$.  If the right critical zero for $a$ occurs in
column 4, then the mask-values are forced and in particular $q$ must
be a plain-zero as shown in Figure~\ref{f:ap}(b).  Hence, we can shave
$a$ and change the mask-value of $q$ to 1, since $q$ is not itself a
defect, nor can any other defect use $q$ to effect a string crossing.
This choice of projection preserves the Deodhar bound.

If the right critical zero for $a$ is located in column $2$ or $3$ and
it is a plain-zero, then we can shave $a$ and change the mask-value of
the right critical zero to 1, since no other defect can use the right
critical zero of $a$ to effect a string crossing.  Hence, we may
assume that the right critical zero of $a$ is itself a zero-defect.

If the right critical zero for $a$ is in column 3, then we have one of the
cases shown in Figure \ref{f:ac3}.
\begin{figure}[h]
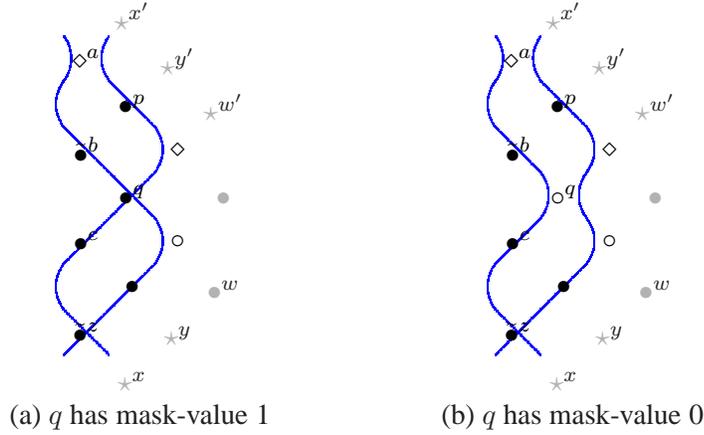

\begin{tabular}{ll}
\heap {
& & {\ha}^{x'} & & {\hs} & & {\hs} & & {\hs} \\
& \StringLR{{\hd}^a} & & {\ha}^{y'} & & {\hs} & & {\hs} \\
\StringR{\hs} & & \StringLX{\hf^p} & & {\ha}^{w'} & & {\hs} & & {\hs} \\
& \StringLX{\tilde{\hf}^b} & & \StringL{\hd} & & {\hs} & & {\hs} \\
{\hs} & & \StringLRX{\hf^q} & & {\hb} & & {\hs} & & {\hs} \\
& \StringRX{\hf^c} & & \StringL{\hz} & & {\hs} & & {\hs} \\
\StringR{\hs} & & \StringRX{\hf} & & {\hb}^{w} & & {\hs} & & {\hs} \\
& \StringLRX{\tilde{\hf}^z} & & {\ha}^{y} & & {\hs} & & {\hs} \\
& & {\ha}^{x} & & {\hs} & & {\hs} & & {\hs} \\
} &
\heap {
& & {\ha}^{x'} & & {\hs} & & {\hs} & & {\hs} \\
& \StringLR{{\hd}^a} & & {\ha}^{y'} & & {\hs} & & {\hs} \\
\StringR{\hs} & & \StringLX{\hf^p} & & {\ha}^{w'} & & {\hs} & & {\hs} \\
& \StringLX{\tilde{\hf}^b} & & \StringL{\hd} & & {\hs} & & {\hs} \\
{\hs} & & \StringLR{\hz^q} & & {\hb} & & {\hs} & & {\hs} \\
& \StringRX{\hf^c} & & \StringL{\hz} & & {\hs} & & {\hs} \\
\StringR{\hs} & & \StringRX{\hf} & & {\hb}^{w} & & {\hs} & & {\hs} \\
& \StringLRX{\tilde{\hf}^z} & & {\ha}^{y} & & {\hs} & & {\hs} \\
& & {\ha}^{x} & & {\hs} & & {\hs} & & {\hs} \\
} \\
(a) $q$ has mask-value 1 & 
(b) $q$ has mask-value 0 \\ 
\end{tabular}
\caption{$a$ has a right-critical zero that is a defect in column 3}
\label{f:ac3}
\end{figure}
Once we choose a mask-value for $q$, the rest of the mask-values are
determined.  In each case, the right critical zero for $a$ cannot be a
defect because its strings do not cross, which is a contradiction.

If the right critical zero for $a$ is in column 2, then we have the
cases shown in Figure \ref{f:ac2}.
\begin{figure}[h]
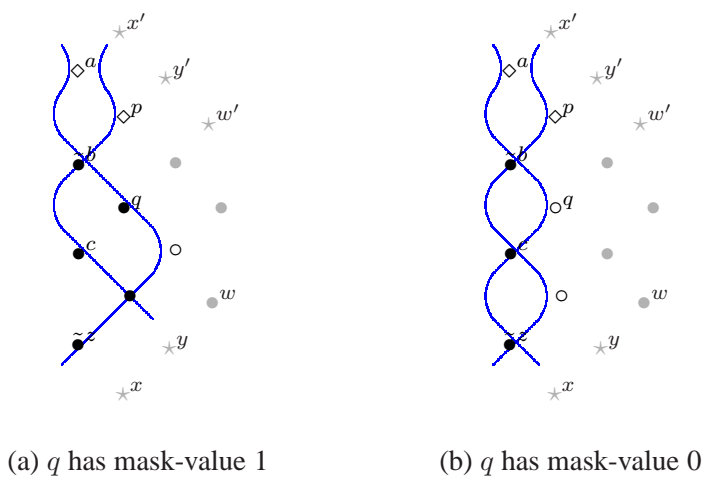

\begin{tabular}{ll}
\heap {
& & {\ha}^{x'} & & {\hs} & & {\hs} & & {\hs} \\
& \StringLR{{\hd}^a} & & {\ha}^{y'} & & {\hs} & & {\hs} \\
\StringR{\hs} & & \StringL{\hd^p} & & {\ha}^{w'} & & {\hs} & & {\hs} \\
& \StringLRX{\tilde{\hf}^b} & & {\hb} & & {\hs} & & {\hs} \\
\StringR{\hs} & & \StringLX{\hf^q} & & {\hb} & & {\hs} & & {\hs} \\
& \StringLX{\hf^c} & & \StringL{\hz} & & {\hs} & & {\hs} \\
{\hs} & & \StringLRX{\hf} & & {\hb}^{w} & & {\hs} & & {\hs} \\
& \StringRX{\tilde{\hf}^z} & & {\ha}^{y} & & {\hs} & & {\hs} \\
& & {\ha}^{x} & & {\hs} & & {\hs} & & {\hs} \\
& {\hs} & & {\hs} & & {\hs} & & {\hs} \\
} &
\heap {
& & {\ha}^{x'} & & {\hs} & & {\hs} & & {\hs} \\
& \StringLR{{\hd}^a} & & {\ha}^{y'} & & {\hs} & & {\hs} \\
\StringR{\hs} & & \StringL{\hd^p} & & {\ha}^{w'} & & {\hs} & & {\hs} \\
& \StringLRX{\tilde{\hf}^b} & & {\hb} & & {\hs} & & {\hs} \\
\StringR{\hs} & & \StringL{\hz^q} & & {\hb} & & {\hs} & & {\hs} \\
& \StringLRX{\hf^c} & & {\hb} & & {\hs} & & {\hs} \\
\StringR{\hs} & & \StringL{\hz} & & {\hb}^{w} & & {\hs} & & {\hs} \\
& \StringLRX{\tilde{\hf}^z} & & {\ha}^{y} & & {\hs} & & {\hs} \\
& & {\ha}^{x} & & {\hs} & & {\hs} & & {\hs} \\
& {\hs} & & {\hs} & & {\hs} & & {\hs} \\
} \\
(a) $q$ has mask-value 1 & 
(b) $q$ has mask-value 0 \\ 
\end{tabular}
\caption{$a$ has a right-critical zero that is a defect in column 2}
\label{f:ac2}
\end{figure}
If $q$ has mask-value 1 then the strings for $a$ cannot cross at $z$,
and if $q$ has mask-value 0, then the right critical zero for $a$
cannot be a defect because its strings do not cross.

Thus, in all cases there is a choice of projection $\pi_{NE}$ or
$\pi_{SE}$ which preserves the non-Deodhar condition.

\bigskip

\noindent \textbf{Case $\mathbf{k(w)=5}$}.  Suppose that $w$ has entries
on 5 or more distinct levels in the first column.  Then, $w$ contains
the pattern $HEX_4$ and so the projections $\pi_{dir}(w)$ contain a hexagon.

\bigskip

We have shown in all cases that there is a choice of projection
$\pi_{NE}(w)$ or $\pi_{SE}(w)$ that remains non-Deodhar if $w$ is
non-Deodhar, concluding the proof.
\end{proof}

\bigskip
\section{Proof of the type $D$ characterization theorem}\label{s:d.proof}

In this section, we complete the classification of the minimal
non-Deodhar embedded factors for type $D$.

\begin{theorem}\label{sbadh.th}
Suppose $w$ is a short braid avoiding, type $D$ element that is not
convex.  Then, $w$ is non-Deodhar if and only if $w$ contains
$[\bar{1} 6 7 8 \bar{5} 2 3 4]$ as a 1-line pattern.
\end{theorem}
\begin{proof}
Since $w$ is not convex, we have that $Heap(w)$ contains a minimal
pair of entries that require a left resolution.

Suppose $w$ is non-Deodhar.  By Lemma~\ref{sbadh.3}, $Heap(w)$ must
have an $s_1 s_{\tilde{1}}$ entry.  By Lemma~\ref{sbadh.4}, $w$
contains a convex non-Deodhar element $\check{w}$ as a 1-line pattern.
Moreover, the construction of $\check{w}$ preserves the $s_1
s_{\tilde{1}}$ entry.  The only convex minimally non-Deodhar heap from
Theorem~\ref{t:minimal.patterns} that has an entry with an $s_1
s_{\tilde{1}}$ entry is:
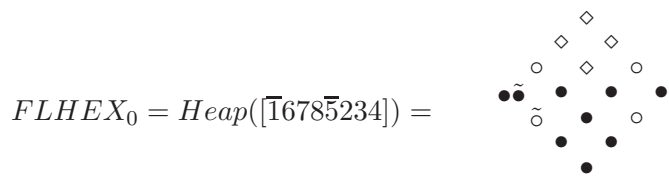
\[
\parbox[t]{2.5in}{ \vspace{.5in} 
$FLHEX_{0}=Heap([\overline{1} 6 7 8 \overline{5} 2 3 4])=$
}
  \xymatrix @=-8pt @! {
{\hs} & & {\hs} & & {\hs} & & {\hs} \\
& {\hs} & & {\hd} & & {\hs} &\\
{\hs} & & {\hd} & & {\hd} & & {\hs} \\
& {\hz} & & {\hd} & & {\hz}\\
{\hf \tilde{\hf} }&  & {\hf} & & {\hf} & & {\hf} & {\hs} \\
& {\tilde{\hz}} & & {\hf} & & {\hz} &\\
{\hs} & & {\hf} & & {\hf} & & {\hs} \\
& {\hs} & & {\hf} & & {\hs} &\\
} 
\]
Hence, $\check{w}$ must contain a $FLHEX_{0}$ factor. Note
the relative order of the 8 strings passing through a $FLHEX_{0}$ gives
rise to a $[\overline{1} 6 7 8 \overline{5} 2 3 4]$ pattern in
$\check{w}$.  The relative order of these strings cannot change due to
multiplication by additional generators above or below in a short
braid avoiding heap since any additional adjacent transposition
applied on the left or right of $[\overline{1} 6 7 8 \overline{5} 2 3
4]$ would create a short braid.  Therefore, both $\check{w}$ and $w$
contain $[\overline{1} 6 7 8 \overline{5} 2 3 4]$ as a 1-line
pattern.

Conversely, suppose that $w$ is a short braid avoiding type $D$
element that contains $[\overline{1} 6 7 8 \overline{5} 2 3 4]$ as a
1-line pattern.  We will reduce to the case where $w_{1}=\overline{1}$
and the values in the pattern representing $678\overline{5}234$ occur
consecutively in the 1-line notation for $w$.  Then, we will be able to
give a reduced factorization for $w$ that contains a $FLHEX_{k}$ which
is known to be non-Deodhar; see Example~\ref{e:infantichain}.

Consider the 1-line notation for $w$ and highlight a pattern instance
for $[\overline{1} 6 7 8 \overline{5} 2 3 4]$ on values
$a_{1}<a_{2}<\ldots < a_{8}$ and in positions $p_{1}< \ldots < p_{8}$:

\[
w=[ \dots \overline{a_{1}} \dots a_{6} \dots a_{7} \dots a_{8} \dots
\overline{a_{5}} \dots a_{2} \dots a_{3} \dots a_{4} \dots ]
\]
Since all 5 of the signed permutations
\[ [\overline{1} \overline{2} \overline{3}],
[\overline{2} \overline{1} \overline{3}],
[\overline{1} \overline{3} \overline{2}],
[\overline{2} \overline{3} \overline{1}],
[\overline{3} \overline{1} \overline{2}] \]
with at least one ascent and all entries negative are forbidden by
Theorem~\ref{t:stem}, we have that all of the entries besides
$\overline{a_{1}}, \overline{a_{5}}$ must be positive.  In addition,
since $[1 \overline{2} \overline{3}]$, $[\overline{2} 1 \overline{3}]$
and $[\overline{2} \overline{3} 1]$ are all forbidden by
Theorem~\ref{t:stem}, we have that $a_1 = 1$.  Hence, if we let
$i=p_{1}$ then $w$ contains $w'=ws_{i-1}s_{i-2}\cdots s_{1}$ as a
factor, and $w'$ also contains the 1-line pattern $[\overline{1} 6 7 8
\overline{5} 2 3 4]$ with $w'_{1}=\overline{1}$.  Therefore, we only
need to consider the case when $w_{1} = a_1 = \overline{1}$.

Next, if there exist entries in $w$ between $\overline{a_{5}}$ and
$a_{2}$, let $t=p_{6}-1$ be the position just to the left of $a_{2}$.
By Theorem~\ref{t:stem}, $a_8 > w_{t} > a_{2}$ would imply a forbidden
$[321]$ pattern, so either
\begin{enumerate}
\item$w_{t}> a_{8}$ in which case $w'=ws_{t}$ also has the pattern and
is a factor of $w$, or 
\item $1<w_{t}<a_{2}$ in which case we can change the highlighted pattern
to choose a pattern instance where the $a_{2}$ is closer to the
$\overline{a_{5}}$.
\end{enumerate}
Therefore, we only need to consider the case when
$\overline{a_{5}}a_{2}$ are consecutive.  Furthermore, by the same
argument in which the role of $a_2$ is replaced by $a_3$ and $a_4$,
respectively, we can assume $\overline{a_{5}}a_{2}a_{3}a_{4}$ are
consecutive in $w$.  Note also that all of the entries in positions $>
p_8$ must have values $> a_4$, for otherwise we obtain a forbidden
$[321]$ instance, with $a_8 > a_4$.

By Theorem~\ref{t:stem}, $w$ avoids the patterns $[2 1 \overline{3}],
[3 1 \overline{2}], [3 2 \overline{1}]$ so the values in $w$ between
the $\overline{1}$ and $\overline{a_{5}}$ are increasing positive
numbers.  By changing the highlighted pattern instance if necessary,
we can assume the values $a_{6} a_{7} a_{8}$ are the biggest three
among these so $a_{6} a_{7} a_{8} \overline{a_{5}}$ are in consecutive
positions.  Consider the value in position $p_{2}-1$, just to the
left of $a_{6}$.  If $a_{5} > w_{p_{2}-1} > a_{2}$ then $w$ contains
the 1-line pattern $[2 \overline{3} 1]$ which is forbidden by
Theorem~\ref{t:stem}.  If $w_{p_{2}-1} > a_5$ then we move the $a_8$
out of the pattern to the right, and change the highlighted pattern so
that $w_{p_{2}-1}$ becomes $a_6$, $a_6$ becomes $a_7$, and $a_7$
becomes $a_8$.  Specifically, if $i = p_4$ then $w$ contains
$w'=ws_{i}s_{i+1}s_{i+2}s_{i+3}$ as a factor and $w'$ also contains
$[\overline{1} 6 7 8 \overline{5} 2 3 4]$ as a 1-line pattern in the
manner described.  Hence, we only need to consider the case when the
entries in $w$ between $\overline{1}$ and $a_{6}$ all have values less
than $a_{2}$.

Summarizing, we can assume the 1-line notation for $w$ is of the form
$[\overline{1} \dots a_{6} a_{7} a_{8} \overline{a_{5}} a_{2} a_{3} a_{4}
\dots]$ with the following conditions:
\begin{enumerate}
\item The elements in the first dotted sequence are increasing, all
with positive values less than $a_{2}$.
\item All entries in the second dotted sequence have value greater
than $a_{4}>0$.
\end{enumerate}

When we draw the string diagram corresponding to these rules, we find
that $w$ contains $FLHEX_k$ as a factor.

Alternatively, we obtain a reduced factorization of $w$ containing a
non-Deodhar factor as follows.  Let $1<y_{1}<\dots< y_{j}<a_{5}$ be
the values that appear to the right of $\overline{a_{5}}$ in $w$.
Note that $a_{2},a_{3}, a_{4}$ are the smallest three of the
$y_{i}$'s.  Let $u$ be the permutation with values $a_{6}, a_{7},
a_{8}$ moved left and consecutively stacked adjacent to $a_{5}$ and
with values $y_{1},\dots, y_{k}$ moved right and consecutively stacked
adjacent to $a_{5}$ , i.e.
\[
u=[123\dotsb a_{2}a_{3}a_{4}y_{4}\dotsb y_{j}a_{5}a_{6}a_{7}a_{8}\dots n].
\]
Next, move the block $y_{4}\dotsb y_{j}$ across the block $a_{5}a_{6}a_{7}a_{8}$ to get 
\[
u'=[123\dotsb a_{2}a_{3}a_{4}a_{5}a_{6}a_{7}a_{8}y_{4}\dotsb y_{j}\dots n].
\]
Then, we have
\[ w = u' \cdot FLHEX_{k} \cdot v \]
for some permutation $v$ that arranges the entries to the right of
$a_{4}$ in their final order.  Here $k=a_{5}-2 -j$ is the number of
strings strictly between 1 and $a_{5}$ that remain between positions
$1$ and $p_{5}$.  The $FLHEX_{k}$ pattern represents the operation 
that starting from $u'$ moves $a_{5}$ to position 2, shifts the block
$a_{2}a_{3}a_{4}$ to the right of $a_{6}a_{7}a_{8}$, applies
$s_{1}s_{\tilde{1}}$ to change the signs on $1$ and $a_{5}$ in
positions $1$ and $2$, and finally moves $\overline{a_{5}}$ back to
position $k+3$.
\end{proof}

\begin{theorem}\label{t:main.d}
In type $D$, we have that $w$ is Deodhar if and only if it avoids the
1-line pattern $[\bar{1} 6 7 8 \bar{5} 2 3 4]$, and the six embedded
factors 
\begin{align*}
& s_1 s_2 s_1 \text{ (short braid) } \\
& s_5 s_6 s_7 s_3 s_4 s_5 s_6 s_2 s_3 s_4 s_5 s_1 s_2 s_3  \text{ (HEX) } \\
& s_3 s_4 s_5 s_6 s_2 s_3 s_4 s_5 s_{\tilde{1}} s_2 s_3 s_4 s_1 s_2 s_3  \text{ ($HEX_2$) } \\
& s_4 s_5 s_6 s_2 s_3 s_4 s_5 s_1 s_2 s_3 s_4 s_{\tilde{1}} s_2 s_1  \text{ ($HEX_{3a}$) } \\
& s_1 s_4 s_5 s_6 s_2 s_3 s_4 s_5 s_{\tilde{1}} s_2 s_3 s_4 s_1 s_2  \text{ ($HEX_{3b}$)} \\
& s_3 s_4 s_5 s_1 s_2 s_3 s_4 s_{\tilde{1}} s_2 s_3 s_1  \text{ ($HEX_4$). } 
\\
\end{align*}
\end{theorem}
\begin{proof}
If $w$ contains a short braid then $w$ is not Deodhar by
Corollary~\ref{c:short.braid}.  If $w$ is convex, then $w$ is Deodhar
if and only if it avoids the patterns $HEX, HEX_{2},HEX_{3a},HEX_{3b},
HEX_{4}$ by Theorem~\ref{t:minimal.patterns}.  If $w$ is short braid
avoiding and not convex, then $w$ is Deodhar if and only if it avoids
$[\bar{1} 6 7 8 \bar{5} 2 3 4]$ by Theorem~\ref{sbadh.th}.
\end{proof}

\bigskip \section{Deodhar elements of exceptional Weyl
groups}\label{s:classification.e}

The minimal non-Deodhar patterns for the other finite Weyl groups are
computed by software that implements a game of K. Eriksson described
in \cite{eriksson} and \cite{b-b}.  We give the minimal lists for each
type that account for Coxeter graph isomorphisms and patterns
contained in parabolic subgroups.  For example, a Deodhar element of
$E_8$ must avoid the short braid on each pair of non-commuting
generators, the hexagon pattern of length 14 contained in the type $A$
parabolic subgroup, the type $D$ patterns, and  all of the
patterns in $E_6$ and $E_7$, the former of which can be embedded in
two different ways.

\begin{theorem}\label{t:main.e}
Below is the complete list of minimal non-Deodhar embedded factor
patterns in the Weyl groups of type $E$.  The only minimal non-Deodhar
elements in types $G_{2}$ and $F_{4}$ are short braids.

\bigskip
\begin{tabular}{|p{0.7in}|p{2in}|p{2.4in}|} 
	\hline
	Lie type & Coxeter graph & Reduced expression patterns \\
	\hline
	$E_6$ & 
	\xymatrix @-1pc {
    & & \gn_{5} &        &                &                            &  \\
    \gn_0 \ar@{-}[r] & \gn_1 & \ar@{-}[l] \ar@{-}[u] \gn_2 \ar@{-}[r] & \gn_3 \ar@{-}[r] & \gn_4 & \\
    }
    & { $s_0 s_1 s_2 s_5 s_3 s_4 s_2 s_3 s_1 s_2 s_5 s_0 s_1$ 
	 $s_5 s_1 s_2 s_3 s_0 s_1 s_2 s_5 s_4 s_3 s_2 s_1 s_0$ 
	 $s_1 s_2 s_5 s_3 s_4 s_2 s_3 s_1 s_2 s_5 s_0 s_1 s_2$ 
	 $s_2 s_5 s_1 s_2 s_3 s_0 s_1 s_2 s_5 s_4 s_3 s_2 s_1$ } \\
	\hline
	$E_7$ & \xymatrix @-1pc {
    & & \gn_{5} &        &                &                            &  \\
    \gn_0 \ar@{-}[r] & \gn_1 & \ar@{-}[l] \ar@{-}[u] \gn_2 \ar@{-}[r] & \gn_3 \ar@{-}[r] & \gn_4 & \gn_6 \ar@{-}[l] \\
    }
    & { $s_0 s_1 s_2 s_3 s_4 s_6 s_5 s_2 s_3 s_4 s_1 s_2 s_3 s_0 s_1$
	$s_3 s_4 s_6 s_1 s_2 s_3 s_0 s_1 s_2 s_5 s_4 s_3 s_2 s_1 s_0$ 
	$s_1 s_2 s_3 s_4 s_6 s_5 s_2 s_3 s_4 s_1 s_2 s_3 s_0 s_1 s_2$
	$s_2 s_3 s_4 s_6 s_1 s_2 s_3 s_0 s_1 s_2 s_5 s_4 s_3 s_2 s_1$
	$s_5 s_2 s_3 s_4 s_6 s_1 s_2 s_5 s_3 s_4 s_2 s_3 s_0 s_1 s_2 s_5$ }  \\
	\hline
\end{tabular}
\bigskip

\end{theorem}
\begin{proof}
These groups are finite, so this list of minimal patterns is
verifiable by computer.   The code used is available at 
\url{http://www.math.washington.edu/~brant/liberikson.html}.
\end{proof}

\bigskip
\section{Patterns of convex elements}\label{s:convex.patterns}

In this section, we prove that the Deodhar property can be
characterized by avoiding finitely many 1-line patterns in types $A$,
$B$ and $D$.  In general, properties characterized by 1-line pattern
avoidance are not equivalent to properties characterized by embedded
factor avoidance, even if the embedded factors are convex, as seen in
the example below.  However, we describe a finite test for when we may
translate between the two types of pattern avoidance on convex
elements.  Although it is not a direct generalization, this idea is
related to a type $A$ result of Tenner \cite{t2}.  In \cite{j2} the
main result of this section has been extended for type $A$ to cases
where the elements may not be convex.  

\begin{example}\label{ex:badp}
Consider the subset of permutations $S_n(p)$ avoiding $p = s_1 s_3 s_2
s_4 = [2 4 1 5 3]$ as an embedded factor.  Then, the element $w = s_4
s_1 s_3 s_5 s_2 = [2 5 1 3 6 4]$ avoids $p$ as an embedded factor, so
$w \in S_n(p)$, yet it contains $p$ as a 1-line pattern.  The string
diagrams below depict how the extra string is added to $p$ to get
$w$:

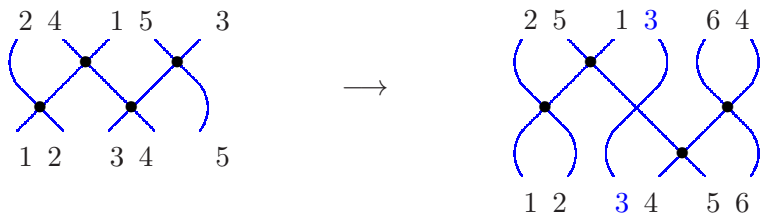
\[
\heap {
  & 2 \ \ 4 &  & 1 \ \ 5 & & 3   \\
  \StringR{\hs} & {\hs} & \StringLRX{\hf} & \hs & \StringLRX{\hf} & {\hs} & {\hs} \\
  &   \StringLRX{\hf} & \hs & \StringLRX{\hf} & \hs & \StringL{\hs} \\
  & 1 \ \ 2 &  & 3 \ \ 4 & & 5   \\
} 
\parbox[t]{1in}{ \vspace{0.3in} \hspace{.3in} $\longrightarrow$ }
\heap {
  & 2 \ \ 5 &  & 1 \ \ {\color{blue} 3} & & 6 \ \ 4 \\
  \StringR{\hs} & {\hs} & \StringLRX{\hf} & \hs & \StringLR{\hs} & {\hs} & \StringL{\hs} \\
  &   \StringLRX{\hf} & \hs & \StringLRX{\hs} & \hs & \StringLRX{\hf} \\
  \StringR{\hs} & {\hs} & \StringLR{\hs} & \hs & \StringLRX{\hf} & &  \StringL{\hs} \\
  & 1 \ \ 2 &  & {\color{blue} 3} \ \ 4 & & 5 \ \ 6 \\
}
\]
\end{example}

Let $W^{C} = \bigcup_{n \geq 1} W^{C}_{n}$ denote the convex elements in
one of the classical families of irreducible Weyl groups, type $A$,
$B$ or $D$.  Let $W^{C}(p)$ be the subset of the convex elements in $W$
that is characterized by avoiding a single embedded factor pattern
$p$.

Let $r(p)$ be the rank of the Weyl group containing $p$, and let
$U^{C}(p)$ be the set of all convex elements in $W_{r(p)}$ that
contain $p$ as a factor, i.e.  $U^{C}(p)$ are the convex elements in
the upper order ideal generated by $p$ in the two-sided weak order on
$W_{r(p)}$.  We will show that when $p$ satisfies an additional
hypothesis called the \em ideal \em condition, avoiding $p$ as an
embedded factor is equivalent to avoiding the elements of $U^{C}(p)$ as
1-line patterns.  To carry this out, we will make frequent use of the
linearized, coalesced heap and string diagrams on connected, convex
elements.

We say that a Coxeter embedding on $W$ \em reverses orientation \em if
the labels on the corresponding linear Coxeter graph are reversed
under the embedding.  Otherwise, we say that it \em preserves
orientation\em.  If $w$ contains an embedded factor $p$ under an
orientation preserving Coxeter embedding then we say that $w$
contains $p$ as an \em orientated embedded factor\em.

\begin{proposition}\label{p:embedded.to.one}
If $w \in W^{C}$ contains $p \in W^{C}$ as an oriented embedded
factor, then $w$ contains an element of $U^{C}(p)$ as a 1-line
pattern.
\end{proposition}
\begin{proof}
By Lemma~\ref{l:ef.heap}, we have that $Heap(w)$ contains a shifted
copy of $Heap(p)$ as a saturated set of lattice points.  Furthermore,
we can build $Heap(w)$ from the shifted copy of $Heap(p)$ by
sequentially adding lattice points that are maximal or minimal with
respect to the intermediate heap as follows.  Suppose that the shifted
copy of $Heap(p)$ occupies columns $i, i+1, \dots, k$ in $Heap(w)$.
Then, we begin with the set of strings $S = \{i, i+1, \dots, k+1\}$
that appear in the shifted copy of $Heap(p)$.  These strings initially
correspond to the 1-line pattern $p$, and we show by induction that
$S$ continues to encode a 1-line pattern from $U_p$ as we add minimal
or maximal lattice points to the heap.

Consider the relative order of the strings in $S$ when we add a
maximal lattice point in column $j$ to the heap.  If the new point
crosses a pair of strings that are both in $S$ then the new string
configuration on $S$ corresponds to an element in the upper order
ideal $U^{C}(p)$.  If the new point crosses a pair $x, y$ of strings
such that at most one is contained in $S$ then the string
configuration on $S$ is unchanged.  Similarly, the relative order of
the strings in $S$ corresponds to an element in $U^{C}(p)$ when we add
a minimal lattice point in column $j$.

Hence, at the end of this inductive construction $Heap(w)$ contains
the 1-line pattern encoded by the strings in $S$, and the element
corresponding to this 1-line pattern contains $p$ as a factor.
\end{proof}

Example~\ref{ex:badp} shows that the converse of
Proposition~\ref{p:embedded.to.one} can fail in general.  However, on
the special patterns defined below a converse can be stated.

\begin{definition}\label{d:ideal}
We say that $p\in W^{C}$ is an \em ideal embedded factor pattern \em
if for every $q \in W^{C}_{r(p)+1}$ containing $p$ as a 1-line pattern,
we have that $q$ contains $p$ as an oriented embedded factor.
\end{definition}

\begin{proposition}\label{p:one.to.embedded}
If $p \in W^{C}$ is an ideal embedded factor pattern and $w \in W^{C}$
contains $p$ as a 1-line pattern, then $w$ contains $p$ as an oriented
embedded factor.
\end{proposition}
\begin{proof}
Consider the case that $w \in W_{r(p)+1}$.  By
Definition~\ref{d:ideal}, we have that $w$ contains $p$ as an oriented
embedded factor.  By Lemma~\ref{l:ef.heap}, this implies that
$Heap(w)$ contains $Heap(p)$ as a saturated subset, so we can
highlight an instance of $Heap(p)$ inside $Heap(w)$.

Now by induction, assume the proposition holds for all convex elements
in $\cup_{k=1}^{n} W^{C}_k$ and let $w \in W^{C}_{n+1}$.  Then if $w$
contains $p$ as a 1-line pattern then $w$ contains some $p' \in
W^{C}_{n}$ that also contains $p$ as a 1-line pattern.  By induction,
$Heap(p')$ contains a shifted copy of $Heap(p)$ as a saturated subset
of lattice points, we want to show that $Heap(w)$ must also contain a
copy of $Heap(p)$. 

The string diagram imposed on $Heap(w)$ can be obtained from the
string diagram on $Heap(p')$ by adding one additional string. The
additional string will add extra points to the heap at each crossing.
This string may cut through the copy $C$ of $Heap(p)$, but since $p$
is ideal, the extra points added along with $C$ must also contain a shifted
copy of $Heap(p)$ as a saturated subset of lattice points by
Definition~\ref{d:ideal}.  Therefore by Lemma~\ref{l:ef.heap}, $w$
contains $p$ as an oriented embedded factor.
\end{proof}

Thus combining Proposition~\ref{p:embedded.to.one} and
Proposition~\ref{p:one.to.embedded}, we have shown the following
result.

\begin{theorem}\label{t:one.iff.embedded.finite}
Suppose $P$ is the subset of $W^{C}$ characterized by avoiding a
finite combination of oriented convex embedded factors $F$ and 1-line
patterns $G$.  If each of the elements in $F' = \bigcup_{p \in F}
U^{C}(p)$ is an ideal pattern, then $P$ is characterized by avoiding
the permutations in $F' \union G$ as 1-line patterns.
\end{theorem}

\begin{corollary}
Under the hypotheses of Theorem~\ref{t:one.iff.embedded.finite}, there
is a polynomial time algorithm available to test an element of $W^{C}$
for membership in $P$.
\end{corollary}

\begin{remark}
Note that the existence of a polynomial time algorithm is not evident
from the embedded factor version of the characterization, because a
typical element can have exponentially many reduced expressions by
\cite{stanley}.
\end{remark}

\begin{remark}
Recall that the fully commutative elements in types $A$ and $B$ are
automatically convex.  In particular, the theorem applies to
properties on $[321]$-avoiding permutations that are characterized by
avoiding finitely many embedded factors.  Similarly, any pattern class
that includes the fully commutative basis elements from
Theorem~\ref{t:stem} together with $[\bar{1} 2 \bar{3}]$ is convex in
type $D$.
\end{remark}

\begin{corollary}\label{c:deodhar.poly}
There exist a finite number of patterns of rank less than 9 in each
Weyl group family of types $A,B,D$ that characterize the Deodhar
elements. Therefore, there exists an $O(n^{8})$ test for the Deodhar
condition in all finite Weyl groups of rank $n$.
\end{corollary}

\begin{proof}
It is straightforward to verify that the embedded factor patterns from
Theorem~\ref{t:minimal.patterns} are ideal using a computer.  The
corollary then follows from Theorem~\ref{t:one.iff.embedded.finite}.
\end{proof}

We have verified by computer that the embedded factor patterns
characterizing the Deodhar elements correspond to 75 type $D$
1-line patterns.  

It is not known how to define a pattern system so that a finite
characterization of the Deodhar condition in other Coxeter groups may
be obtained, but Example~\ref{e:infantichain} shows that something
stronger than embedded factor containment is required in general.

\begin{question}\label{conj: rootsubsst}
Can Deodhar elements in other Coxeter groups be characterized by
avoiding a finite number of root subsystem patterns?
\end{question}


\bigskip \section{Toward Enumerating Deodhar Elements}\label{s:open}

Stankova and West \cite{s-w} found a homogeneous linear recurrence
relation with constant coefficients that gives the number of
321-hexagon avoiding permutations.  

\begin{theorem}\cite{s-w}
The number $c_n$ of 321-hexagon-avoiding permutations in $S_n$ satisfies the
recurrence
\[ c_{n+1} = 6 c_{n} - 11c_{n-1} + 9c_{n-2} - 4c_{n-3} - 4c_{n-4} + c_{n-5} \]
for all $n \geq 6$ with initial conditions $c_1 = 1$, $c_2 = 2$, $c_3 = 5$,
$c_4 = 14$, $c_5 = 42$, $c_6 = 132$.
\end{theorem}

Also, Vatter \cite{vatter} has obtained an enumeration scheme for
these elements automatically using the WilfPlus Maple package.  Do
elegant enumerative formulas such as this exist for counting the
Deodhar elements in type $D$?

Figure~\ref{f:enum} shows the number of Deodhar elements for the
finite Weyl groups as a fraction of the fully commutative elements.
The latter were enumerated by \cite{s2}.
\begin{figure}[h]
\begin{center}
\begin{tabular}{|p{0.8in}|p{.2in}|p{0.2in}|p{0.25in}|p{.4in}|p{0.4in}|p{0.5in}|p{0.5in}|} 
	\hline
Type/Rank & 	2 &	 3 &	4 &	5 &	6 &	7 & 8 \\
	\hline
$A$/$B$/$G$/$F$   & 	5 &	 14 &	42 &	132 &	429 &	1426/ 1430 & 4806/ 4862\\
	\hline
$D$               &	   &  &	48 &	167 &	575/ 593 &	1976/ 2144 &	6791/ 7864 \\
	\hline
$E$               &	   &  &	 &	 &	642/ 662 &	2341/ 2670 &	8305/ 10846 \\
\hline
\end{tabular}
\end{center}
\caption{Enumeration of Deodhar elements}
\label{f:enum}
\end{figure}

It would be interesting if one could find an enumerative formula for
the number of Deodhar elements in type $D$.  More generally,
what can be said about enumerating families avoiding a list of
embedded factor patterns? 

It is possible to show that in type $A$, the embedded factor pattern
classes satisfy a Stanley--Wilf bound using a theorem of Tenner.

\begin{theorem}\cite{t2}\label{t:tenner}
If $p \in S_k$ avoids $[2143]$ and $w \in S_n$ contains $p$ as a
1-line permutation pattern, then $w$ contains $p$ as an oriented
embedded factor.
\end{theorem}

We denote the set of elements from $S_n$ that avoid $p$ as an
embedded factor by $S_n(p)$, and the set of elements from $S_n$ that
avoid $[p]$ as a 1-line permutation pattern by $S_n[p]$.

\begin{corollary}
For all permutations $p$, there exists a constant $c = c_p$ such that
$|S_n(p)| \leq c^n$.
\end{corollary}

\begin{proof}
We first show how to construct a permutation $q \in S_k$ that
contains $p$ as a factor, but avoids $[2143]$.
Suppose $p \in S_k$ contains a $[2143]$ instance.  Begin by choosing
the leftmost position $i$ in the 1-line notation for $p$ from the
set of all positions that play the role of $4$ in any $[2143]$
instance of $p$.  Then, we can multiply $p$ on the right by the
adjacent transposition $s_{i-1}$ to move $p_i$ one position to the
left.  Note that $p_{i-1} < p_{i}$ or else we could have chosen
$p_{i-1}$ to play the role of $4$ in any $[2143]$ instance in which
$p_i$ participates, contradicting that $p_i$ was chosen to be
leftmost.

By continuing to move the entry $p_i$ to the left in a reduced
fashion, we can eventually move it past the leftmost entry that plays
the role of $1$ in any $[2143]$ instance in which $p_i$ plays the role
of $4$.  Having removed all of the $[2143]$ instances where the entries that
play the role of $4$ occur in positions $\leq i$, we choose the next
leftmost position that plays the role of $4$ in some $[2143]$
instance and repeat the argument.  The resulting element $q$ contains
$p$ as a factor, and contains no $[2143]$ instances.

Hence, if $w$ contains $q$ as an embedded factor, then it contains $p$
as an embedded factor, so contrapositively $S_n(p) \subset S_n(q)$ and
by Theorem~\ref{t:tenner} we have $S_n(q) \subset S_n[q]$.  We can
apply the Marcus--Tardos Theorem \cite{m-t} to $S_n[q]$ since it is
expressed as a 1-line permutation pattern class, and we obtain the
upper bound.
\end{proof}


\bigskip
\section*{Acknowledgments}

We wish to thank Arkady Berenstein, Hugh Denoncourt, Richard Green,
Matt Kahle, Elizabeth Kelly, William McGovern, Yuval Roichman, John
Stembridge, Bridget Tenner, and Greg Warrington for helpful
conversations while we were working on these results.  We would also
like to thank the anonymous referee for providing many useful
suggestions.



\end{document}